\documentclass[preprint,12pt]{elsarticle}

\usepackage[hidelinks]{hyperref}
\usepackage{graphicx}
\usepackage{subcaption}
\usepackage{amssymb}
\usepackage{amsmath}
\usepackage{multirow}
\usepackage[utf8]{inputenc}
\usepackage{algorithm}
\usepackage[noend]{algpseudocode}
\usepackage{listings}
\usepackage{textcomp}
\usepackage{courier}
\usepackage{tikz}
\usepackage{booktabs}

\usepackage{xcolor}
\usepackage{xargs}
\usepackage[colorinlistoftodos,textsize=footnotesize]{todonotes}

\newcommandx{\unsure}[2][1=]{\todo[linecolor=red,backgroundcolor=red!25,bordercolor=red,#1]{#2}}
\newcommandx{\change}[2][1=]{\todo[linecolor=blue,backgroundcolor=blue!25,bordercolor=blue,#1]{#2}}
\newcommandx{\info}[2][1=]{\todo[linecolor=OliveGreen,backgroundcolor=OliveGreen!25,bordercolor=OliveGreen,#1]{#2}}

\newcommand{\ten}[1]{\ensuremath{\mathbf{#1}}}

\usepackage{lineno}

\journal{}

\begin{document}

\begin{frontmatter}

  \title{An efficient, open source, iterative ISPH scheme}
  \author[IITB]{Abhinav Muta\corref{cor1}}
  \ead{abhinavm@aero.iitb.ac.in}
  \author[IITB]{Prabhu Ramachandran}
  \ead{prabhu@aero.iitb.ac.in}
  \author[IITB]{Pawan Negi}
  \ead{pawan.n@aero.iitb.ac.in }
\address[IITB]{Department of Aerospace Engineering, Indian Institute of
  Technology Bombay, Powai, Mumbai 400076}

\cortext[cor1]{Corresponding author}

\begin{abstract}

  In this paper a simple, robust, and general purpose approach to implement
  the Incompressible Smoothed Particle Hydrodynamics (ISPH) method is
  proposed. This approach is well suited for implementation on CPUs and GPUs.
  The method is matrix-free and uses an iterative formulation to setup and
  solve the pressure-Poisson equation. A novel approach is used to ensure
  homogeneous particle distributions and improved boundary conditions. This
  formulation enables the use of solid wall boundary conditions from the
  weakly-compressible SPH schemes. The method is fast and runs on GPUs without
  the need for complex integration with sparse linear solvers. We show that
  this approach is sufficiently accurate and yet efficient compared to other
  approaches. Several benchmark problems that illustrate the robustness,
  performance, and wide range of applicability of the new scheme are
  demonstrated. An open source implementation is provided and the manuscript
  is fully reproducible.

\end{abstract}


\end{frontmatter}


\section*{Program summary}
\noindent
\textit{Program title}: SISPH

\noindent
\textit{Licensing provisions}: BSD 3-Clause

\noindent
\textit{Programming language}: Python

\begin{sloppypar}
  \noindent
  \textit{External routines/libraries}: PySPH
  (\url{https://github.com/pypr/pysph}), scipy
  (\url{https://pypi.org/project/scipy/}), matplotlib
  (\url{https://pypi.org/project/matplotlib/}), compyle
  (\url{https://pypi.org/project/compyle/}), automan
  (\url{https://pypi.org/project/automan/}).
\end{sloppypar}
\noindent
\textit{Nature of problem}: Incompressible SPH solvers solve a sparse system
of equations to obtain a solution to the pressure-Poisson equation. This is
often implemented by solving a sparse matrix, and can be complex to implement
on GPUs. SPH accuracy deteriorates when particles become disordered. Ensuring
homogeneous particle distributions requires particle shifting algorithms which
requires tuning of parameters. We propose the use of matrix-free algorithm and
an accurate and general purpose way to ensure particle homogeneity.

\noindent
\textit{Solution method}: We use an explicit iterative approach to solve the
pressure-Poisson equations which makes it easy to implement in parallel on
GPUs. This also makes it easy to implement the boundary conditions. We propose
an accurate way to ensure homogeneous particle distribution by using a
background pressure.

\noindent
\textit{Additional comments}: The source code for this repository can be found at
\url{https://gitlab.com/pypr/sisph}.

\section{Introduction}
\label{sec:intro}

The Smoothed Particle Hydrodynamics method was originally proposed to simulate
astrophysical hydrodynamic problems by \citet{lucy77} and
\citet{monaghan-gingold-stars-mnras-77}. It is a mesh-free, Lagrangian,
particle-based method that has since been used to simulate incompressible
fluids. The Weakly-Compressible SPH (WCSPH)~\cite{sph:fsf:monaghan-jcp94}
scheme was proposed to simulate incompressible flows by treating the fluid as
weakly compressible and using a stiff equation of state. This method has been
used to simulate fluid flows with a free-surface. One significant difficulty
with the WCSPH method is the fact that it introduces an artificial sound speed
which is typically around 10 times the maximum speed of the flow. This
introduces severe timestep limitations in the scheme.

The Incompressible SPH (ISPH) schemes have their origin in the work of
\citet{sph:psph:cummins-rudman:jcp:1999} who proposed a projection method
where a pressure-Poisson equation (PPE) is solved to ensure a divergence free
velocity field. In this, the Laplacian of the pressure is related to the
divergence of the velocity field. \citet{isph:shao:lo:awr:2003} proposed a
slightly different variant which relates the Laplacian of the pressure to a
change in density. Later,~\citet{isph:hu-adams:jcp:2007} proposed a method
that combines both of these approaches to enforce incompressibility. The
significant advantage with the class of ISPH schemes is that they do not
involve the sound speed and hence can use timesteps that are an order of
magnitude larger than the WCSPH schemes.

The ISPH schemes have the disadvantage that they require the solution to large
(but sparse) linear systems. This makes implementing them in parallel and on
GPUs fairly difficult. The WCSPH implementations are generally more popular as
they are much easier to implement and parallelize. The recent work of
\citet{chow:isph:cpc:2018} discusses many of the challenges in implementing
traditional ISPH schemes for large scale computing on GPUs.

An explicit version of the ISPH (EISPH) was first proposed
by~\citet{hosseini-2007} which focuses on simulating non-Newtonian flows. This
approach removes the requirement to solve a linear system and instead sets up
an explicit (non-iterative) equation to compute the pressure. A similar
explicit  approach was used by~\citet{rafiee-2009} for solving fluid-structure
interaction problems. Barcarolo~\cite{bar13,barcarolo-2014} validates and
explores the issues with EISPH and compare the scheme with the WCSPH in the
context of internal and free surface flows. \citet{nomeritae-eisph-2016}
assess the performance of EISPH with a comparison with WCSPH and $\delta$-SPH
schemes~\cite{marrone-deltasph:cmame:2011}. \citet{basser-2017} simulate
multi-fluids with porous media using EISPH and compared their results with
experimental data. While~\citet{hosseini-2007}, \citet{rafiee-2009}, and
\citet{barcarolo-2014} solve for the pressure by satisfying a constant density
condition~\cite{isph:shao:lo:awr:2003}, \citet{nomeritae-eisph-2016} and
\citet{basser-2017} use the divergence free
condition~\cite{sph:psph:cummins-rudman:jcp:1999} to obtain the pressure.  It
is important to note that in all these cases, the authors perform a single
step to solve the PPE and do not perform any iterations thus making the method
a truly explicit ISPH method.

In the graphics community, the Implicit ISPH (IISPH)
scheme~\cite{iisph:ihmsen:tvcg-2014} has been developed that provides an
iterative solution to the ISPH formulation. This formulation is tied to the
way in which the pressure forces are approximated between pairs of particles.
The method does not work well with negative pressures. Our own implementation
of this scheme demonstrates some sensitivity to changes in the timestep or
choice of smoothing kernels. The method is however, matrix free, and very fast.

In the present work, we consider the projection-based ISPH scheme of
\citet{sph:psph:cummins-rudman:jcp:1999}. We use an iterative approach to
solve  the PPE as constrasted with the EISPH approach. This makes the
formulation completely matrix-free, and easy to derive and implement. This
approach is not new and is similar to the approach used in the
EISPH~\cite{hosseini-2007,bar13,nomeritae-eisph-2016,basser-2017} with the key
difference being that we iteratively solve the same PPE.  This makes a
significant difference when we simulate problems at higher Reynolds numbers.
We use a successive-over-relaxation (SOR) procedure to accelerate the
convergence of our iterations.  Our implementation is efficient, and works
comfortably on a GPU architecture. This reformulation also allows us to
implement boundary relatively easily as is done in WCSPH schemes and this is
harder to do with the traditional ISPH. We note that despite the existence of
the similar EISPH for many years, the traditional approach of solving the
system of equations using a separate linear solver continues to be
popular~\cite{chow:isph:cpc:2018}.

Both the WCSPH and ISPH schemes suffer from inaccuracies when the particles
become disordered. In flows with significant shear, the particles can become
significantly disordered  leading to poor accuracy and particle clumping in
extreme cases. The ISPH community has developed a few particle shifting
techniques~\cite{acc_stab_xu:jcp:2009,diff_smoothing_sph:lind:jcp:2009,fickian_smoothing_sph:skillen:cmame:2013}
that add a small amount of particle motion to ensure uniform distributions.
These improve the particle distribution at the cost of a negligible amount of
diffusion. It does requires some tuning to handle free surfaces.

For the WCSPH schemes, the Transport Velocity Formulation
(TVF)~\cite{Adami2013} and the Generalized TVF~\cite{zhang_hu_adams17} provide
a slightly different and more accurate way of generating homogeneous particle
distributions. The method relies on the fact that the naive SPH gradient of a
constant pressure field will generate forces that tend to push particles into
a uniform configuration. Hence, the particles are moved with a ``transport
velocity'' which is due to the fluid forces along with this constant pressure
force. The background pressure is subtracted from the main momentum equation
through an artificial stress term. The Generalized
TVF~\cite{zhang_hu_adams17}, extends the TVF to free-surface flows. In this
paper we apply the GTVF to work with ISPH schemes. Since the timesteps are
larger, we use a small modification to move the particles during the
regularization.

Although~\citet{barcarolo-2014} explore internal flows and flow past a
circular cylinder, they do so for low Reynolds numbers where the chances of
particle disorder are less. In fact, in \cite{bar13}, it is found that there
is some particle voiding present for the case of a lid driven cavity at a
Reynolds number of 3200. Moreover, in EISPH, the PPE is solved using a single
step. This approach does not work well for high Reynolds number cases and we
demonstrate these issues in the current work. Previous works employing the
EISPH do not explore the performance on GPUs where the scheme has a high
potential for rapid, accurate simulation. Furthermore, they do not simulate
the Taylor-Green problem which typically fails to work unless a careful
implementation of particle shifting is included. Our implementation, on the
other hand works well for this problem and also demonstrates good performance
on GPUs.

In order to handle solid boundaries accurately, we modify the method proposed
by~\citet{Adami2012} to work with the ISPH schemes. The original method does
not work too well with the ISPH due to the larger timesteps that the method
allows. We provide a simple way to compute the normals on a solid body with
just the point information, and use this to improve the solid boundary
condition. The resulting method works effectively for a variety of test cases.

In summary, the proposed method is easy to implement, accurate, matrix-free,
fast, and works on GPUs. We apply the transport velocity formulation to ensure
that the particle distribution is uniform and suggest an improved solid wall
boundary condition.  We demonstrate the method with several standard
benchmarks. These include internal flows, external free-surface problems,
wind-tunnel type problems requiring an inlet and an outlet and also
demonstrate good performance on a GPU. We demonstrate the performance of the
new scheme and compare this with the performance of the traditional
matrix-based ISPH. We also demonstrate the performance achieved on different
GPUs as well as multi-core CPUs.

Our entire implementation is open source. We employ the open source
PySPH~\cite{PR:pysph:scipy16,pysph} framework for the implementation of the
scheme. In the interest of reproducible research, our entire manuscript is
reproducible. Every figure presented in the results section of this manuscript
is automated~\cite{pr:automan:2018} and the code for the computations is made
available at \url{https://gitlab.com/pypr/sisph}.

\section{The ISPH method}
\label{sec:sph}

The basic formulation is that of~\citet{sph:psph:cummins-rudman:jcp:1999}. For
an incompressible fluid, the continuity equation is,
\begin{equation}
  \label{eq:continuity}
  \nabla \cdot \ten{u} = 0
\end{equation}
where $\ten{u}$ is the velocity field, and the momentum equation is given by,
\begin{equation}
  \label{eq:mom}
  \frac{d \ten{u}}{dt} = - \frac{1}{\rho} \nabla p + \nu \nabla \cdot (\nabla
  \ten{u}) + \ten{f}
\end{equation}
where $\frac{d\ten{u}}{dt}$ is the material derivative of the velocity field,
$\rho$ is density, $p$ is pressure, $\nu$ is kinematic viscosity and $\ten{f}$
is the external force.

Pressure is obtained by the projection method, which uses Hodge decomposition
to project any velocity field, $\ten{u}^{*}$ into a divergence-free component, and a
curl-free component,
\begin{equation}
  \ten{u}^* = \ten{u} + \frac{\Delta t}{\rho} \nabla p
\label{eq:isph:proj}
\end{equation}
where $\ten{u}$ is the divergence free velocity field, $\nabla p$ is the
curl-free component. the pressure is then obtained by taking divergence of
equation~\eqref{eq:isph:proj}, which reads,
\begin{equation}
  \nabla \cdot (\frac{1}{\rho}\nabla p) = \frac{\nabla \cdot \ten{u}^*}{\Delta t}
\label{eq:isph:ppe}
\end{equation}
where $\ten{u}^*$ is the intermediate velocity obtained by integrating the
momentum equation without considering the pressure gradient
terms~\cite{sph:psph:cummins-rudman:jcp:1999}. After solving for the pressure, $p$
(curl-free component), the divergence free velocity $\ten{u}$ is obtained
by subtracting the gradient of pressure from the intermediate velocity
$\ten{u}^*$.

The ISPH like most SPH methods suffers from particle disorder. We consider the
Generalized Transport Velocity Formulation (GTVF) by~\citet{zhang_hu_adams17}
to remedy the particle disorder. GTVF extends the Transport Velocity
Formulation (TVF) by~\citet{Adami2013} to free surface flows and solid
dynamics by using variable background pressure instead of a constant global
background pressure. In the GTVF formulation, positions and velocity are
advected using the transport velocity, which introduces an additional
artificial stress term to the momentum equation. The transport velocity is
regulated by variable background pressure, which keeps the particles in a
uniform configuration. The momentum equation for GTVF is,
\begin{equation}
  \label{eq:mom:gtvf}
  \frac{\tilde{d} \ten{u}}{dt} = - \frac{1}{\rho} \nabla p
+ \nu \nabla \cdot (\nabla \ten{u}) + \frac{1}{\rho} \nabla \cdot \ten{A} +
\ten{f}
\end{equation}
where $\frac{\tilde{d}\ten{u}}{dt} = \frac{\partial
\ten{u}}{\partial t} + \tilde{\ten{u}} \cdot \nabla \ten{u}$ is the material
derivative of velocity field which is advected by the
transport velocity, $\tilde{\ten{u}}$, and $ \ten{A} = \rho \ten{u} \otimes
(\tilde{\ten{u}} - \ten{u})$ is the artificial stress term.

\subsection{SPH discretization}

We apply SPH discretization to the scheme discussed above. In all the
simulations we do not keep the density constant, density is computed by the
summation density,
\begin{equation}
  \label{eq:summation_density}
  \rho_i = \sum_{j \in N(i)} m_j W_{ij}
\end{equation}
where the subscript $i$ denotes the $i^{\text{th}}$ particle, $W_{ij} =
W(|\ten{r}_{ij}|, h_{ij})$ is the SPH smoothing kernel, where $\ten{r}_{ij} =
\ten{r}_i - \ten{r}_j$, $h_{ij} = (h_i + h_j)/2$, $\ten{r}$ is the position of
the particle and $h$ the smoothing length of the kernel. The summation is
carried over $N(i)$ neighbors, of the $i^{\text{th}}$ particle. We use the
density in order to provide a better estimate of the particle volume given by
$m/\rho$.

The SPH discretization of the momentum equation~\eqref{eq:mom:gtvf} is given
by,
\begin{equation}
  \label{eq:mom:sph}
  \frac{\tilde{d} \ten{u}_i}{dt} =  \ten{f}_{\text{p}} + \ten{f}_{\text{visc}} +
  \ten{f}_{\text{b}} + \ten{f}_{\text{avisc}} + \ten{f}_{\text{astress}}
\end{equation}
where $ \ten{f}_{\text{p}}$ is the force due to pressure gradient,
$\ten{f}_{\text{visc}}$ is the force due to viscous forces,
$\ten{f}_{\text{b}}$ is the body force, $\ten{f}_{\text{avisc}}$ is the force
due to artificial viscosity and $\ten{f}_{\text{astress}}$ is the force due
to artificial stress which is a consequence of the GTVF formulation.

In the literature two ways of computing pressure gradient are encountered. One
is a symmetric, momentum-preserving form~\cite{isph:hu-adams:jcp:2007,
  sph:psph:cummins-rudman:jcp:1999} and the other is asymmetric and does not
preserve momentum
\cite{acc_stab_xu:jcp:2009, fickian_smoothing_sph:skillen:cmame:2013,
  diff_smoothing_sph:lind:jcp:2009}. In the present work the symmetric form is
labelled as ``symm'' and the asymmetric form as ``asymm''. We have used
both the forms in this work and report the differences between them in our
results. The SPH discretization of the symmetric form~\cite{monaghan92} is,
\begin{equation}
  \label{eq:mom:sph:press:symm}
  \ten{f}_{\text{p, symm}} = - \sum_{j \in N(i)} m_j
  \left( \frac{p_i}{\rho_i^2} + \frac{p_j}{\rho_j^2} \right) \nabla W_{ij}
\end{equation}
and the asymmetric form~\cite{monaghan92} is,
\begin{equation}
  \label{eq:mom:sph:press:pij}
  \ten{f}_{\text{p, asymm}} = - \sum_{j \in N(i)}\frac{m_j}{\rho_i\rho_j}(p_j
    - p_i) \nabla W_{ij}
\end{equation}
where $\nabla W_{ij}$ is the gradient of the smoothing kernel. Viscous forces
are discretized as~\cite{sph:psph:cummins-rudman:jcp:1999},
\begin{equation}
  \label{eq:mom:sph:visc}
  \ten{f}_{\text{visc}} = \sum_{j \in N(i)} m_j \frac{4 \nu}{(\rho_i + \rho_j)}
  \frac{\ten{r}_{ij} \cdot \nabla W_{ij}}{(|\ten{r}_{ij}|^{2} + \eta h_{ij}^{2})}
  \ten{u}_{ij}
\end{equation}
where $\ten{u}_{ij} = \ten{u}_i - \ten{u}_j$, $\nu$ is the kinematic
viscosity, $\eta = 0.01$. Artificial viscosity~\cite{monaghan-review:2005} is
added to the momentum equation wherever necessary given by,
\begin{equation}
  \label{eq:mom:sph:av}
  \ten{f}_{\text{avisc}} = \sum_{j \in N(i)} \Pi_{ij} \nabla W_{ij}
\end{equation}
where, $\Pi_{ij}$ is computed by,
\begin{align}
  \label{eq:mom:sph:av:params}
  \Pi_{ij} =
  \begin{cases}
    \frac{-\alpha h_{ij} \bar{c}_{ij} \phi_{ij}}{\bar{\rho}_{ij}},
    & \ten{u}_{ij}\cdot \ten{r}_{ij} < 0 \\
    0, & \ten{u}_{ij}\cdot \ten{r}_{ij} \ge 0 .\\
  \end{cases}
\end{align}
The artificial stress due to the GTVF formulation~\cite{zhang_hu_adams17} is
evaluated using,
\begin{equation}
  \label{eq:mom:sph:gtvf:stress}
  \ten{f}_{\text{astress}} = \sum_{j \in N(i)} \left(\frac{\ten{A}_i}{\rho_i^2} +
  \frac{\ten{A}_j}{\rho_j^2} \right)
\end{equation}
where, ${A}_{i} = \rho_{i} \ten{u}_{i} \otimes(\tilde{\ten{u}}_{i} -
\ten{u}_{i})$. The transport velocity, $\tilde{\ten{u}}$ due to GTVF~\cite{zhang_hu_adams17} is
updated by adding an additional background pressure to the simulation,
\begin{equation}
  \label{eq:int:vhatnp1}
  \tilde{\ten{u}}^{n+1}_i = \ten{u}_i + \Delta t \left(\frac{\tilde{d} \ten{u}_i}{dt}
  + \ten{f}_{i, \text{GTVF}}\right)
\end{equation}
where, $\frac{\tilde{d} \ten{u}_i}{dt}$ is obtained from
equation~\eqref{eq:mom:sph} and
\begin{equation}
  \label{eq:mom:sph:gtvf}
  \ten{f}_{i, \text{GTVF}} = - p^0_i \sum_{j
    \in N(i)} \frac{m_j}{\rho_j^2} \nabla W(\ten{r}_{ij}, \tilde{h}_{ij}).
\end{equation}
On numerical investigation, we suggest $\tilde{h} = 0.5 h_{ij}$ for external
and free surface flows, and $\tilde{h}_{ij} = h_{ij}$ for internal flows,
$p^{0}_i = \min(10|p_i|, p_{\text{ref}})$ for external flows, and $p^0_i =
p_{\text{ref}}$ is kept constant for internal flows. The choice of reference
pressure, $p_{\text{ref}}$, is typically taken to be $\rho c^2$ for external
flows and $2\max(p_i)$, i.e., twice the maximum pressure in the entire flow
after the first iteration for internal flows. Here, we assume that $c = 10
|\ten{u}|_{\text{max}}$ and this is simply a reference pressure.



The pressure-Poisson equation that is solved is,
\begin{equation}
  \label{eq:ppe}
  \nabla \cdot {\left( \frac{1}{\rho} \nabla p \right)}_i =
  \frac{\nabla \cdot \ten{u^*}_i}{\Delta t}
\end{equation}
where $\Delta t$ is the timestep. The approximate projection form for the PPE
operator in equation \eqref{eq:ppe} given
by~\citet{sph:psph:cummins-rudman:jcp:1999} is used, resulting in the
discretized form as,
\begin{align}
  \label{eq:sph-ppe}
  \sum_j \left( \frac{4 m_j}{\rho_i (\rho_i + \rho_j)} \right)
  \frac{p_{ij} \ten{r}_{ij} \cdot \nabla_i W_{ij}}{|r^2_{ij}| + \eta h^2_{ij}}
  &= \sum_j -\frac{m_j}{\rho_j \Delta t} \ten{u^*}_{ij} \cdot \nabla_i W_{ij},
\end{align}
where $p_{ij} = p_i - p_j$.

\subsection{Time Integration}

Due to its simplicity we use the time integration as given
by~\citet{sph:psph:cummins-rudman:jcp:1999}. However, the higher order method
proposed in~\citet{nair2015} could also be used. Time integration is performed
by first calculating the intermediate position using the transport velocity at
the current time,
\begin{equation}
  \label{eq:int:rnp1}
  \ten{r}^{*}_i = \ten{r}^{n}_i + \Delta t \tilde{\ten{u}}^{n}_i
\end{equation}
The intermediate velocity is obtained by integrating the momentum
equation neglecting the force due to the pressure gradient,
\begin{equation}
  \ten{u}^{*}_i = \ten{u}^{n}_i + \Delta t ( \ten{f}^{n}_{i, \text{visc}} +
    \ten{f}^{n}_{i, \text{avisc}} + \ten{f}^{n}_{i, \text{astress}} +
    \ten{f}^{n}_{i, \text{body}}).
\label{eq:int:vint}
\end{equation}
The pressure-Poisson equation given by \eqref{eq:sph-ppe} is solved to obtain
the pressure. The divergence free velocity is then found by correcting the
intermediate velocity as,
\begin{equation}
  \label{eq:int:vnp1}
  \ten{u}^{n+1}_i = \ten{u}^{*}_i + \Delta t \ten{f}^{n}_{i, \text{p}}.
\end{equation}
We use modified GTVF to redistribute particles. The GTVF formulation is first
used in the context of weakly compressible flows, where the time step is much
lower than ISPH schemes. The timesteps are large in the ISPH and applying the
GTVF correction in one timestep introduces stability issues. We modify the
GTVF formulation by applying the GTVF force ($\ten{f}_{\text{GTVF}}$,
equation~\eqref{eq:mom:sph:gtvf}) to shift the particles iteratively in one
time step. We then compute the transport velocity using the initial and final
positions after the iteration.

The GTVF acceleration is a function of positions $\ten{r}$, background
pressure $p_{\text{ref}}$ and smoothing length $h_{ij}$. In all our
simulations $h_{ij}$ is fixed and for a given time the background pressure
remains unchanged. We use $K$ sub-timesteps to integrate the GTVF
accelerations without recomputing any neighbors such that $\Delta \tau =
\Delta t/K$. If we denote $\ten{r}^k, \ten{\tilde{u}}^k$ as the position and GTVF
velocity of a particle at the $k$'th sub-iterate, we first start with
$\ten{r}^0$ set to the current particle position, and set $\ten{\tilde{u}}^0=0$.
\begin{equation}
  \label{eq:int:gtvf:pos}
  \ten{r}^{k} = \ten{r}^{(k-1)} + \Delta \tau
  \tilde{\ten{u}}^{k-1} + \frac{1}{2} (\Delta \tau)^2 \ten{f}({\ten{r}^{k-1}})_{\text{GTVF}}
\end{equation}
\begin{equation}
  \label{eq:int:gtvf:vel:tmp}
  \tilde{\ten{u}}^{k}  =  \tilde{\ten{u}}^{(k-1)}
  + \Delta \tau \ten{f}({\ten{r}^{(k-1)}})_{\text{GTVF}}
\end{equation}
After the GTVF iterative step is done, we compute the transport velocity by
looking at the inital and final positions
\begin{equation}
  \label{eq:int:vhatnp2}
  \tilde{\ten{u}}^{n+1}_i = \ten{u}^{n+1}_i + \frac{\ten{r}^{K} - \ten{r}^0}{\Delta t}
\end{equation}
particles are then advected to the next time step using transport velocity of
current and previous time steps given by,
\begin{equation}
  \label{eq:int:rnp2}
  \ten{r}^{n+1}_i = \ten{r}^{n}_i + \Delta t
  \left(\frac{\tilde{\ten{u}}^{n+1}_i + \tilde{\ten{u}}^{n}_i}{2} \right)
\end{equation}

\section{Iterative formulation}
\label{sec:isph}

\newcommand{\OD}{\ensuremath{\text{OD}}}

Traditionally, equation~\eqref{eq:sph-ppe} is solved by setting up a sparse
matrix corresponding to the coefficients of $p_i$ and $p_j$ in the equation
and computing the right-hand side. This system is then solved using any fast,
sparse matrix solver. In the
EISPH~\cite{hosseini-2007,bar13,nomeritae-eisph-2016,basser-2017} a simple
explicit formula is provided to update the pressure.

We use this approach to write this system in terms of a Jacobi iteration we
introduce a time index, $k$, and write $p^k$ as the pressure at the $k$th
iteration. We rewrite equation~\eqref{eq:sph-ppe} in the form,

\begin{equation}
  \label{eq:iterative-p}
  D_{ii} p_i^{k+1} = \text{RHS}_i - \sum_j \OD_{ij}p_j^{k},
\end{equation}
where $\OD_{ij}$ represents the off-diagonal coefficients which are multiplied
by pressure at the earlier iteration $p^k$ and $\text{RHS}_i$ is the
right-hand-side of equation~\eqref{eq:sph-ppe}. We can easily write the
following expressions for the coefficients in equation~\eqref{eq:iterative-p}
as,
\begin{align}
  \label{eq:coeff}
  D_{ii} &= \sum_j \left( \frac{4 m_j}{\rho_i (\rho_i + \rho_j)} \right)
           \frac{ \ten{r}_{ij} \cdot \nabla_i W_{ij}}{|r^2_{ij}| + \eta^2} \\
  \OD_{ij} &= \left( \frac{-4  m_j}{\rho_i (\rho_i + \rho_j)} \right)
             \frac{ \ten{r}_{ij} \cdot \nabla_i W_{ij}}{|r^2_{ij}| + \eta^2}
             \label{eq:coeff1}
\end{align}

To begin with, we set $p^{k=0}(t + \Delta t) = p(t)$, that is our first
iteration starts with the existing value of pressure of the particle. We find
the pressure for the next iteration using Successive-Over-Relaxation (SOR) as,
\begin{align}
  \label{eq:p-sor}
  p_i^{k+1} = \omega \frac{(\text{RHS}_i - \sum_j \OD_{ij}p_j^k)}{D_{ii}} +
  (1 - \omega) p_i^{k}
\end{align}
The default value of $\omega=0.5$. We note that the entire term $\sum_j
\OD_{ij}p_j^{k}$, needs to be accumulated into a single term and hence only
requires a single number per particle for storage. Thus, per particle we
require storing three terms, viz.~$D_{ii}, \sum_j \OD_{ij}p_j^{k}$, and
$\text{RHS}_i$.
We note that if a particle has no neighbors, then $D_{ii}$ can be zero, in
which case we set the pressure of the particle to zero. Similarly, to handle
free surfaces we compute the summation density of the particles and if the
ratio $\rho_i/\rho_0 < 0.8$, where $\rho_0$ is the reference or rest density,
we treat the particle as a free-surface particle and set its pressure to zero.
Note that during these pressure iterations we do not move the particles.
We continue to iterate until the relative change in the pressure is less than
some user-specified amount,
\begin{align}
  \label{eq:convergence}
  \sum_i \frac{|p_i^{k+1} - p_i^k|}{\max(\Lambda, \sum_i |p_i^{k+1}|)} < \epsilon,
\end{align}
where $\Lambda = \max(RHS_i/D_{ii})$. When the mean pressure in the flow is
less than one, we use the largest value of $RHS_i/D_{ii}$ as a measure of the
mean pressure and this conveniently avoids zeros in the denominator while also
using a relative scale for the error. The minimum number of iterations
performed is set to two.

We vary the tolerance parameter, $\epsilon$, in our numerical studies. We
consider the suitable choice of this tolerance parameter next.

\subsection{Choice of the tolerance parameter}
\label{sec:eps-analysis}

While iterating for convergence, we note that using very small values of
$\epsilon$ is not necessary. We first note that the discretized terms in
equation~\eqref{eq:coeff} involve errors of at best $O(h^2)$ so there is no
major advantage to solving the linear system to high accuracy.

Furthermore, for many timesteps, the pressure difference at a point between
two timesteps is a small quantity and since we start with an accurate
pressure, we expect that obtaining the new pressure will take very few
iterations. By setting a tolerance of say $0.001$, we ensure that if there are
any changes to the pressure during the Jacobi iterations that are less than
0.1\% that we stop iterations. In practice, we find that most often only a few
iterations (typically 2--10) are necessary. It is important to note that
Chorin~\cite{chorin1968}, in his original fractional step work, makes a
similar suggestion. In addition as mentioned in the introduction, the
EISPH~\cite{hosseini-2007,bar13,nomeritae-eisph-2016,basser-2017} only uses a
single iteration and obtains acceptable results for a wide variety of
problems. However, in the present work, we find that a single iteration is
insufficient especially in the case of high-Reynolds number problems and we
explore this in section~\ref{sec:tg-eisph}. As such, we recommend choosing the
$\epsilon$ to be of $O(h)$ or $O(\Delta t)$.

In Table~\ref{table:iters:avg}, we show the average number of Jacobi
iterations for different problems and tolerances.  These results seem to show
that using Jacobi iterations is very effective.
\begin{table}[!h]
\centering
\begin{tabular}{lrr}
\toprule
      Problem & Tolerance &  Jacobi (avg no of. iterations) \\
\midrule
       Cavity &     0.001 &                             3.0 \\
 Dam-break 2D &     0.001 &                             5.0 \\
       Cavity &      0.01 &                             2.0 \\
 Dam-break 2D &      0.01 &                             2.0 \\
\bottomrule
\end{tabular}

\caption{Average number of iteration taken for Jacobi for various tolerances
  and problems.}
\label{table:iters:avg}
\end{table}
Although we do not discuss this in detail, our repository contains an
implementation of a BiCGStab~\cite{Vor92} iterative algorithm to solve
the PPE. A BiCGStab is needed since when we apply the boundary conditions, the
resulting matrix may become non-symmetric. We find that implementing the
BiCGStab correctly is much more complex as compared to the Jacobi iterations.
Further, even though it converges faster, we find that the iterations do not
seem to produce significantly better results. We find that in practice, it is
much faster to use the Jacobi approach since most of the time we only need a
few iterations to converge. This suggests that using more accurate iterative
sparse linear solvers may be unnecessary.

We therefore recommend the use of the Jacobi iterations approach for its
efficiency and simplicity.

\subsection{The Algorithm}

\begin{algorithm}[t]
\caption{Simple Iterative ISPH algorithm}\label{alg:sisph}
\begin{algorithmic}[1]
\While{$t < t_{\text{final}}$}
\State{Compute $\rho_i$ using equation~\eqref{eq:summation_density}.}
\State{Predict position $\ten{r}^{*}_i \leftarrow \ten{r}^{n}_i
  + \Delta t \tilde{\ten{u}}^i$}
\State{Compute $\ten{f}_{i, \text{visc}}$\,, $\ten{f}_{i, \text{avisc}}$\,,
  $\ten{f}_{i, \text{b}}$\,, and $\ten{f}_{i, \text{astress}}$.}
\State{Predict velocity $\ten{u}_i^{*} \leftarrow \ten{u}_i^{n} + \Delta t
  (\ten{f}_{i, \text{visc}} + \ten{f}_{i, \text{b}} + \ten{f}_{i,
    \text{avisc}} + \ten{f}_{i, \text{astress}})$}
\State{Compute $\nabla \cdot \ten{u}^{*}_i$ from r.h.s.\ of
  equation~\eqref{eq:sph-ppe}.}
\While{$\sum_i \frac{|p_i^{k+1} - p_i^k|}{\max(1, \sum_i |p_i^{k+1}|)} < \epsilon$}
\State{Compute $D_{ii}$ and $OD_{ij}$\,.}
\State{Compute $p^{k+1}_i$ using equation~\eqref{eq:p-sor}.}
\EndWhile{}
\State{Compute $\ten{f}_{i, \text{p}}$ either from
equation~\eqref{eq:mom:sph:press:pij} or~\eqref{eq:mom:sph:press:symm}.}
\State{Correct velocity $\ten{u}_i^{n+1} \leftarrow \ten{u}_i^{*} + \Delta
  t \ \ten{f}_{i, \text{p}}$.}
\State{Set $k \leftarrow 0$, $\ten{r}^{k=0} \leftarrow \ten{r}^{*}$,
  $\tilde{\ten{u}}^{0} \leftarrow 0$}
\While{$k < K$}
\State{Compute $\ten{f}_{i, \text{GTVF}}(\ten{r}^{k})$ using
  equation~\eqref{eq:mom:sph:gtvf}.}
\State{Compute $\ten{r}^{k+1}$ using~\eqref{eq:int:gtvf:pos}}
\State{Compute $\tilde{\ten{u}}^{k+1}$ using~\eqref{eq:int:gtvf:vel:tmp}}
\EndWhile{}
\State{Update GTVF velocity to $\tilde{\ten{u}}_i^{n+1} \leftarrow
  \ten{u}_i^{n+1} + \frac{\ten{r}^{K} - \ten{r}^0}{\Delta t}$.}
\State{Update the positions $ \ten{r}^{n+1}_i = \ten{r}^{n}_i + 0.5 \Delta t
  (\tilde{\ten{u}}^{n+1}_i + \tilde{\ten{u}}^{n}_i)$}
\EndWhile{}
\end{algorithmic}
\end{algorithm}
For clarity, the Algorithm~(\ref{alg:sisph}) shows the steps involved in
the proposed scheme.
This simple formulation works very well in practice as borne out by our
numerous benchmarks. The formulation allows us to satisfy boundary conditions
that are traditionally not easy to do with a matrix-based formulation.

The time steps are chosen based on the following criteria, the CFL criterion
is,
\begin{equation}
  \Delta t_{\text{cfl}} = 0.25 \frac{h_{ij}}{|\ten{U}|}
\end{equation}
where $|\ten{U}|$ is the magnitude of maximum velocity in the simulation, the
viscous criterion is,
\begin{equation}
  \Delta t_{\text{visc}} = 0.25 \frac{h_{ij}^2}{\nu}
\end{equation}
the condition due to external force (if applicable) is,
\begin{equation}
  \Delta t_{\text{force}} = 0.25 \sqrt{\frac{h_{ij}}{g}}
\end{equation}
time steps are then chosen to be the minimum of above conditions,
\begin{equation}
  \Delta t = \min{(\Delta t_{\text{cfl}}, \Delta t_{\text{visc}}, \Delta
    t_{\text{force}})}
\end{equation}

Adaptive time steps can be used in the following way. At each time step the
maximum velocity in the simulation is calculated, and used in the CFL
condition,
\begin{equation}
  \Delta t^{n+1}_{\text{cfl}} = 0.25 \frac{h_{ij}}{\max{|{\ten{u}^{n}_i}|}}
\end{equation}
where $\max{|\ten{u}^{n}_i|}$ is the maximum magnitude of velocity in the
domain at the current time step. Similarly, the time step restriction due to
force is calculated as,
\begin{equation}
  \Delta t^{n+1}_{\text{force}} = 0.25\sqrt{\frac{h_{ij}}{\max{|\ten{f}^{n}_i|}}}
\end{equation}
where $\ten{f}^{n}_i = \ten{f}^{n}_{\text{p}} + \ten{f}^{n}_{\text{visc}} +
\ten{f}^{n}_{\text{body}} + \ten{f}^{n}_{\text{avisc}} +
\ten{f}^{n}_{\text{astress}}$. The restriction due to viscosity has no effect
as viscosity is held constant in all our simulations. Hence the time step for
the next iteration is given by,
\begin{equation}
  \Delta t^{n+1} = \min{(\Delta t^{n+1}_{\text{cfl}}, \Delta t^{n+1}_{\text{force}})}
\end{equation}

\subsection{Implementation details}
\label{sec:pysph}
Here we outline the procedure in solving the iterative formulation using the
PySPH framework \cite{PR:pysph:scipy16,pysph}. PySPH is a Python framework
which implements various SPH formulations. The user code is written in Python
from which high performance serial (OpenMP) or parallel (OpenCL or CUDA C)
code is automatically generated. PySPH also supports multi-CPU execution using MPI.

Listing~\ref{lst:equations} shows the Python code used to define the
inter-particle interactions for a Taylor-Green problem. This is done in two
separate stages. Each stage is defined as a list of equations. Each equation
has a \texttt{dest} and \texttt{sources} keyword argument. The \texttt{dest}
refers to the particles on which the equation is to be solved and
\texttt{sources} is the list particles that influence the \texttt{dest}
particles. \texttt{stage1} is a list of equations that are to be solved before
the first integration step similarly \texttt{stage2} are solved before the second
integration step. A \texttt{Group} is a PySPH construct which updates all the
particles with a given set of equations. A \texttt{Group} is solved using the
current particle properties. For example, in the \texttt{stage1}, the
summation density is evaluated and
updated density is found for all the particles before moving on to the next
group of equations. In the Listing~\ref{lst:equations}, we see that the first
stage performs the steps 2-4 in Algorithm~\ref{alg:sisph}. The second stage
performs the remaining computations. Step 5 in the algorithm is performed by
the \texttt{VelocityDivergence} equation, next the equations of steps 7 and 8
are performed in an iterated group. The iterated group iterates until the
convergence criterion is satisfied or if predefined maximum iterations are completed (in
practice this limit is never reached). Finally, the last group performs steps
9 and 10. A suitable integrator updates the velocities and positions. More
details on PySPH are available in~\cite{PR:pysph:scipy16}.

An example equation is shown in Listing~\ref{lst:coeff_equation}. This shows
how the equation for the pressure coefficients,
i.e.~equations~\eqref{eq:coeff} and \eqref{eq:coeff1}, are written. The
\texttt{initialize} method of the class initializes the variables to zero,
\texttt{d\_idx} is the index of \texttt{dest} particles. The \texttt{loop} method is called
for every pairwise interaction with the neighbors, where \texttt{s\_idx} is
the index of the source (neighbor) particles. The source particle properties are
prefixed with \texttt{s\_} and the destination with \texttt{d\_}, for example
\texttt{s\_pk} is an array of the pressure of the source particles at the
$k$\textsuperscript{th} iteration, $p_k$, and \texttt{d\_diag} is the diagonal
term, $D_{ii}$ of the destination particle. Various quantities like
\texttt{DWIJ} which is the gradient vector for the current particle, and
\texttt{XIJ} is the distance between destination particle and source particle
are available to each method. The listings demonstrate the ease with which
the new scheme can be implemented in the PySPH framework.

\lstset{caption={Algorithm in PySPH for a Taylor-Green simulation, showing two
    stages which are executed before the integration steps are performed.
    }} \lstset{basicstyle=\footnotesize\ttfamily}
\begin{lstlisting}[label={lst:equations},frame=lines,language=Python,upquote=True]
stage1 = [
    Group(equations=[
        SummationDensity(dest='fluid', sources=['fluid'])
    ]),
    Group(equations=[
        LaminarViscosity(dest='fluid', sources=['fluid'],
                            nu=0.01,
        ),
        MomentumEquationArtificialStress(
            dest='fluid', sources=['fluid'], dim=2
        )
    ])]

stage2 = [
    Group(equations=[
        VelocityDivergence(dest='fluid', sources=['fluid'])
    ]),
    Group(equations=[
        Group(equations=[
            PressureCoeffs(dest='fluid', sources=['fluid']),
            PPESolve(dest='fluid', sources=['fluid'],
                     omega=0.5, tolerance=0.01)
        ],
        iterate=True, max_iterations=1000, min_iterations=2)
    ]),
    Group(equations=[
        MomentumEquationPressureGradient(
            dest='fluid', sources=['fluid']
        ),
        GTVFAcceleration(dest='fluid', sources=['fluid'],
                         pref=100.0)
    ])]
\end{lstlisting}
\lstset{caption={An example of equation in PySPH framework showing the
        implementation of diagonal and off-diagonal terms in the pressure
        solver.}}
\lstset{basicstyle=\footnotesize\ttfamily}
\begin{lstlisting}[label={lst:coeff_equation},language=Python,frame=lines]
class PressureCoeffs(Equation):
    def initialize(self, d_idx, d_diag, d_odiag):
        d_diag[d_idx] = 0.0
        d_odiag[d_idx] = 0.0

    def loop(self, d_idx, s_idx, s_m, d_rho, s_rho, d_diag,
             d_odiag, s_pk, XIJ, DWIJ, R2IJ, EPS):
        rhoij = (s_rho[s_idx] + d_rho[d_idx])
        rhoij2_1 = 1.0/(d_rho[d_idx]*rhoij)

        xdotdwij = (XIJ[0]*DWIJ[0] + XIJ[1]*DWIJ[1]
                    + XIJ[2]*DWIJ[2])

        fac = 4.0*s_m[s_idx]*rhoij2_1 * xdotdwij / (R2IJ + EPS)

        d_diag[d_idx] += fac
        d_odiag[d_idx] += -fac * s_pk[s_idx]
\end{lstlisting}

\subsection{Solid wall boundary conditions}
\label{sec:bc}

In order to satisfy solid wall boundary conditions, the method proposed by
\citet{Adami2012} is slightly modified to work with the ISPH scheme. This
method uses 3--4 layers of dummy (or ghost) solid particles and then performs
a Shepard interpolation of the pressure and also accounts for any fluid
acceleration to set the solid wall pressure. As discussed in the original
paper the pressure of the dummy wall particles are set using,
\begin{align}
  \label{eq:p-shepard}
  p_w = \frac{\sum_f p_f W_{wf} + (\ten{g} - \ten{a}_w) \cdot
  \sum_f \rho_f \ten{r}_{wf} W_{wf}}{\sum_f W_{wf}},
\end{align}
where the subscript $w$ denotes a wall particle and $f$ a fluid particle,
$\ten{a}_w$ denotes the acceleration of the wall and $\ten{g}$ the
acceleration due to gravity if relevant. In the case of free-surface flows, we
ensure that the wall does not have any negative pressures as this often causes
particles to stick to and sometimes penetrate the body. Hence if the pressure
on the wall is negative, we set it to zero.

In our iterative scheme to solve for the pressure, during each iteration, we
perform a Shepard interpolation of the fluid pressure using
equation~\eqref{eq:p-shepard} and use those values in evaluating the diagonal
and off-diagonal terms.  Thus the solid particles are treated just as fluids
but with their pressure evaluated based on the extrapolated fluid pressure
from the previous iteration.
The dummy particle velocities are normally set by first calculating a
Shepard-extrapolated fluid velocity on the ghost particles using,
\begin{align}
  \label{eq:vf}
  \ten{\tilde{u}}_i =
  \frac{\sum_j \ten{u}_j W_{ij}(h_{ij}/2)}{\sum_j W_{ij}(h_{ij}/2)},
\end{align}
where $h_{ij} = (h_i + h_j)/2$. Note that we use a smaller kernel radius than
the original scheme to ensure that the closest fluid particles have the most
influence on the wall velocity. The velocity of the dummy particles is then
set by,
\begin{equation}
  \label{eq:vg-adami}
  \ten{u}_w = 2\ten{u}_i - \ten{\tilde{u}}_i,
\end{equation}
where $\ten{u}_i$ is the physical velocity of the wall itself. If a no-slip
boundary condition is satisfied, the above wall velocity is used for the wall
particles when calculating the viscous forces. For the ISPH method, due to the
larger allowed timesteps, this basic approach generally leads to some amount
of leakage. We modify this boundary condition and ensure that the resulting
ghost velocity $\ten{u}_w$ has no component into the solid itself. This is
done by using the normal of the solids. We compute the normal of the solid
bodies and if the component of the ghost velocity in the direction of the
normal is into the solid, we subtract that component,
\begin{equation}
  \label{eq:vg-no-normal}
  \ten{u}_w = \ten{u}_w - (\ten{u}_w \cdot \ten{\hat{n}}) \ten{\hat{n}}
\end{equation}
Note that this is not done when the direction of the dummy particle velocity
is away from the solid body.

When solving the PPE, the boundary condition on $\ten{u}^*$ is typically set
such that the walls do not have any slip. However, we find that when
simulating high-Reynolds number flows this introduces unnecessary divergence
on the fluid particles leading to noisy results. This can be mitigated by
using a slip velocity for the $\ten{u}^*$. This is illustrated for the case of
the lid-driven cavity problem in section~\ref{sec:cavity}.

\subsubsection{Computation of normals}
\label{sec:normals}

The normal is itself computed purely using the particle positions in a novel
way by first computing the following quantity for the solid particles,
\begin{equation}
  \label{eq:normal-approx}
  \ten{n}^*_i = \sum_j -\frac{m_j}{\rho_j} \nabla_i W_{ij}
\end{equation}
If the magnitude of the resulting vector is less than $\frac{1}{4h_i}$, then
the $\ten{n}^*$ is set to zero otherwise we normalize the vector by its
magnitude.  Then, we smooth these normals using an SPH approximation,
\begin{equation}
  \label{eq:normal}
  \ten{n}_i = \sum_j \frac{m_j}{\rho_j}\ten{n}^*_j W_{ij}
\end{equation}
These normal vectors are then made into unit normals. This produces smooth
normals with the positions of the particle alone.

\subsection{Inlet and outlet boundary conditions}
\label{sec:io-bc}

It is relatively straightforward to simulate wind-tunnel type problems with
our scheme as well.  This requires that the inlet and outlet boundary
conditions be set appropriately.

For the inlet, we set the prescribed velocity, this is often a constant or a
particular velocity profile.  The pressure of the inlet region is also solved
for, based on the iterative pressure-Poisson equation.

For the outlets we use a modified ``do-nothing'' type condition as elaborated
in \citet{pawan:inlet-outlet:2019}. As a fluid particle enters the outlet
region, its properties $\ten{u}$ and $p$ are frozen. The particle is advected
in the outlet region using a Shepard extrapolated velocity from the current
fluid particles along the direction of the outlet. For example if the outlet
is oriented perpendicular to the $x$-axis, then we take the fluid's
$x$-component of velocity $\ten{u}$ and perform a Shepard interpolation of
that on the fluid at each timestep and move the outlet particles. This simple
approach works very well as seen from the results presented later.

\section{Results and discussion}
\label{sec:results}

We simulate several standard benchmarks illustrating the proposed method.
These benchmarks include internal flows with exact solutions, external flows
with a free-surface, two and three dimensional cases and finally the flow past
a circular cylinder. We also demonstrate the performance of the new scheme on
a GPU.~Where possible we compare our results with those produced using the
traditional ISPH implementations, and with weakly-compressible SPH using the
TVF~\cite{Adami2013} or EDAC-SPH~\cite{edac-sph:cf:2019} formulations. Our
implementation of the scheme is made using the PySPH
framework~\cite{PR:pysph:scipy16,pysph}. In the interest of reproducible
research, all the code for producing the results in this manuscript are made
available, furthermore every figure is automatically generated using a simple
automation framework~\cite{pr:automan:2018}. The code is available at
\url{https://gitlab.com/pypr/sisph}.

\subsection{Taylor Green Vortex}
\label{sec:tgv}

The Taylor-Green Vortex problem in two-dimensions is periodic and has an exact
solution given by,
\begin{align}
  \label{eq:tgv}
  u &= - U e^{bt} \cos(2 \pi x) \sin(2 \pi y) \\
  v &=   U e^{bt} \sin(2 \pi x) \cos(2 \pi y) \\
  p &=  -U^2 e^{2bt} (\cos(4 \pi x) + \cos(4 \pi y))/4,
\end{align}
where $U=1m/s$ and $b=-8\pi^2/Re$, $Re=U L/\nu$ and $L=1m$.

We simulate this problem by setting the initial condition at $t=0$ for a given
$Re$.  We choose $Re=100$ for our initial tests and compare the results with
the exact solution.  We use the quintic spline with $h/\Delta x = 1.0$.

The decay rate of the velocity is computed by plotting the maximum velocity
$|\ten{u}_{\max}|$ at each time. We compute the $L_1$ error in the velocity
magnitude as,

\begin{equation}
  \label{eq:tg:l1}
  L_1 = \frac{\sum_i |\ten{u}_{i, computed}| - |\ten{u}_{i, exact}|}
  {\sum_i |\ten{u}_{i, exact}|},
\end{equation}

where $\ten{u}_{i, exact}$ is found at the position of the $i$'th particle.
We compute the $L_1$ error in the pressure using,
\begin{equation}
  \label{eq:tgv_p_l1}
  p_{L_1} = \frac{\sum_i |p_{i, computed} - p_{i, avg} - p_{i, exact}|}
  {\text{max}_{i} (p_{i, exact})}.
\end{equation}
Here $p_{i,avg}$ is the average pressure due to the particle and its
neighbors. This is done to avoid any constant pressure bias in any of the
schemes.

Since the Taylor-Green problem has an exact solution, we also take the
opportunity to test various parameters in the new scheme. In particular we
explore variations in the following,
\begin{itemize}
\item Comparison with the original matrix-based formulation of ISPH, with and
    without the use of shifting, vs.\ GTVF and perturbation.
\item Different tolerance value, $\epsilon$.
\item The use of the symmetric form of the pressure gradient.
\item Two different resolutions with different Reynolds numbers.
\end{itemize}

We then compare the results with the WCSPH scheme, the IISPH scheme and the
EDAC scheme. We add a small amount of perturbation (of at most $\Delta x/5$),
in the initial position of the particles for these schemes to ensure uniform
particle distribution. But, the mass and density remain unchanged.

\subsubsection{Comparison with ISPH}

We compare the proposed scheme hereinafter called as SISPH, to the original
matrix form of ISPH~\cite{sph:psph:cummins-rudman:jcp:1999}, we use a particle
shifting technique~\cite{acc_stab_xu:jcp:2009} to maintain uniform particle
distribution, as without particle shifting the ISPH method is unstable for
higher Reynolds number due to particle disorder~\cite{acc_stab_xu:jcp:2009}.
Here the tolerance is set to, $\epsilon = 10^{-2}$. The ``asymm'' form of
pressure gradient is used. For the pressure coefficient matrix in ISPH we
construct a sparse matrix, and use biconjugate gradient stabilized method from
SciPy~\cite{scipy} sparse matrix library to solve the pressure-Poisson matrix.

\begin{figure}[!h]
  \centering
  \begin{subfigure}{0.48\textwidth}
    \centering
    \includegraphics[width=1.0\textwidth]{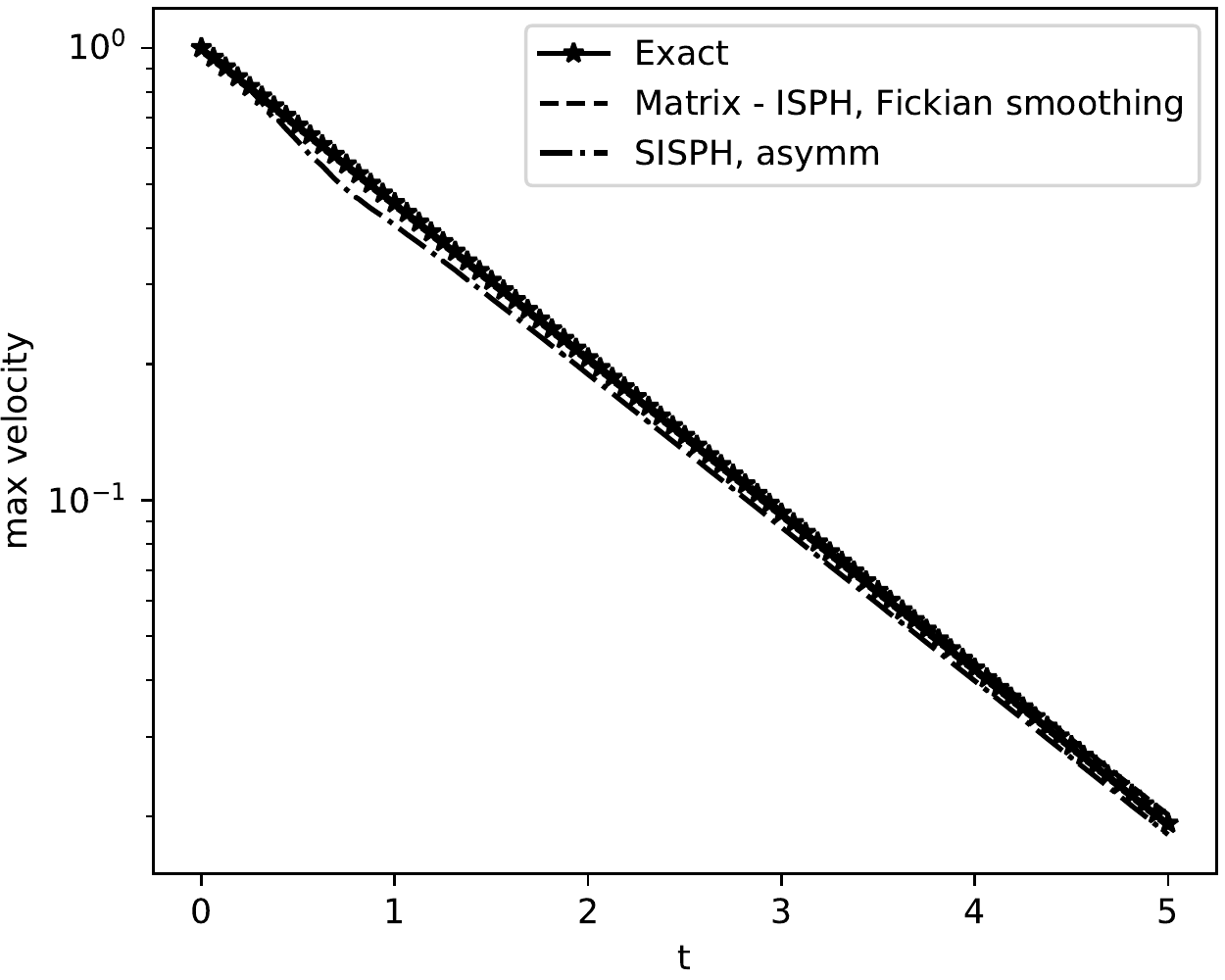}
    \subcaption{Maximum velocity decay vs time.}\label{fig:tg:gtvf:decay}
  \end{subfigure}
  \begin{subfigure}{0.48\textwidth}
    \centering
    \includegraphics[width=1.0\textwidth]{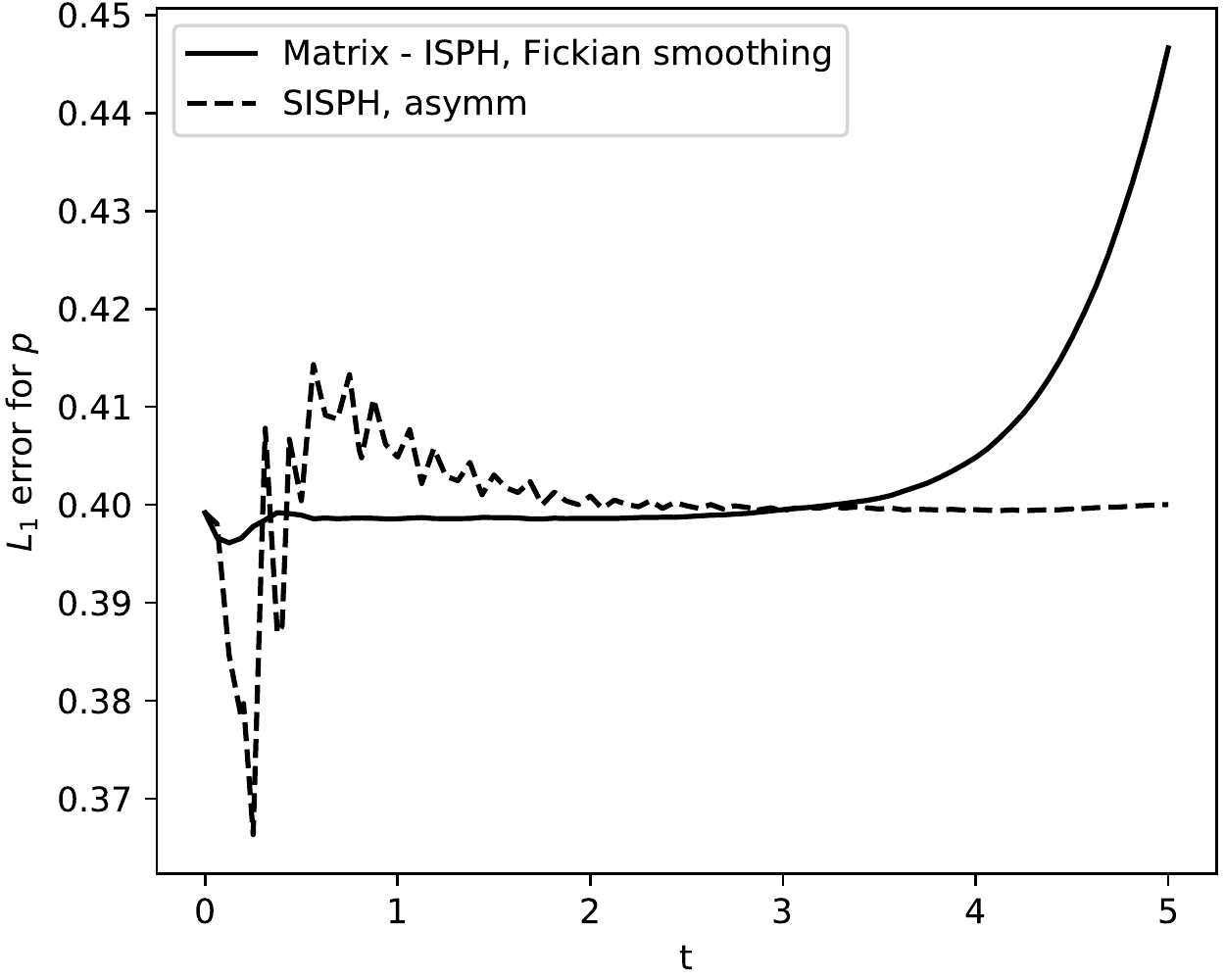}
    \subcaption{${p_{L_1}}$ vs time.}\label{fig:tg:gtvf:l1}
  \end{subfigure}
  \caption{SISPH scheme compared with Matrix based ISPH method with shifting.}
\label{fig:tg:gtvf}
\end{figure}

\begin{figure}[!h]
  \centering
    \includegraphics[width=1.0\textwidth]{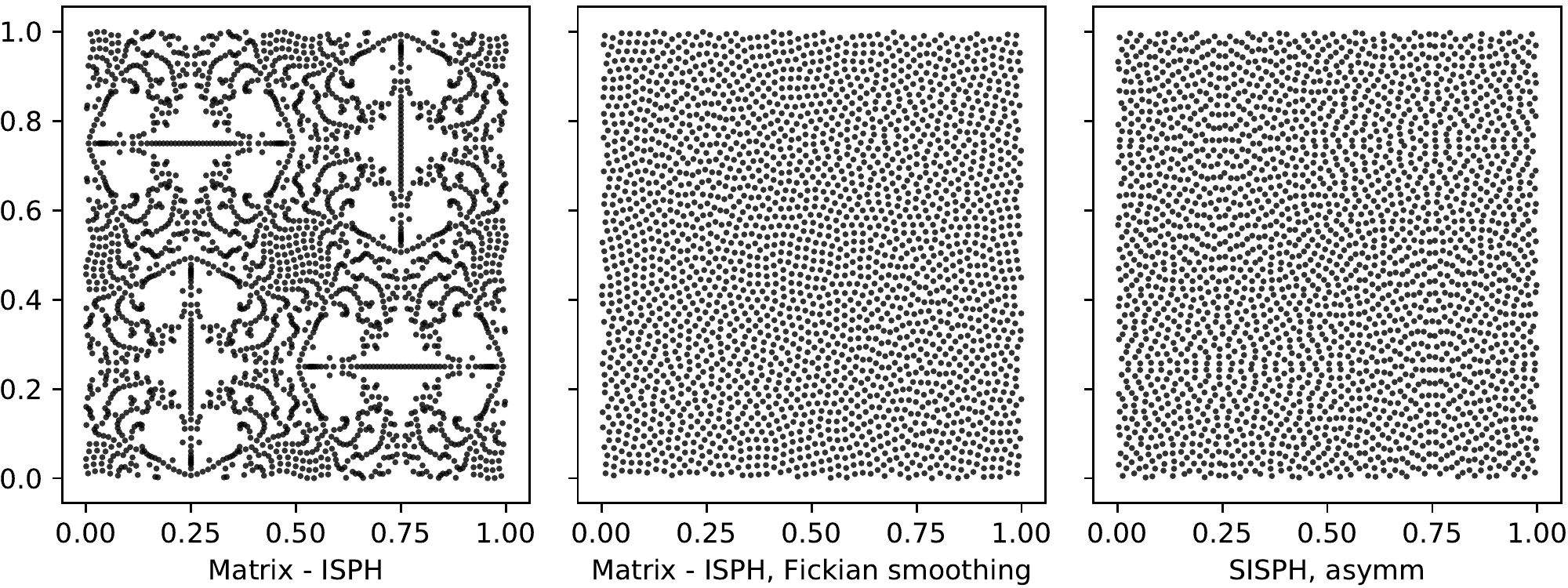}
    \caption{Particle plots for ISPH without shifting, ISPH with shifting, and
        SISPH scheme at $t = 1.2s$. ISPH without shifting is unstable, whereas
        shifting, and SISPH shows uniform distribution.}
\label{fig:tg:pplots}
\end{figure}

In Fig.~\ref{fig:tg:gtvf}, the decay of the velocity and the error in the
pressure are plotted. As can be seen, the decay in velocity of the new scheme
is close to the exact solution. The pressure plots also reveal that the new
scheme performs very well. Fig.~\ref{fig:tg:pplots} shows the particle plots
for ISPH, ISPH with Fickian diffusion based
shifting~\cite{diff_smoothing_sph:lind:jcp:2009,
  fickian_smoothing_sph:skillen:cmame:2013}, and SISPH schemes. As said before
ISPH is unstable without shifting~\cite{acc_stab_xu:jcp:2009}. ISPH with
Fickian based shifting and SISPH both result in a uniform distribution of
particles. This shows that the new scheme performs better than the original
ISPH and is comparable with the case of ISPH along with shifting.


\subsubsection{Change in form of pressure gradient.}

\begin{figure}[!h]
  \centering
  \begin{subfigure}{0.48\textwidth}
    \centering
    \includegraphics[width=1.0\textwidth]{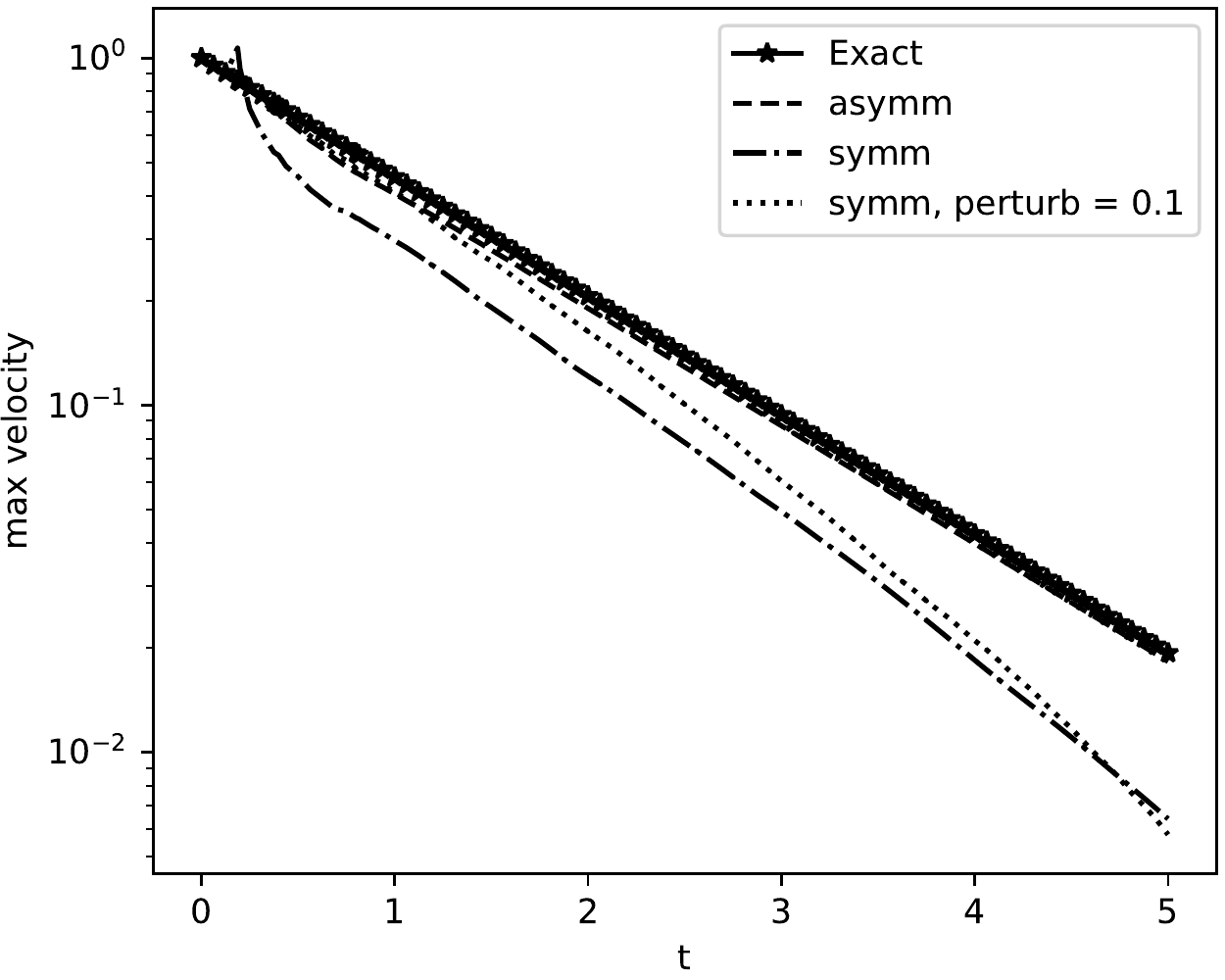}
    \subcaption{Maximum velocity decay vs time.}\label{fig:tg:sym:decay}
  \end{subfigure}
  \begin{subfigure}{0.48\textwidth}
    \centering
    \includegraphics[width=1.0\textwidth]{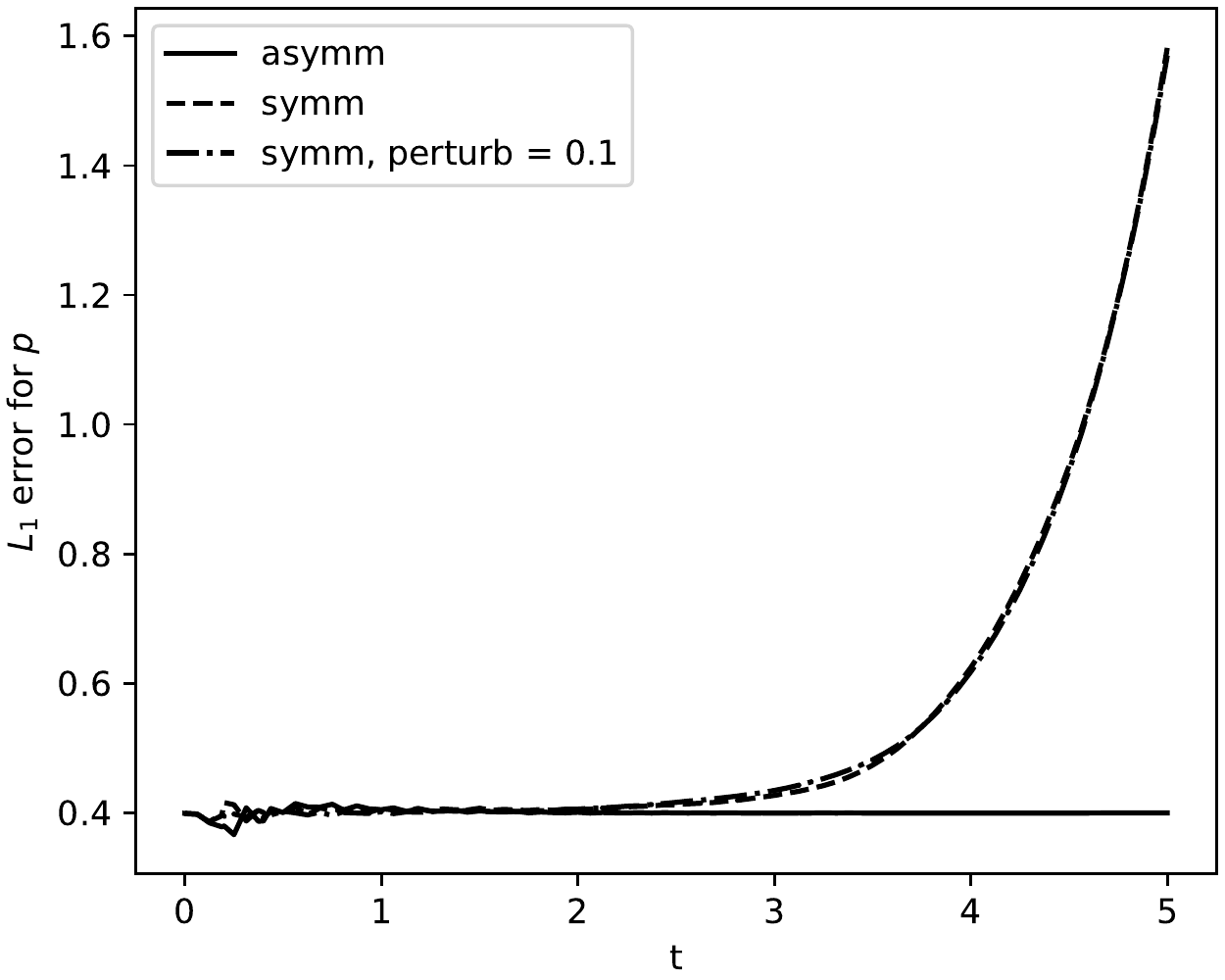}
    \subcaption{$p_{L_1}$ vs time.}\label{fig:tg:sym:l1}
  \end{subfigure}
  \caption{The use of ``symm'' form of pressure gradient vs ``asymm'' form
    of pressure gradient, in this simulation ``asymm'' is better than
    ``symm''.}
\label{fig:tg:sym}
\end{figure}

Here we observe the effect of the two ways of computing the pressure gradient.
The tolerance used for iteration is $\epsilon =10^{-2}$. We also note that the
``symm'' form of pressure gradient works well even without the use of the GTVF
formulation, but the ``asymm'' pressure gradient does not work without the
addition of the GTVF. However, the ``asymm'' pressure gradient with the GTVF
produces very good results, as can be seen in Fig.~\ref{fig:tg:sym}. These
result are better than the ``symm'' form for the pressure and decay rates with
and without an initial perturbation. Hence for the Taylor-Green problem we use
the ``asymm'' form to compute pressure gradient.

\subsubsection{Change in tolerance}

\begin{figure}[!h]
  \centering
  \begin{subfigure}{0.48\textwidth}
    \centering
    \includegraphics[width=1.0\textwidth]{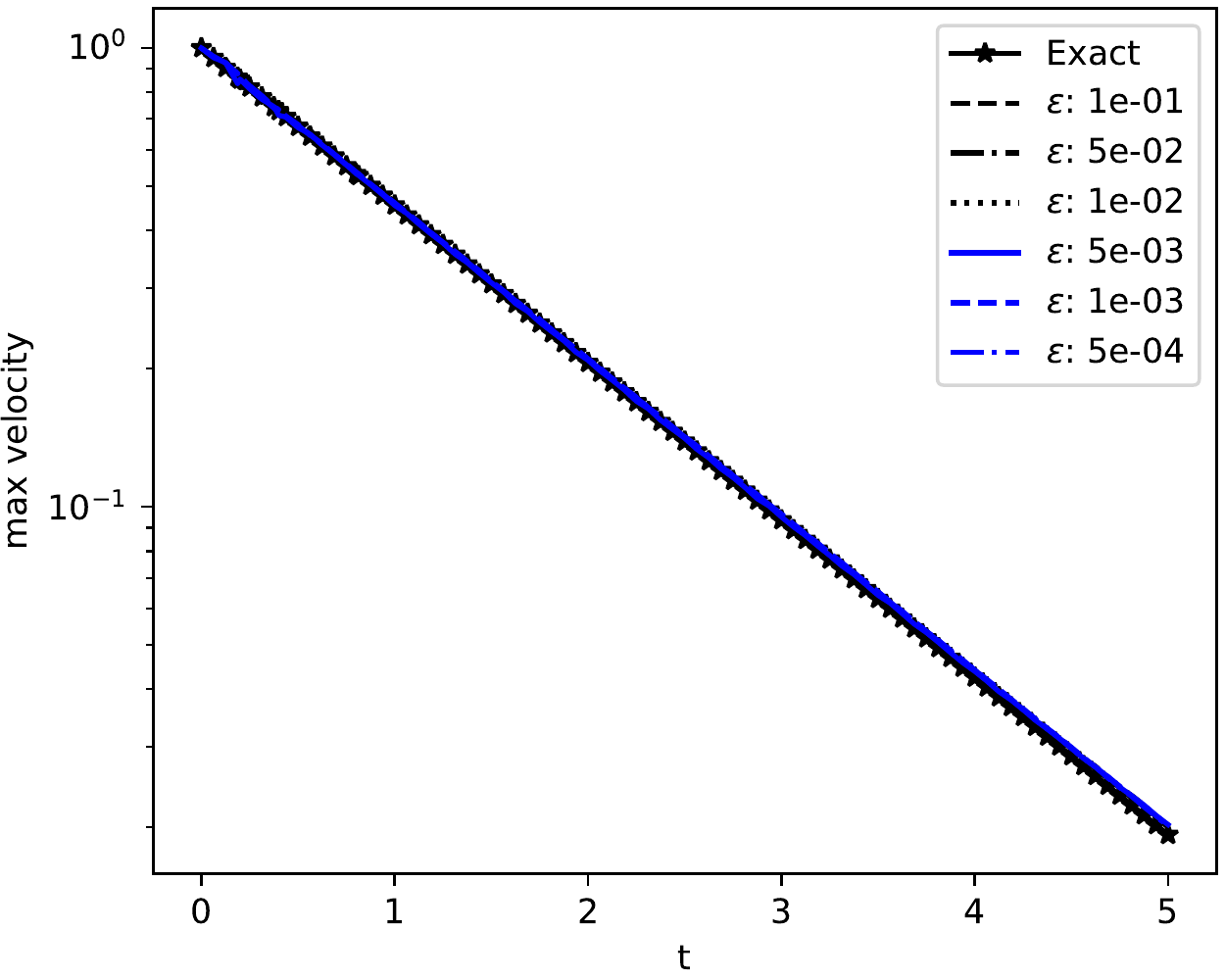}
    \subcaption{Maximum velocity decay vs time.}\label{fig:tg:tol:decay}
  \end{subfigure}
  \begin{subfigure}{0.48\textwidth}
    \centering
    \includegraphics[width=1.0\textwidth]{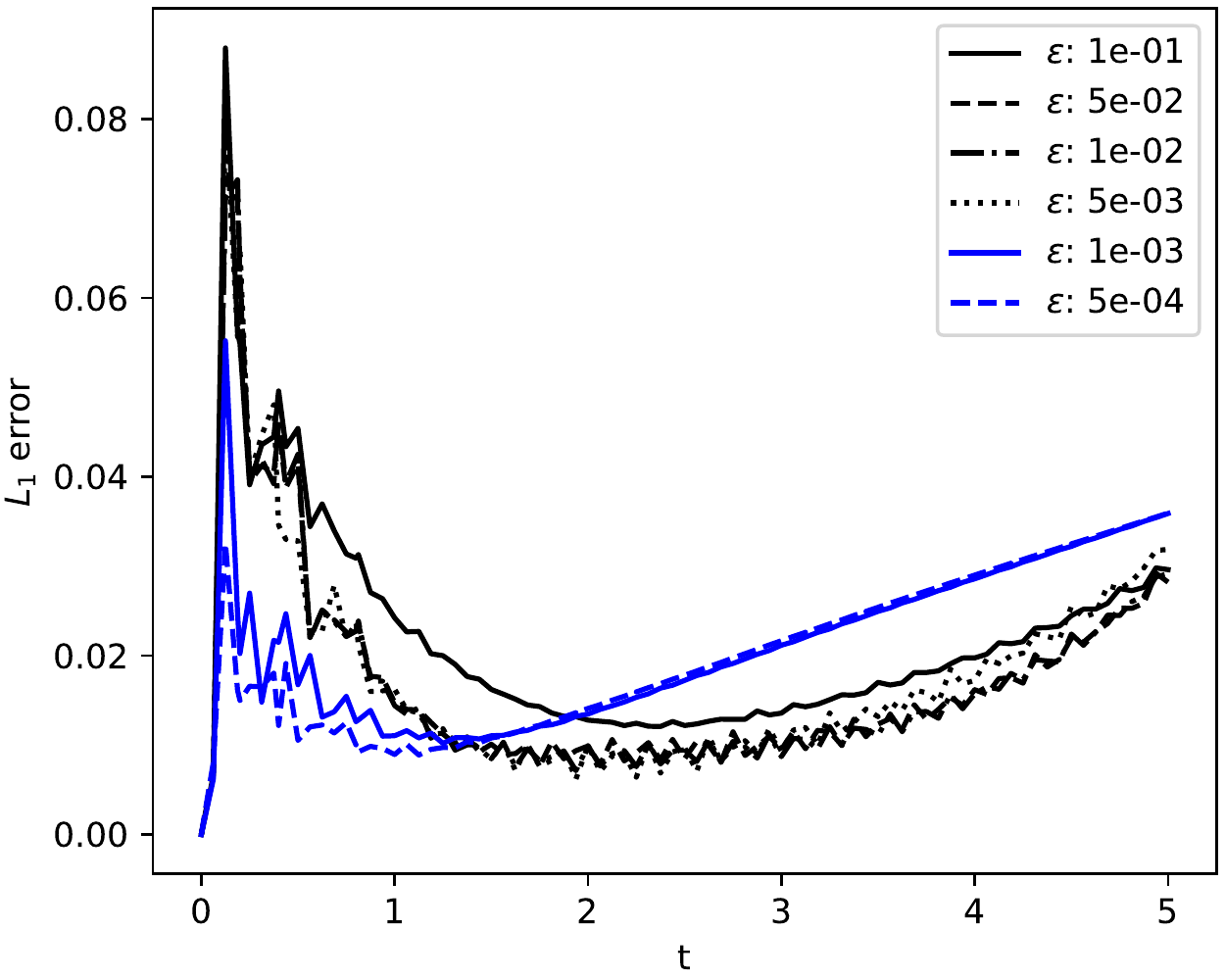}
    \subcaption{$L_1$ error of velocity magnitude vs t.}\label{fig:tg:tol:l1}
  \end{subfigure}
 \caption{Change in tolerance, $\epsilon$, for Taylor-Green problem simulated
  for t = 5s using SISPH formulation.}
\label{fig:tg:tol}
\end{figure}

We next study the variation in the iteration tolerance for the SISPH scheme.
We vary the tolerance in the range
$\epsilon = [0.1, 0.05, 10^{-2}, 5 \times 10^{-3}, 10^{-3}, 5 \times
10^{-4}]$. The ``asymm'' pressure gradient form with a $100 \times 100$ grid
of particles is used. As can be seen in Fig.\ref{fig:tg:tol}, we see that the
results are the same even as we vary the tolerance showing that it is not
necessary to use very low tolerance when the spatial resolution does not
demand it.

\subsubsection{Change in Re with varying resolution.}

\begin{figure}[!h]
  \centering
  \begin{subfigure}{0.48\textwidth}
    \centering
    \includegraphics[width=1.0\textwidth]{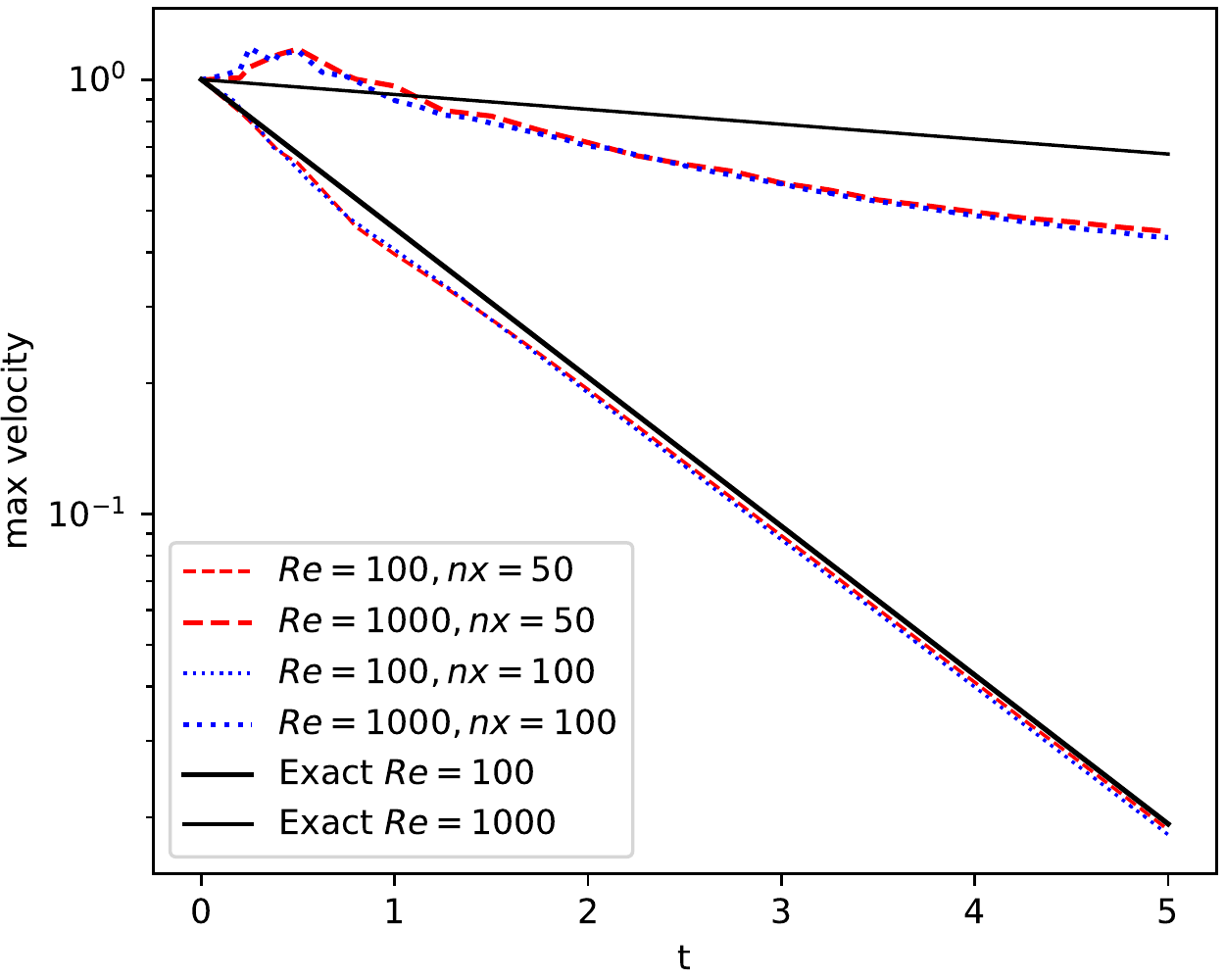}
    \subcaption{Maximum velocity vs time.}\label{fig:tg:re:decay}
  \end{subfigure}
  \begin{subfigure}{0.48\textwidth}
    \centering
    \includegraphics[width=1.0\textwidth]{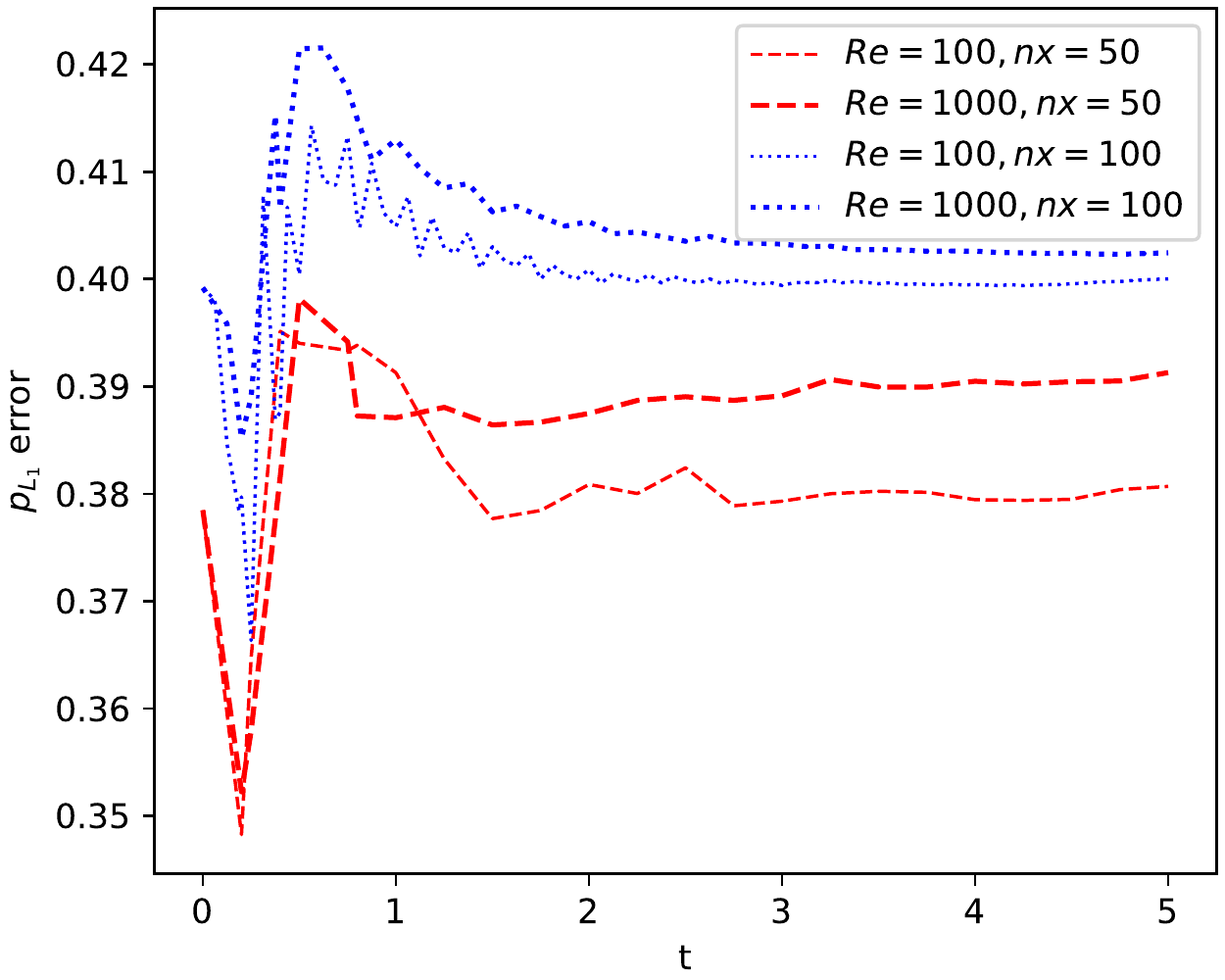}
    \subcaption{$p_{L_1}$ with time.}\label{fig:tg:re:l1}
  \end{subfigure}
  \caption{Taylor-Green problem simulated for t = 5s, with two different
    Reynolds numbers $Re = [100, 1000]$ and initial particle distribution of
    $50 \times 50$ and $100 \times 100$, using SISPH with ``asymm'' pressure
gradient form. }
\label{fig:tg:re}
\end{figure}

We next compare the Taylor-Green problem with different Reynolds numbers, $Re
= 100, 1000$, with different resolutions of $50 \times 50$ and $100 \times
100$. The tolerance is, $\epsilon = 10^{-2}$. ``asymm'' pressure gradient is
used here. Fig.~\ref{fig:tg:re} shows the results. These indicate that the
scheme does work well for higher Reynolds numbers.

\subsubsection{Comparison with other schemes}
Here we compare with other established schemes such as WCSPH, IISPH and EDAC.
All the schemes have a perturbed initial condition where the particles are
displaced by a maximum of $\Delta x/5$ randomly about their original position
with a uniform distribution. All the schemes use the same initial particle
distribution. We use a tolerance of $\epsilon = 10^{-2}$ for the SISPH scheme,
``asymm'' form of pressure gradient with a $100 \times 100$ grid of initial
particles.

\begin{figure}[!h]
  \centering
  \begin{subfigure}{0.48\textwidth}
    \centering
    \includegraphics[width=1.0\textwidth]{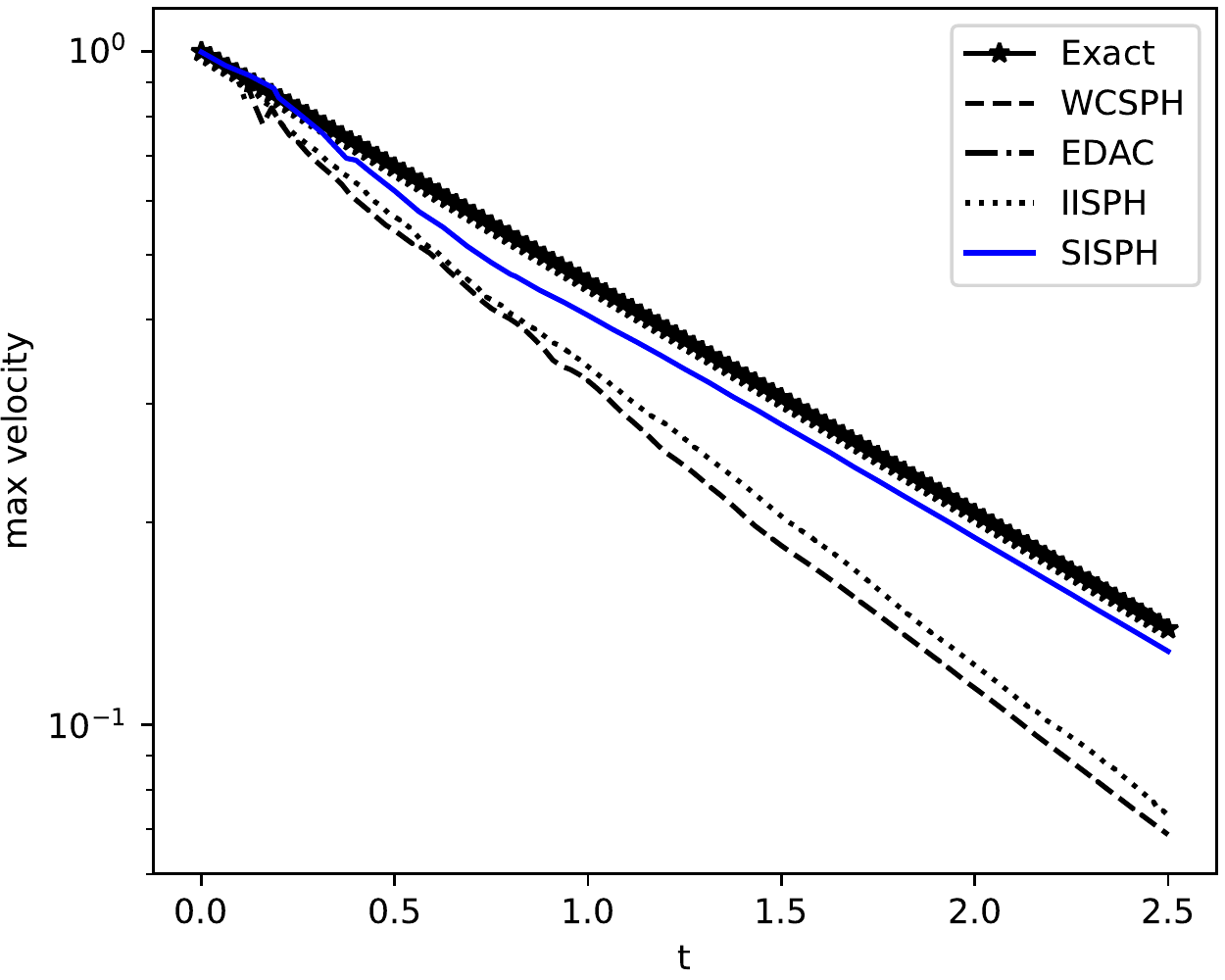}
    \subcaption{Maximum velocity decay vs time.}\label{fig:tg:scheme:decay}
  \end{subfigure}
  \begin{subfigure}{0.48\textwidth}
    \centering
    \includegraphics[width=1.0\textwidth]{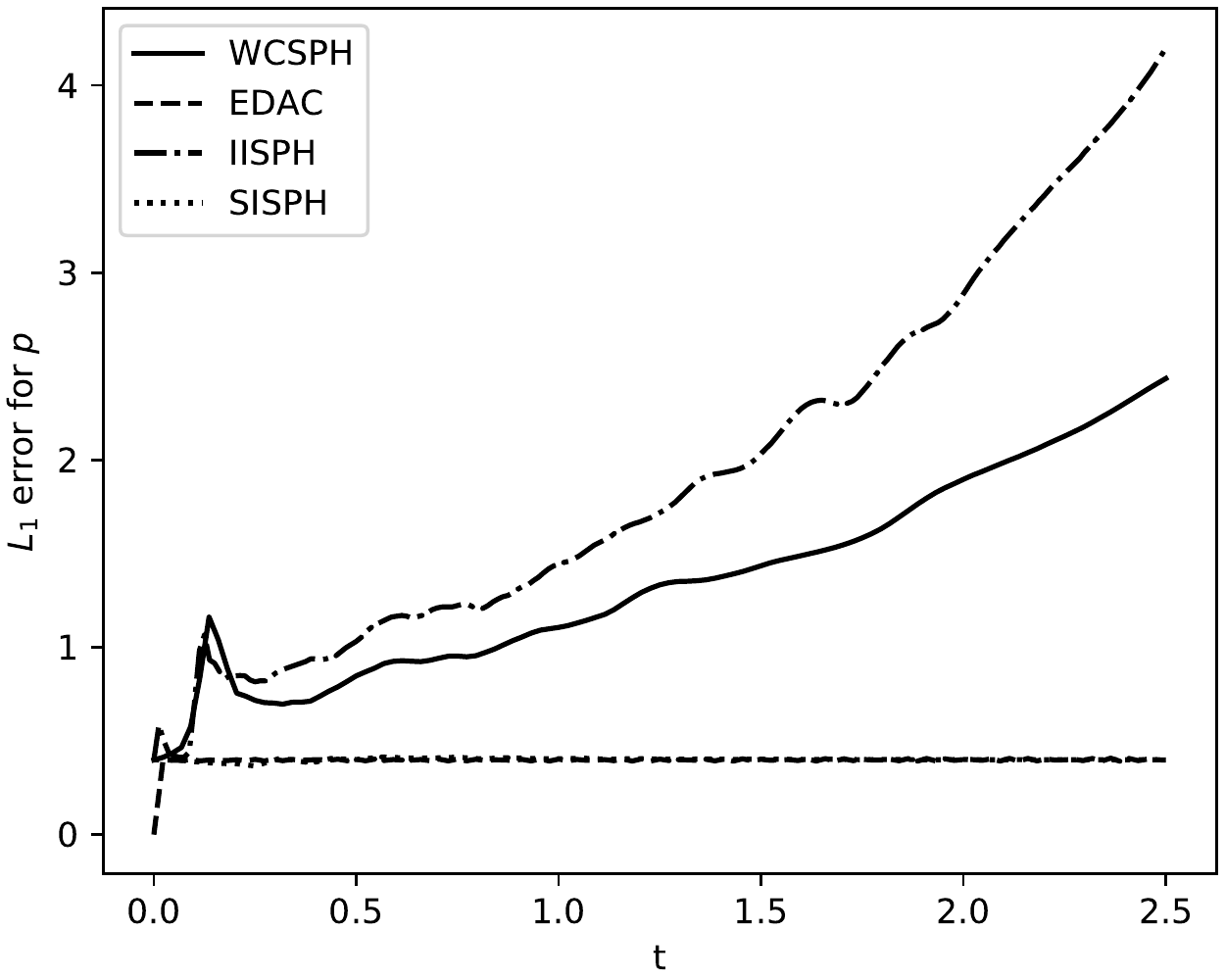}
    \subcaption{$p_{L_1}$ vs time.}\label{fig:tg:scheme:l1}
  \end{subfigure}
  \caption{Comparison of SISPH using ``asymm'' form
    of pressure gradient with WCSPH, EDAC and IISPH schemes. The initial
    configuration is perturbed by a uniformly distributed random number scaled
    by $0.2\Delta x$ for all the schemes and simulated for t = 2.5s.}
\label{fig:tg:scheme}
\end{figure}

In Fig.~\ref{fig:tg:scheme:decay} we note that EDAC and SISPH formulation
matches the exact decay rate, and in the Fig.~\ref{fig:tg:scheme:l1} the
errors in pressure, $p_{L_1}$, are much better compared to WCSPH and IISPH.

In conclusion SISPH with a tolerance of $\epsilon = 10^{-2}$ works well with
the Taylor-Green problem and hence we choose this for simulating the
subsequent problems.

\subsection{Comparison with EISPH}
\label{sec:tg-eisph}

Here we compare the Explicit ISPH
(EISPH)~\cite{nomeritae-eisph-2016,hosseini-2007,rafiee-2009,barcarolo-2014}
with SIPSH. In EISPH, the PPE is solved with an explicit equation which is
equivalent to one single iteration of the Jacobi iterations in each time step.
In SISPH the PPE is solved iteratively to reach a desired tolerance per every
time step. This is important in the case of high Reynolds number simulations.
We compare the results using both the methods for the Taylor-Green problem at
$Re = 10000$. We perturb the initial particle positions by a maximum of
$\Delta x/ 10$ without changing the mass or density.

\begin{figure}[!h]
  \centering
  \includegraphics[width=1.0\textwidth]{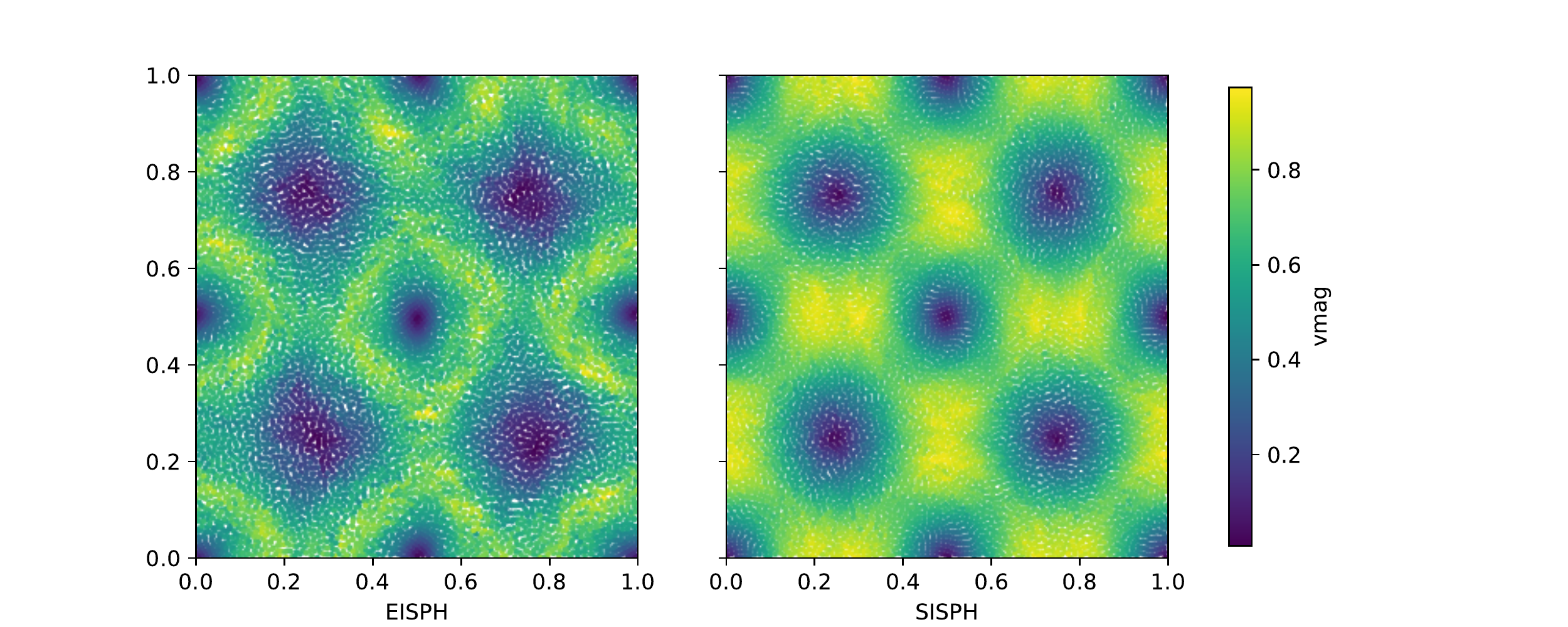}
  \caption{Comparison of EISPH with SISPH for $Re = 10000$, using
    $100 \times 100$ initial particle distribution. Particles are perturbed
    from their initial position by a maximum of $\Delta x/10$. }\label{fig:tg:eisph:pplots}
\end{figure}

Figure~\ref{fig:tg:eisph:pplots} shows the particle plots of EISPH vs SISPH
where the color indicates velocity. This shows that the using EISPH approach
leads to an inaccurate velocity. This shows the importance of multiple
iteration per time step for high Reynolds number flows.

\subsubsection{Performance of SISPH.}

\begin{figure}[!h]
  \centering
  \includegraphics[width=0.6\textwidth]{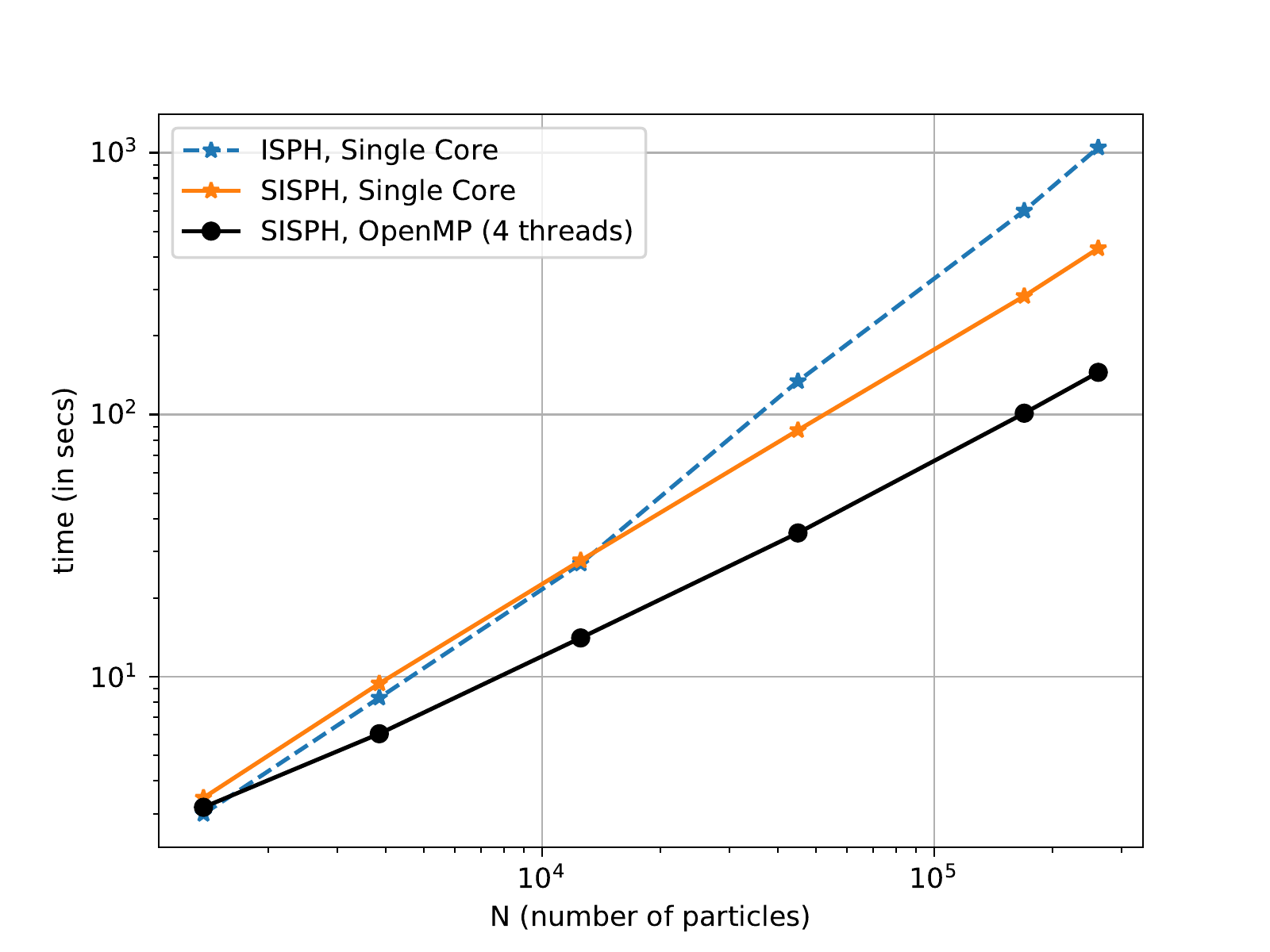}
  \caption{Log-log plot showing time taken vs number of particles for SISPH on
    single and multi-core CPUs, and matrix based ISPH on single-core
    CPU.}
\label{fig:taylor_green:perf}
\end{figure}

\begin{figure}[!h]
  \centering
  \includegraphics[width=0.6\textwidth]{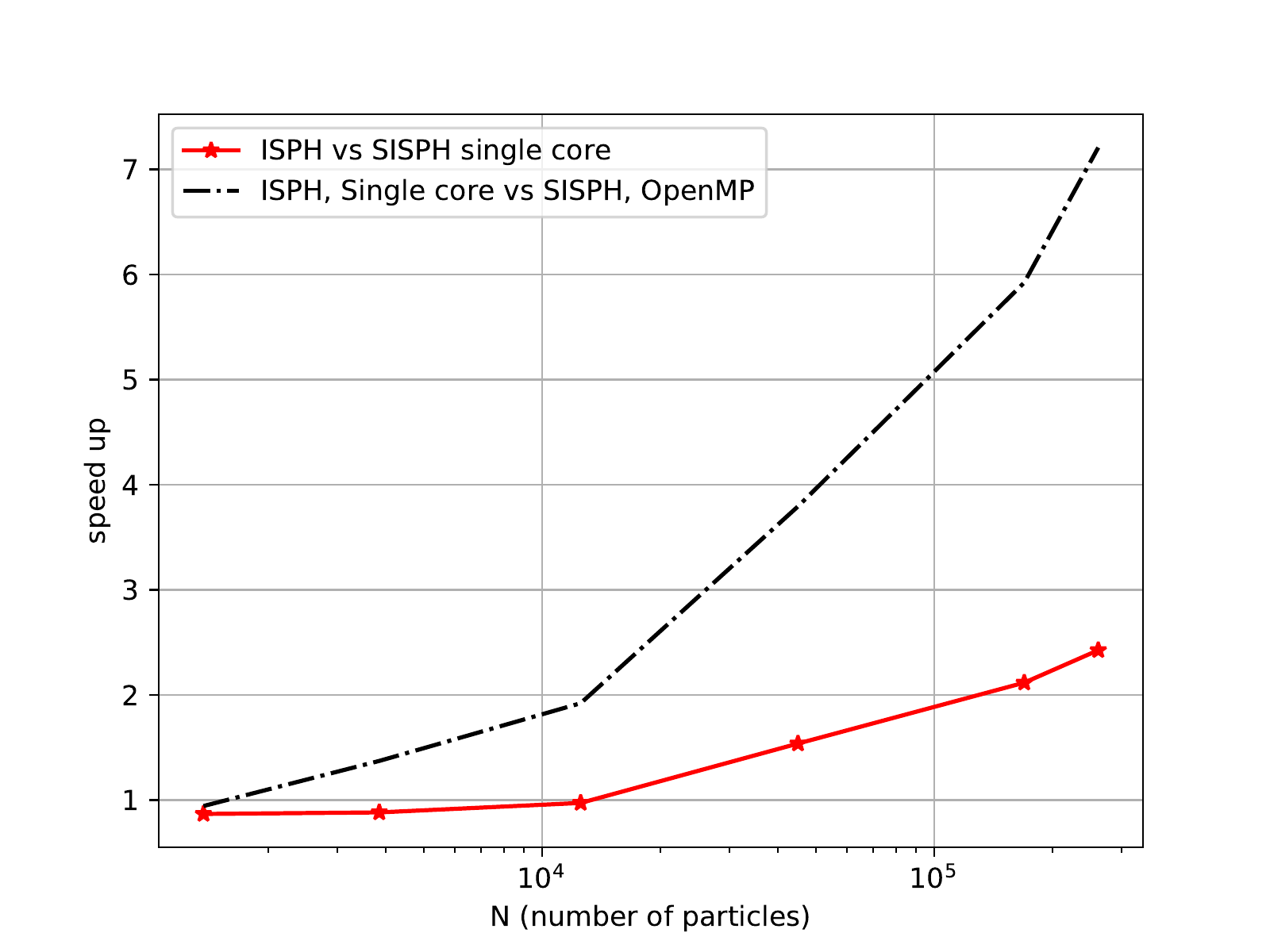}
  \caption{Performance speed up observed for SISPH scheme on single and
      multi-core with ISPH on single-core.}\label{fig:taylor_green:speedup}
\end{figure}

Here we compare the performance of SISPH scheme with ISPH scheme on a CPU for
a single-core and show the performance of SISPH on
multi-cores. Fig.~\ref{fig:taylor_green:perf} shows the time taken versus the
number of particles, $N$, for the different cases.  When $N$ is large the
SISPH method is an order of magnitude faster than the matrix based ISPH
method. Fig.~\ref{fig:taylor_green:speedup} compares the scale up of SISPH
with ISPH scheme on single core, the single core performance of SISPH is
increase to 4 times, and the multi-core (with 4 cores) is nearly 12 times as
fast as single core ISPH for large $N$. The matrix ISPH cannot be run on GPUs
without a significant amount of programming effort whereas SISPH can be
executed on the GPUs very easily as demonstrated for a 3D problem later.

\subsection{Lid Driven Cavity}
\label{sec:cavity}

The next problem considered is the classical lid-driven cavity problem. This
problem allows us to test the solid wall boundary condition and the ability of
the new scheme to handle the solid wall boundaries.

The problem involves a unit square vessel filled with incompressible fluid.
The upper boundary over the fluid is made to move rightward with a unit speed,
$U$. The other walls are treated as no-slip walls. For this problem the
Reynolds number is defined as $Re=UL/\nu$, where $L$ is the length of a side.
We simulate the problem at $Re=100$ with two grid sizes, $50 \times 50$ and
$100 \times 100$. The results are compared with those of~\citet{ldc:ghia-1982}
where the velocity profile along a line passing horizontally and vertically
through the center is shown after the simulation has stabilized. We also
compare the results with those from the TVF~\cite{Adami2013} scheme. We use
the parameter $h/dx = 1.0$ and employ the quintic spline kernel, the tolerance
for SISPH is, $\epsilon = 10^{-2}$, and since this is an internal flow the
background pressure in GTVF formulation is fixed to a constant reference
pressure which is $2 \max(p_i)$ after the first iteration. The ``symm''
form of pressure gradient is used.
\begin{figure}[!h]
  \centering
  \begin{subfigure}[b]{0.6\linewidth}
    \centering
    \includegraphics[width=\linewidth]{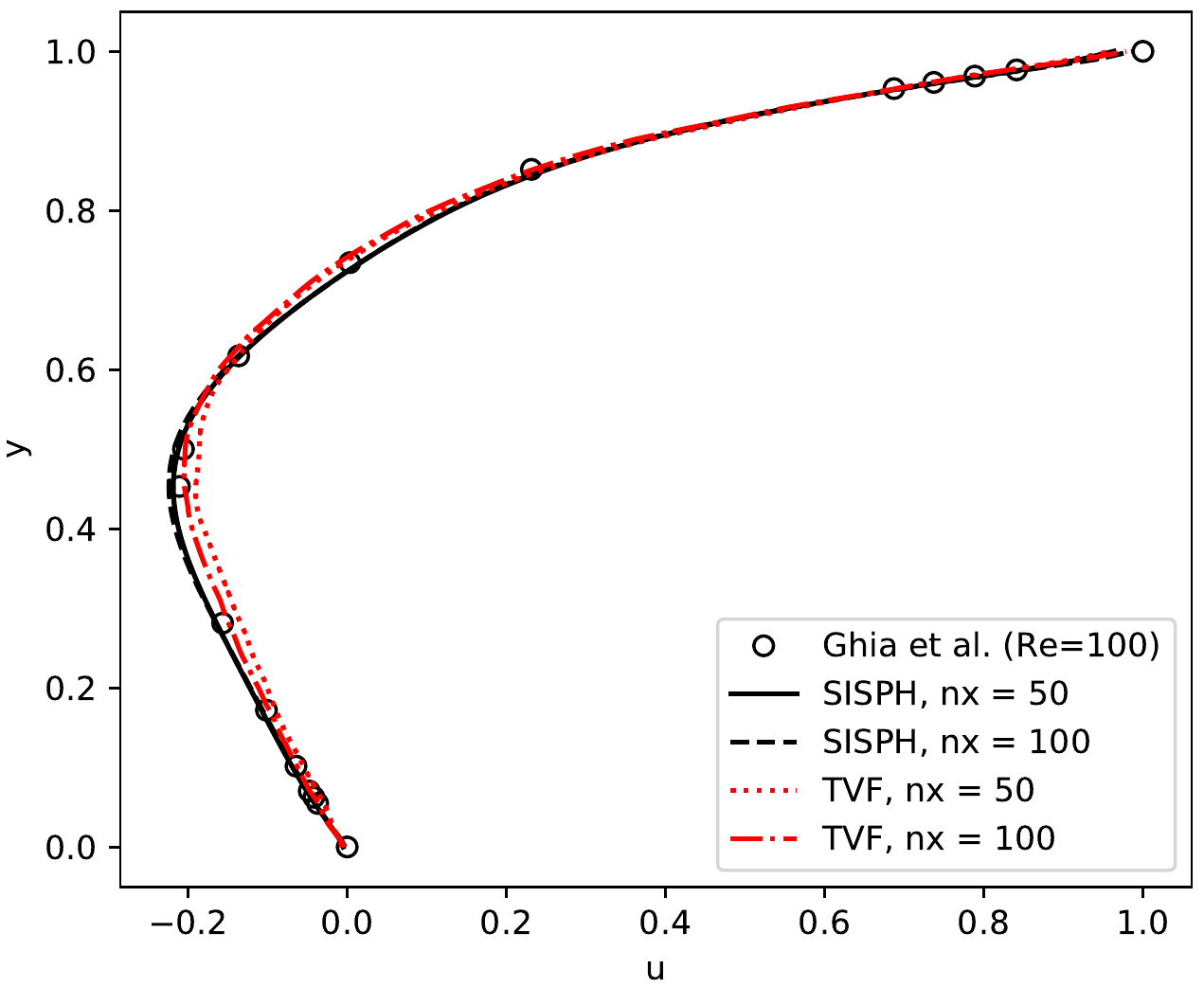}
  \end{subfigure}
\\
  \begin{subfigure}[b]{0.6\linewidth}
  \centering
    \includegraphics[width=\linewidth]{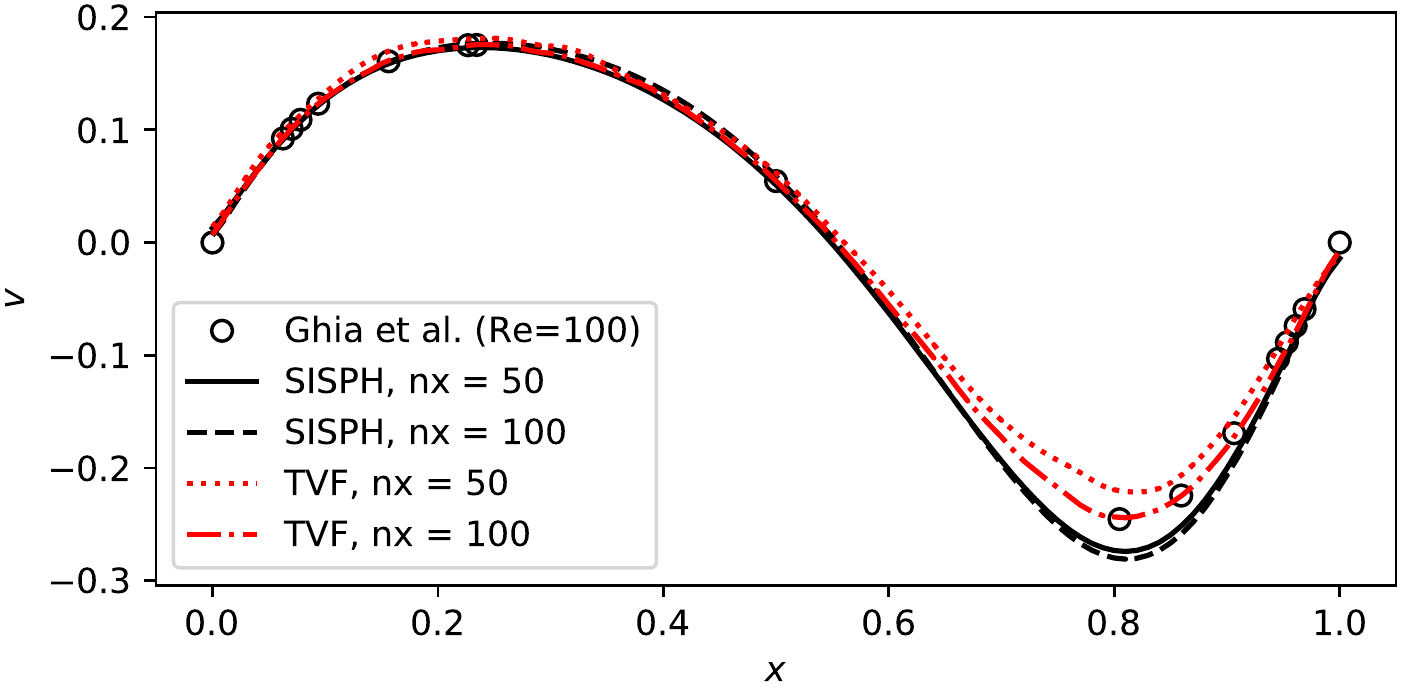}
  \end{subfigure}
  \caption{Velocity profiles $u$ vs.\ $y$ and $v$ vs.\ $x$ for the
    lid-driven-cavity problem at $Re=100$. Two particle discretizations of
    $50 \times 50$ and $100 \times 100$ are shown. We compare the results with
    those of~\citet{ldc:ghia-1982} and with the TVF scheme.}
\label{fig:ldc:uv_re100}
\end{figure}
Fig.~\ref{fig:ldc:uv_re100}, plots the distribution of $u$ and $v$ along the
centerline along the horizontal and vertical. The results are in very good
agreement. In this simulation it is observed that, either use
of ``symm'' form of pressure gradient or the ``asymm'' form doesn't affect
the solution significantly, but the ``symm'' form works well even when no GTVF
is used whereas the ``asymm'' form fails to work without the GTVF.

Due to the larger time steps the ISPH method allows, the new scheme is about 5
times faster than the TVF simulation. The $Re=100$ case with $100\times 100$
fluid particles takes only around 87 seconds to execute for 10 seconds of
simulation time (on a quad core Intel i5-7400) with the new scheme. With the
TVF scheme the same problem takes about 433 seconds. This shows that the new
method works well and is much more efficient.

\begin{figure}[!h]
  \centering
  \includegraphics[width=\linewidth]{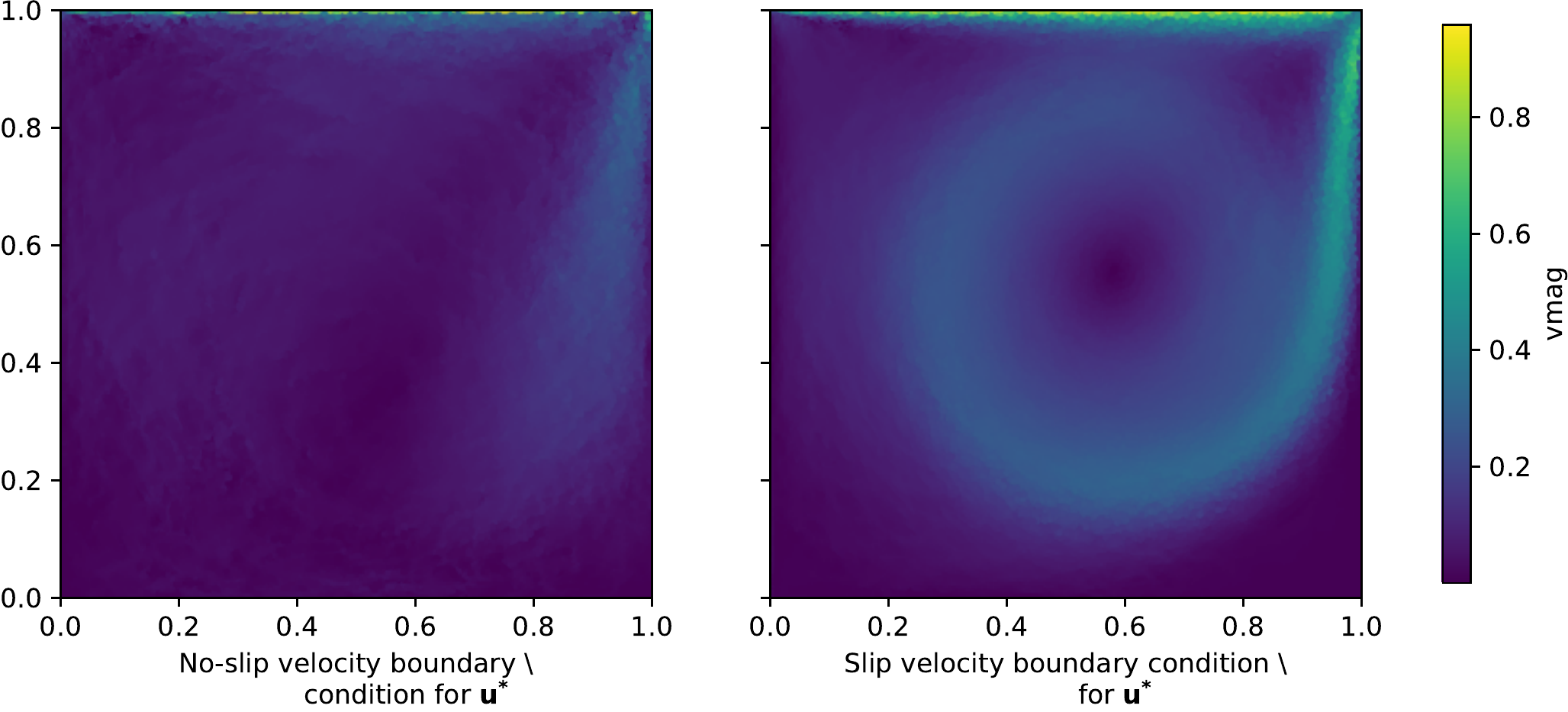}
  \caption{Particle plots indicating the velocity magnitude for the
    lid-driven-cavity problem at $Re = 10000$. Two different boundary
    conditions are used in the plots. On the left the no-slip wall boundary
    condition is used on the intermediate velocity while solving the PPE. On
    the right we use the slip velocity boundary condition.}
\label{fig:ldc:good-vs-bad}
\end{figure}

We next compare the results for $Re=10000$. The velocity profiles show a good
agreement with~\citet{ldc:ghia-1982}. It should however be noted that the
simulation of~\citet{bar13} shows voiding with $Re = 3200$. In our case the
GTVF and an accurate pressure leads to very good results.
\begin{figure}[!h]
  \centering
  \begin{subfigure}[b]{0.6\linewidth}
    \centering
    \includegraphics[width=\linewidth]{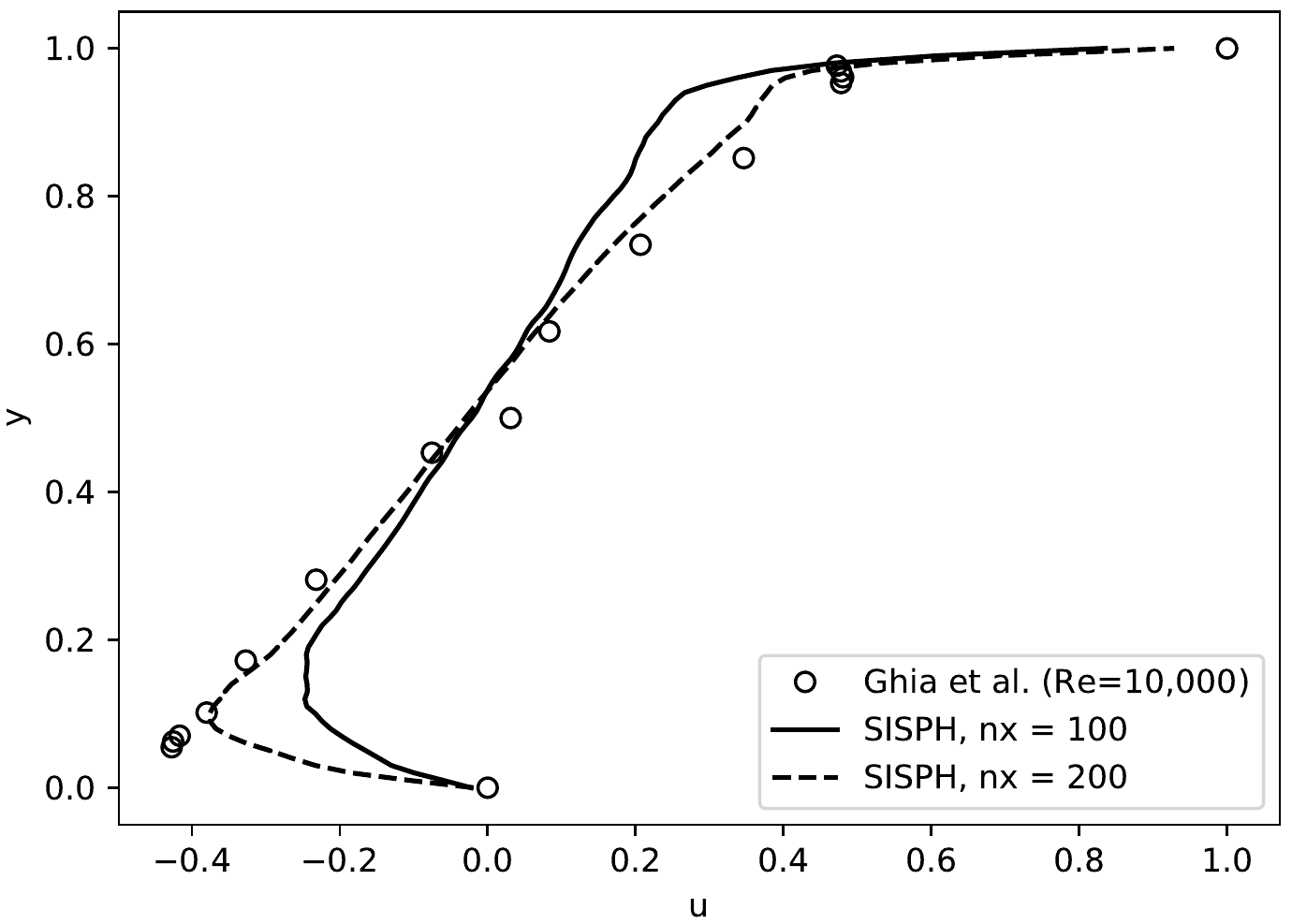}
  \end{subfigure}
\\
  \begin{subfigure}[b]{0.6\linewidth}
  \centering
    \includegraphics[width=\linewidth]{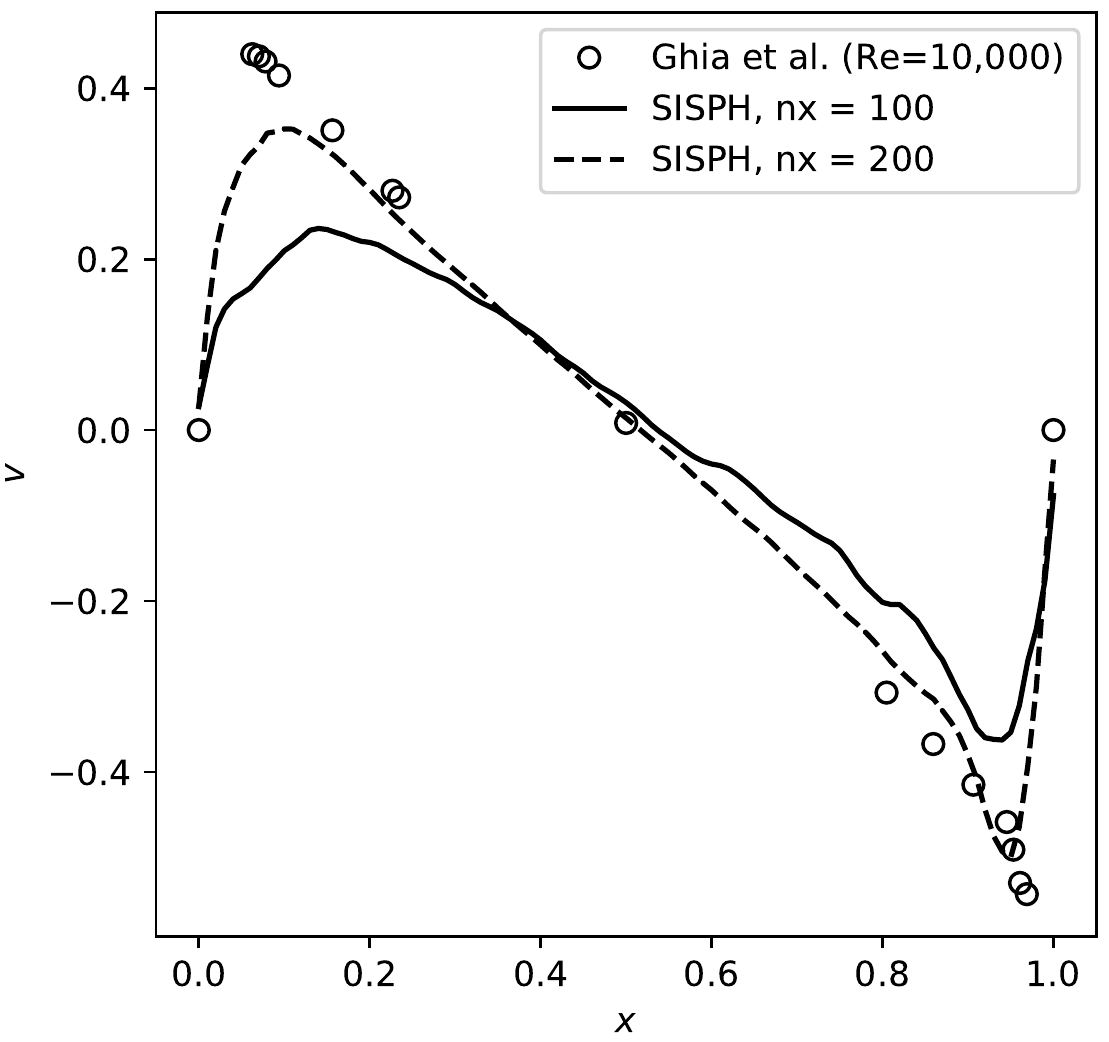}
  \end{subfigure}
  \caption{Velocity profiles $u$ vs.\ $y$ and $v$ vs.\ $x$ for the
    lid-driven-cavity problem at $Re=10,000$. Two particle discretizations of
    $100 \times 100$ and $200 \times 200$ are shown. We compare the results with
    those of~\citet{ldc:ghia-1982}.}
\label{fig:ldc:uv_re10000}
\end{figure}
Fig.~\ref{fig:ldc:good-vs-bad}, plots the velocity magnitude of the cavity
problem at $Re = 10000$. It highlights the difference between the boundary
condition that is applied to the intermediate velocity ($\ten{u}^{*}$) when
solving the PPE.  In the plot on the left, the boundary condition on
$\ten{u}^{*}$ is set to ensure a no-slip wall, whereas, on the right, the
boundary condition is set to ensure a slip wall. The slip wall boundary
condition on the intermediate velocity reduces the noise introduced in the
velocity field due to the unnecessary divergence created at the wall.
\subsection{Evolution of a Square Patch}
\label{sec:square-patch}

The simulation of an initially square patch is a free-surface benchmark
introduced by~\citet{colagrossi-phdthesis:2005}. Here, an initially square
patch of fluid having length $L$ is given the following initial conditions,
\begin{equation}
  \begin{aligned}
    u_{0}(x, y) &= \omega y \\
    v_{0}(x, y) &= -\omega x
  \end{aligned}
\label{eq:sqp:vel}
\end{equation}
\begin{equation}
  p_{0}(x, y) = \rho \sum_m^\infty \sum_n^\infty -\frac{32\omega^2/(m n \pi^2)}
  {\left[{\left(\frac{n\pi}{L}\right)}^2 + {\left(\frac{m \pi^2}{L}\right)}^2
    \right]} \sin\left(\frac{m \pi x^*} {L}\right)\sin\left(\frac{n \pi y^*}
    {L}\right) m, n \in \mathbb{N}_{odd}
\label{eq:sqp:p}
\end{equation}
where $X^* = x + L/2$ and $y^* = y + L/2$.

This problem is simulated for 3 seconds using the new scheme, and the WCSPH
scheme, for the discussion on the results of IISPH scheme
see~\citet{prabhu:dtsph:2019}. The quintic spline kernel with $h/\Delta x
=1.3$ is used for all schemes. An artificial viscosity coefficient of
$\alpha=0.15$ is used for the WCSPH and the new scheme. The tolerance chosen for
the new scheme is $\epsilon=10^{-2}$. Initial particle distribution of $100
\times 100$ is used for all the schemes.

\begin{figure}[!h]
  \centering
  \includegraphics[width=1.0\textwidth]{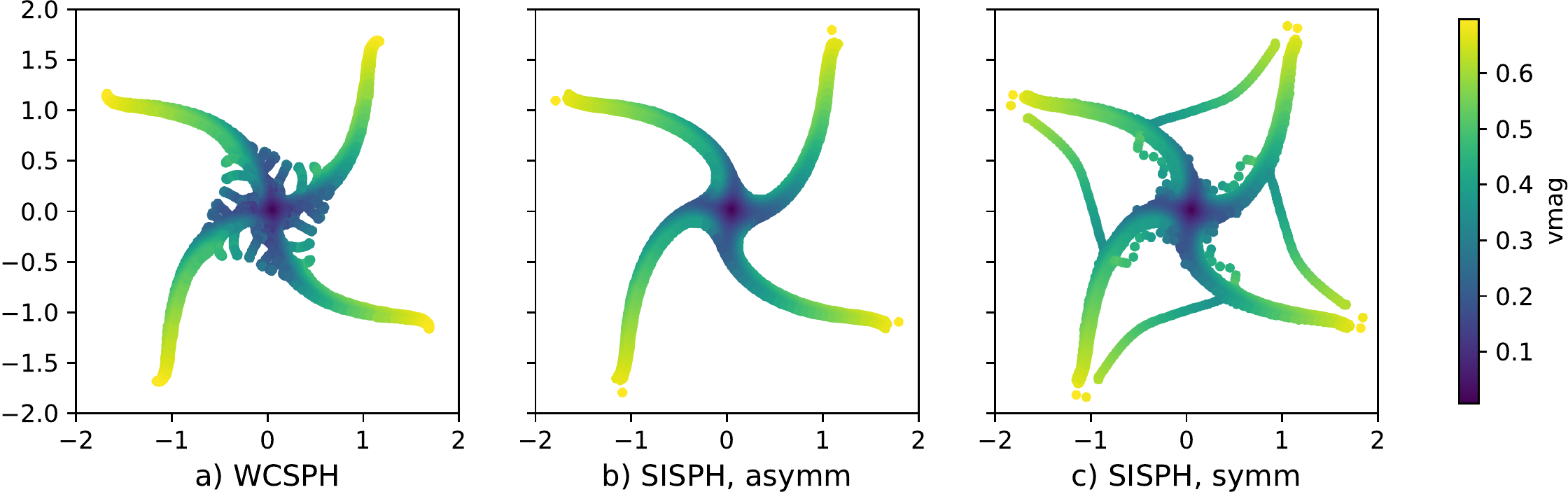}
  \caption{Distribution of particles at $t=3s$ for the initially square patch,
    with a $100 \times 100$ initial particle, distribution. SISPH
    scheme is compared with WCSPH scheme. Effect of different forms
    of pressure gradient is shown. In the ``symm'' pressure gradient particles
    along the free surface are deviated further.}
\label{fig:sqp:particle_plots}
\end{figure}

A noticeable difference is observed when using the different forms of pressure
gradient, the ``asymm'' gradient of pressure simulates the problem
significantly better than the other, where even though the core structure
remains, the particles along the free surface are deviated further. We also
observed that computing the density after each time step using summation
density is very crucial for this problem, without it i.e., using a constant
density through out makes the simulation unstable. A comparison with WCSPH
scheme is made, and the effect of pressure gradient is shown in the
Fig.~\ref{fig:sqp:particle_plots}.  It is clear that the new scheme with the
``asymm'' form of pressure gradient performs very well.

\subsection{Dam break in two-dimensions}
\label{sec:db2d}

The dam-break problem in two-dimensions is a classic problem involving
free-surfaces and solid walls. We simulate the standard problem as described
in~\cite{lee_violeau:db3d:jhr2010}. This involves a vessel of width 4m with a
block of fluid of width 1m and height 2m on the left side. The block of water
is released from rest and falls under the influence of gravity ($g=9.81m/s$).

The problem is simulated with the proposed scheme using a quintic spline
kernel. Artificial viscosity of $\alpha=0.05$ is used for all the schemes.~The
fluid is not treated as viscous. The $h/\Delta x = 1.3$ used for all the
schemes, the default particle spacing, $\Delta x = 0.01$, is used which
results in 37k particles. The time step is fixed and chosen as $\Delta t =
0.125 h/\sqrt{2 g h_w}$, where $g$ is the acceleration due to gravity and
$h_w$ is the initial height of the water which is 2m. The tolerance value for
the new scheme is, $\epsilon = 10^{-2}$. The ``symm'' form of pressure
gradient is used, although no significant difference is observed with either
of the forms.

\begin{figure}[!h]
  \centering
  \begin{subfigure}[b]{0.48\linewidth}
    \includegraphics[width=1.0\linewidth]{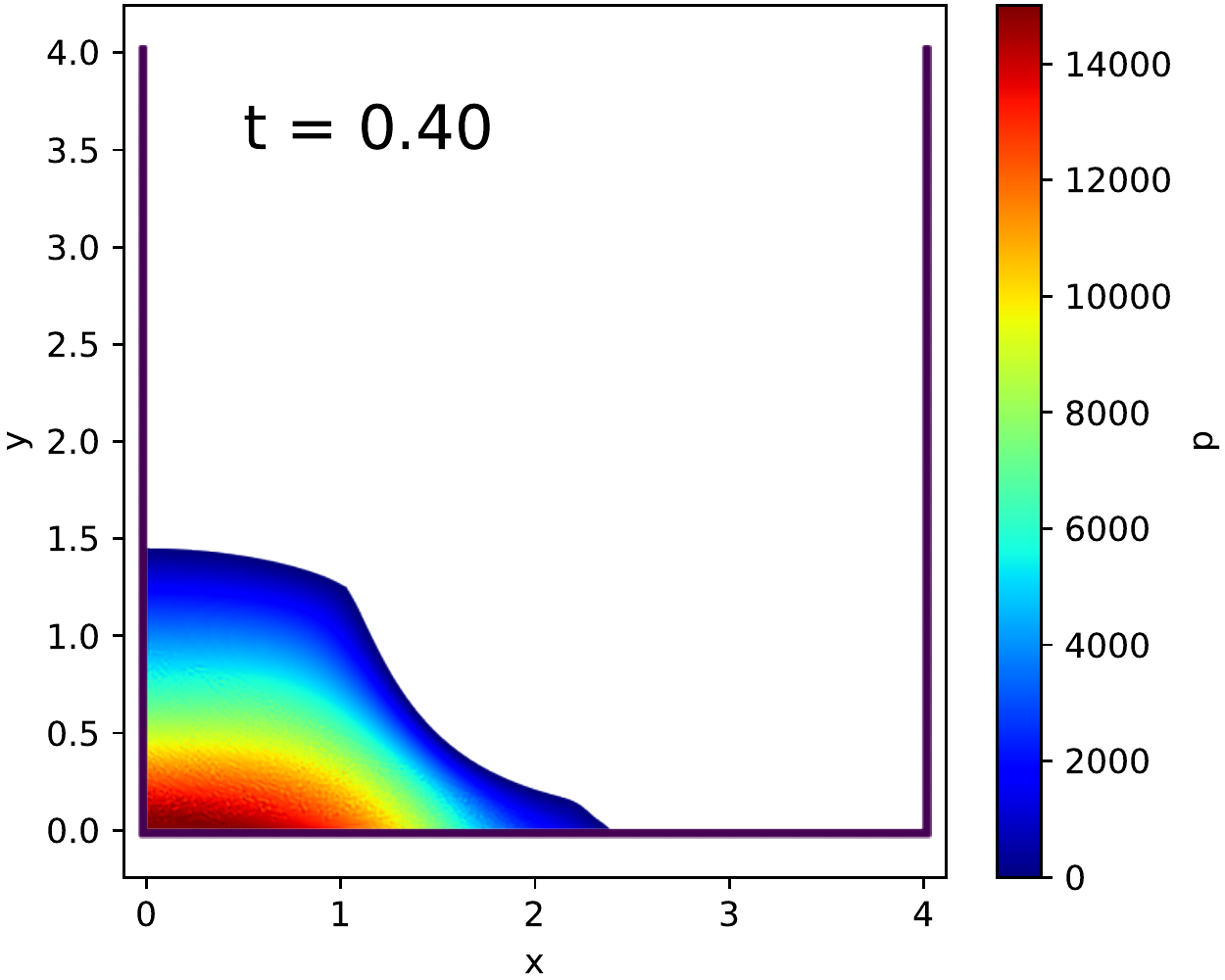}
  \end{subfigure}
  \begin{subfigure}[b]{0.48\linewidth}
    \includegraphics[width=1.0\linewidth]{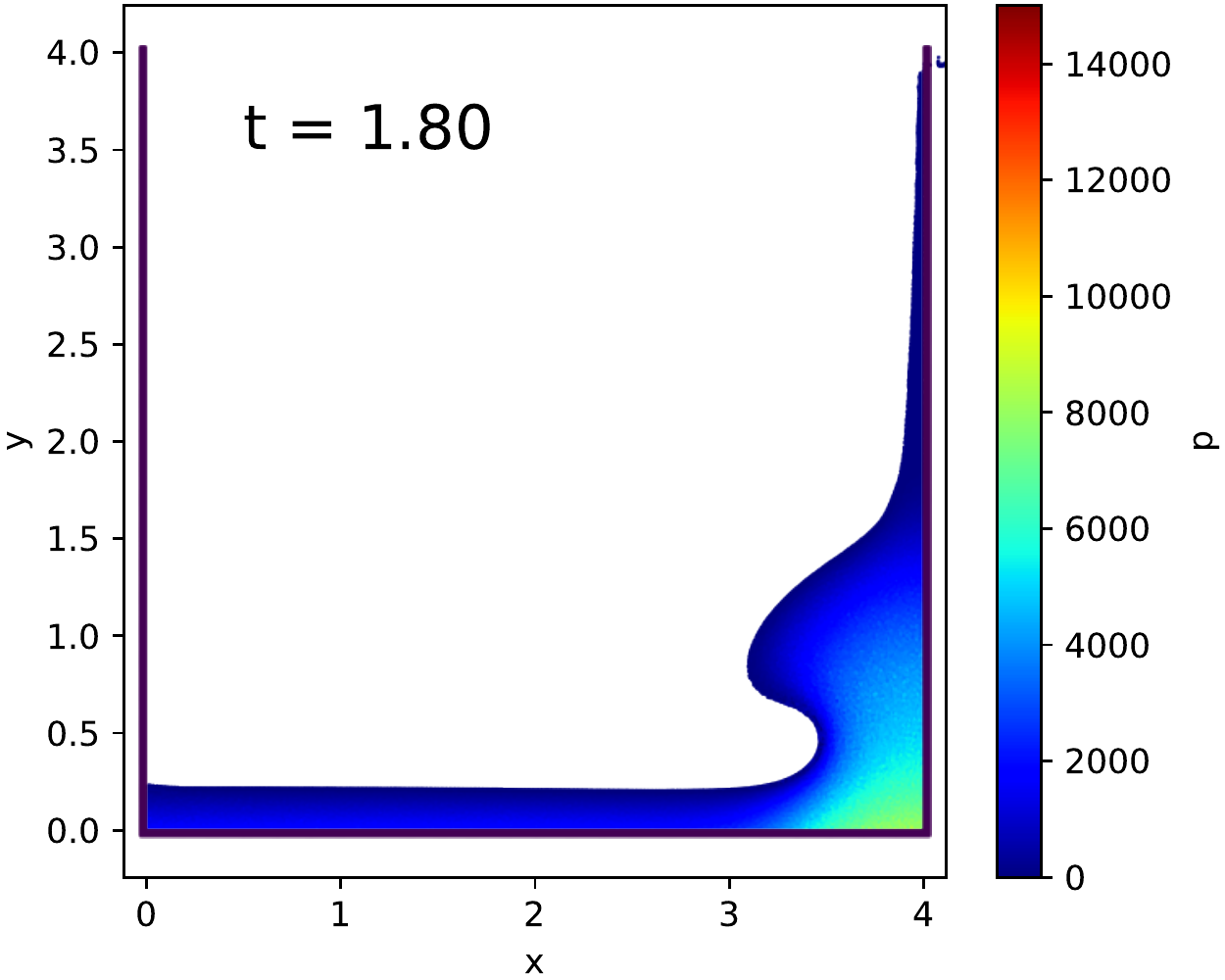}
  \end{subfigure}
  \begin{subfigure}[b]{0.48\linewidth}
    \includegraphics[width=1.0\linewidth]{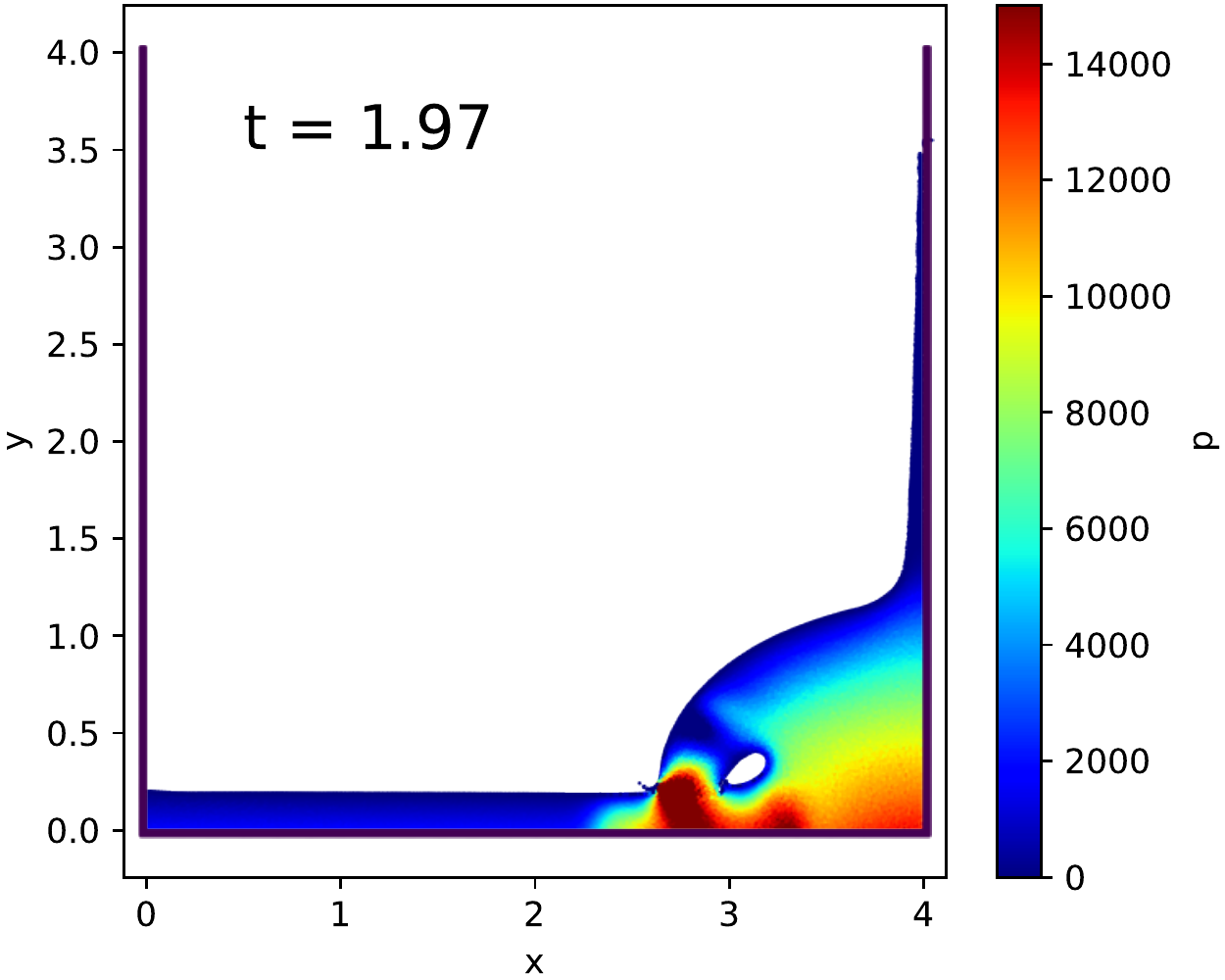}
  \end{subfigure}
  \begin{subfigure}[b]{0.48\linewidth}
    \includegraphics[width=1.0\linewidth]{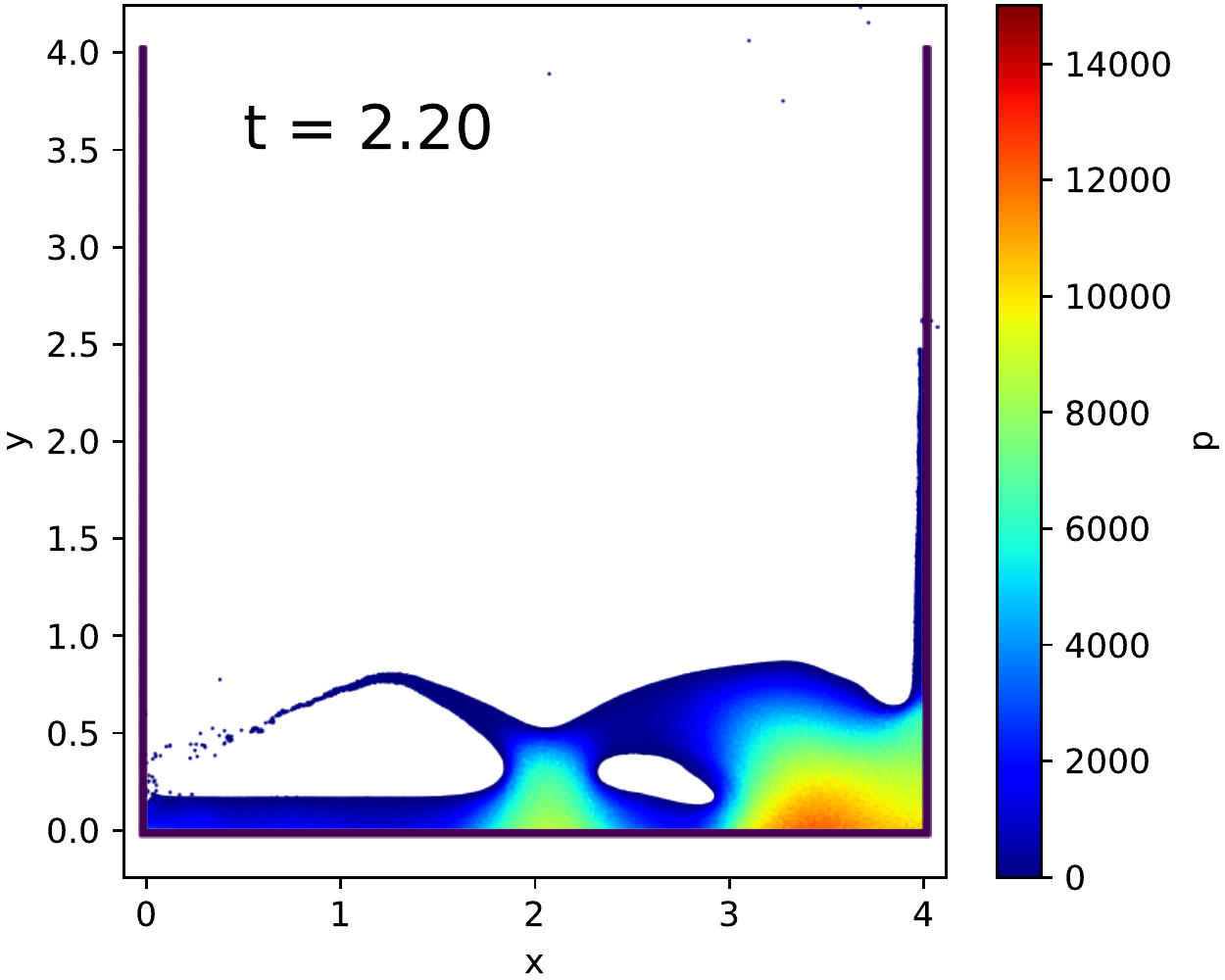}
  \end{subfigure}
  \caption{Particle plots of Dam break in 2D showing pressure at various times
    using SISPH.}
\label{fig:db:particle_plots:pressure}
\end{figure}

\begin{figure}[!h]
  \centering
  \includegraphics[width=0.6\linewidth]{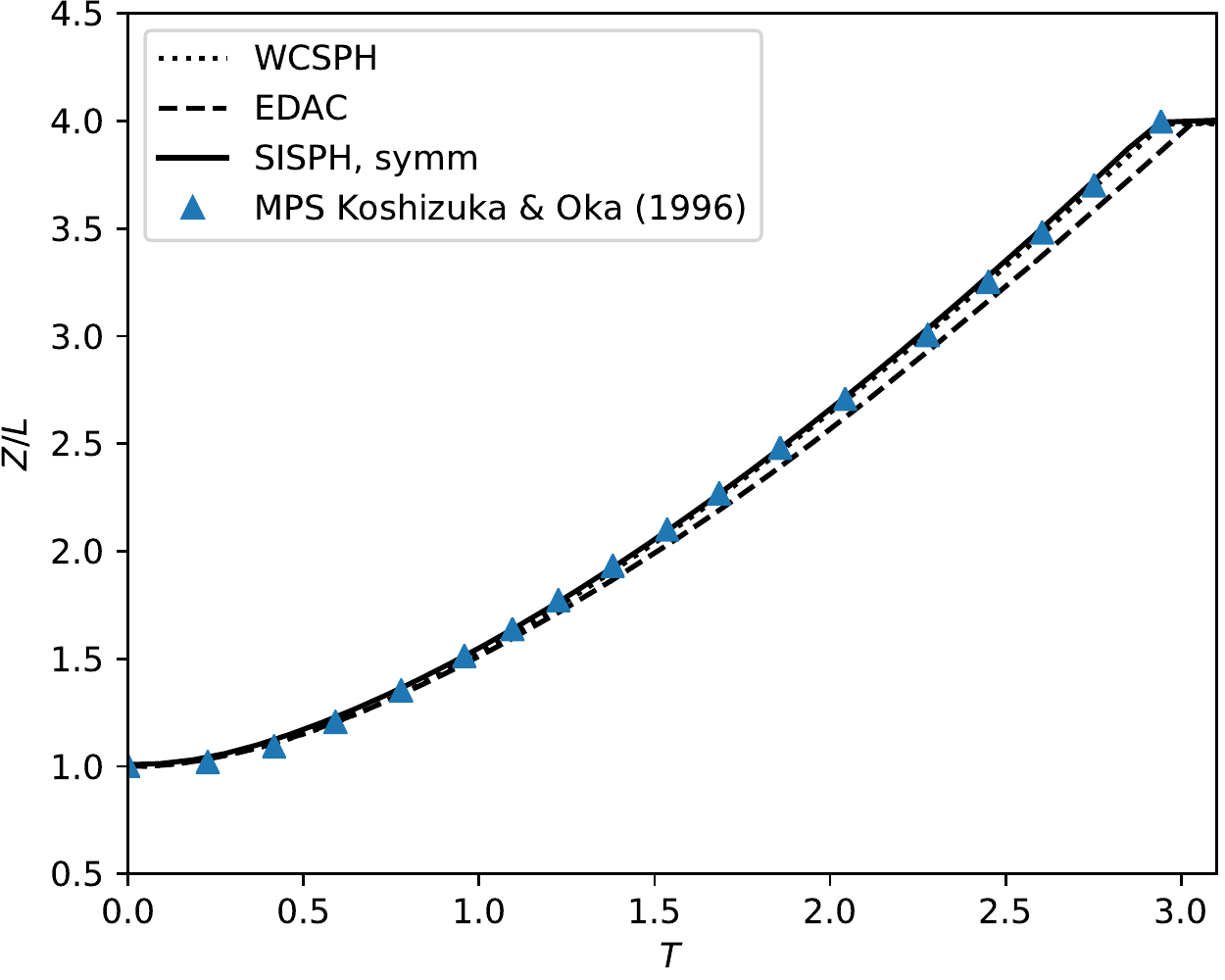}
  \caption{X-coordinate of the toe of the dam versus time as computed with the
    new scheme, WCSPH, EDAC and the simulation
    of~\cite{koshizuka_oka_mps:nse:1996}. $Z$ is the distance of toe of the
    dam from the left wall and $L$ is the initial width of the dam}
\label{fig:db2d:toe_vs_t}
\end{figure}

The particle distribution of the flow at different times is plotted in
Fig.~\ref{fig:db:particle_plots:pressure}. The results show that the scheme
works very well. Fig.~\ref{fig:db2d:toe_vs_t} plots the toe of the breaking
water versus time. The results are compared with the WCSPH scheme, EDAC scheme
and the numerical simulations in~\citet{koshizuka_oka_mps:nse:1996}. As can be
seen, the results are in very good agreement.

\subsection{Dam break in three-dimensions}
\label{sec:db3d}

\begin{figure}[!h]
  \centering
  \begin{subfigure}{0.32\linewidth}
    \includegraphics[width=1.0\linewidth]{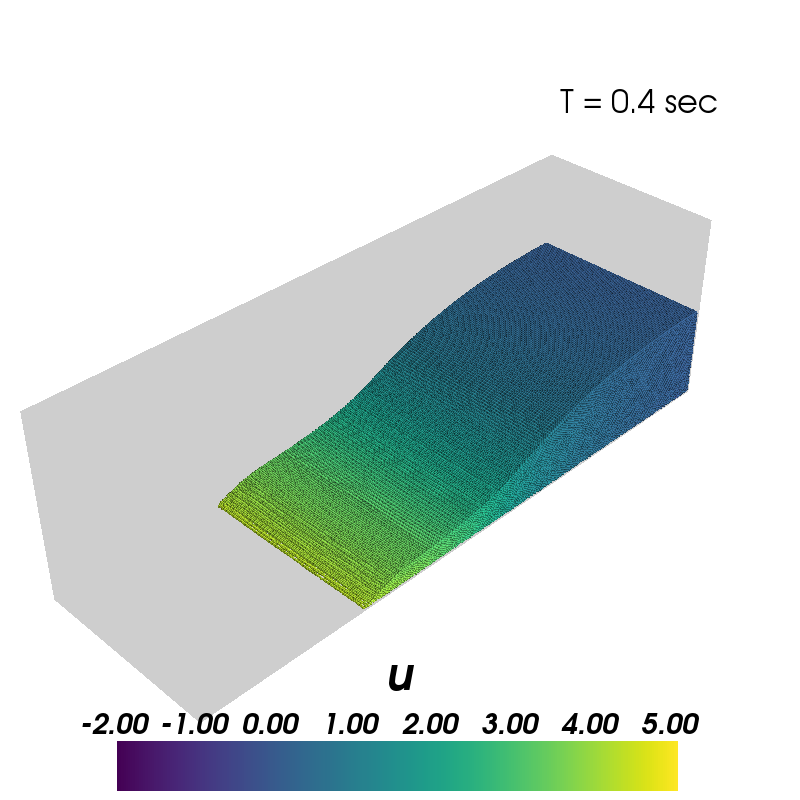}
  \end{subfigure}
  \begin{subfigure}{0.32\linewidth}
    \includegraphics[width=1.0\linewidth]{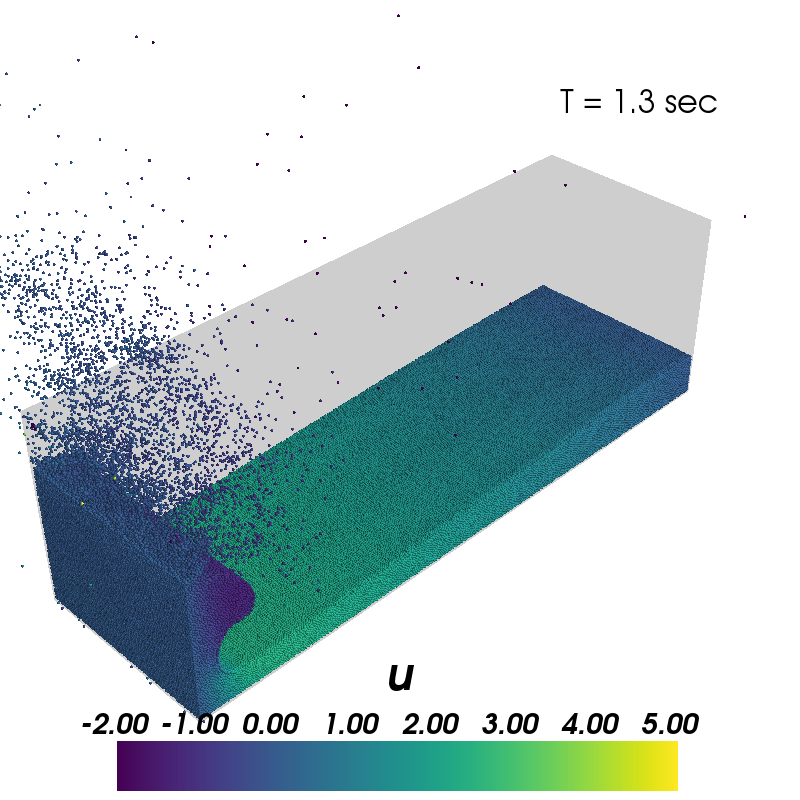}
  \end{subfigure}
  \begin{subfigure}{0.32\linewidth}
    \includegraphics[width=1.0\linewidth]{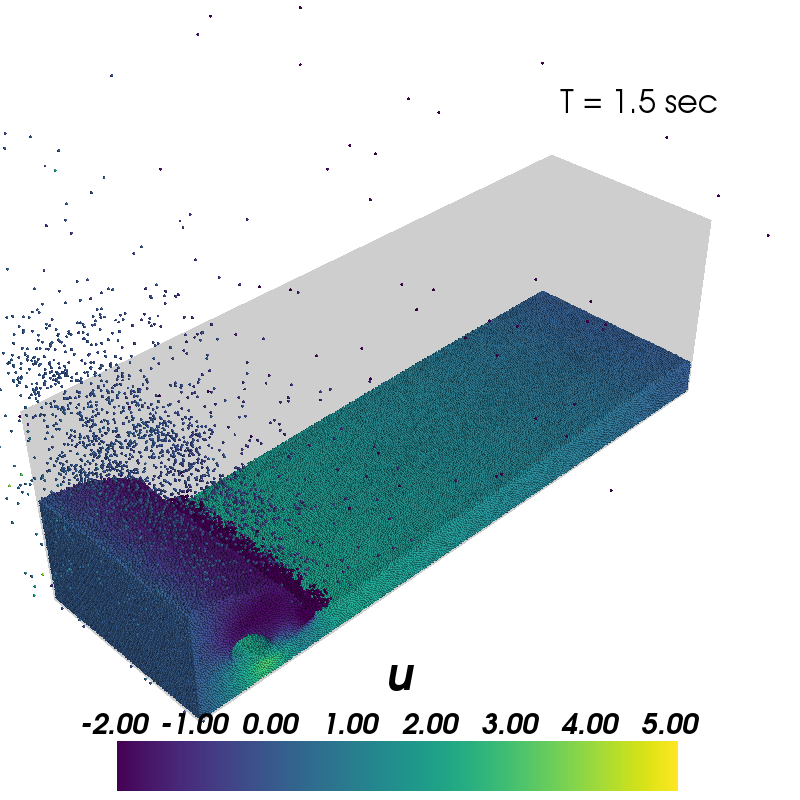}
  \end{subfigure}
\\
  \begin{subfigure}{0.32\linewidth}
    \includegraphics[width=1.0\linewidth]{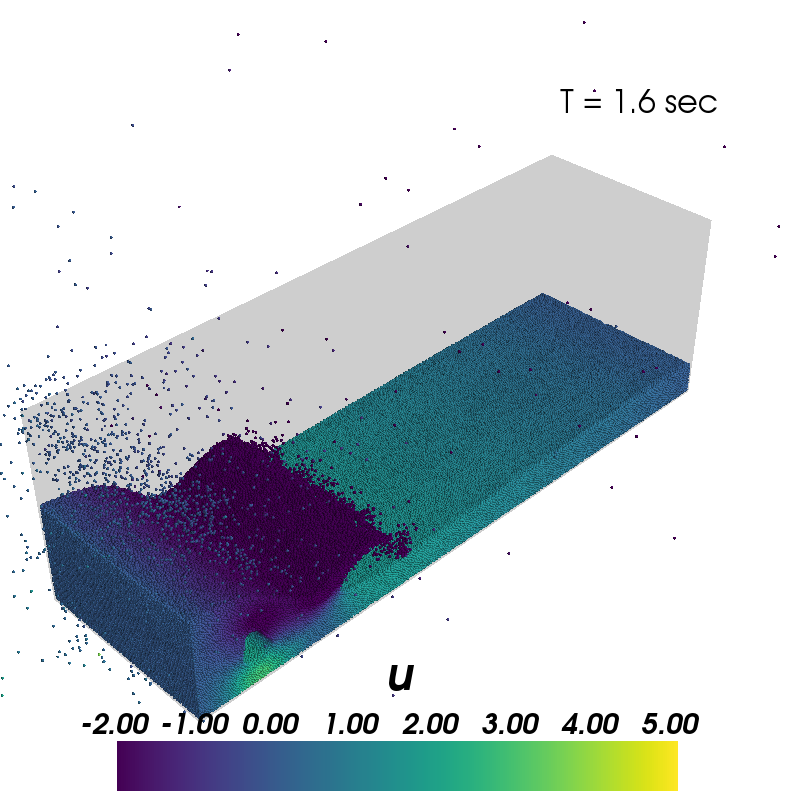}
  \end{subfigure}
  \begin{subfigure}{0.32\linewidth}
    \includegraphics[width=1.0\linewidth]{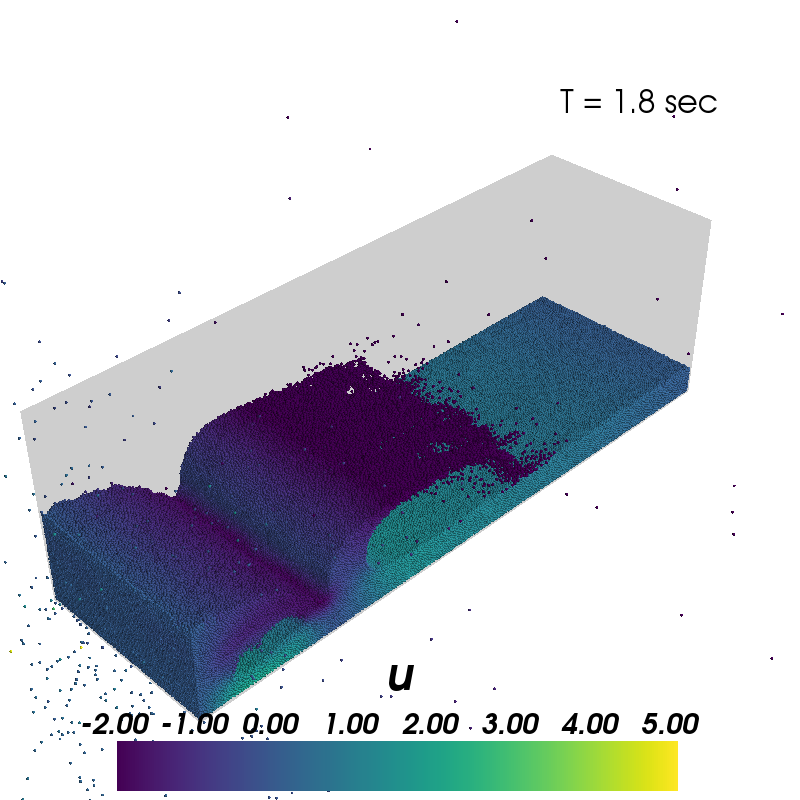}
  \end{subfigure}
  \begin{subfigure}{0.32\linewidth}
    \includegraphics[width=1.0\linewidth]{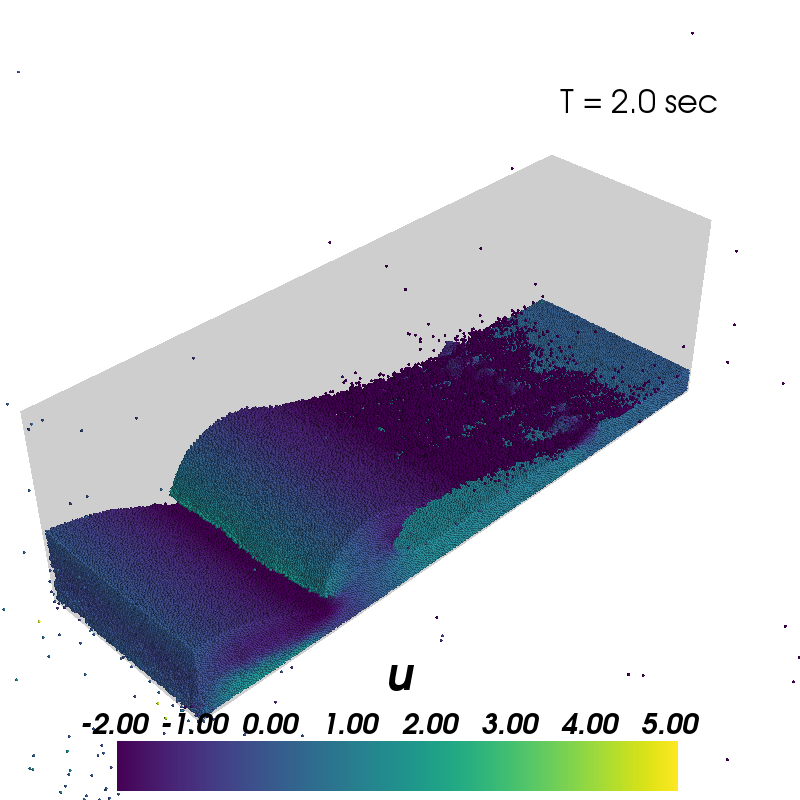}
  \end{subfigure}
  \caption{Dam break in three dimensions using SISPH
shown at times t = $(0.4, 1.3, 1.5, 1.6, 1.8, 2.1)$secs, color indicates $u$
velocity.}
\label{fig:db3d}
\end{figure}

A three-dimensional dam break is considered in order to show that the scheme
works in three dimensions and demonstrates the performance of the scheme on a
GPU. The problem considered is a small modification of the problem described
in~\cite{lee_violeau:db3d:jhr2010}, wherein an obstacle is placed in the path
of the flow. We do not use an obstacle in our example. A rectangular block of
fluid having height $h_w = 0.55m$, width of $1m$ and length of $1.228m$ is at
rest at one end of a container. The container is a long rectangular container
of length $3.22m$ and open at the top. Four layers of boundary particles are
used to enforce the no-penetration boundary condition. A quintic spline kernel
is used with $h/\Delta x = 1.0$. The reference velocity for the flow is
$v_{ref}=\sqrt{2 g h_w}$, where $g=9.81m/s^2$, and the tolerance is,~$\epsilon
= 0.01$. The particle spacing is chosen as $\Delta x = 0.01m$, this results in
around 664k fluid particles and around 627k solid wall particles. The time step
is chosen to be $h/(8 v_{ref})$. An artificial viscosity is used with the
parameter $\alpha=0.25$. The results at different times are shown in
Fig.~\ref{fig:db3d}. There is a large amount of splashing initially but the
scheme produces good results and are similar to those of
\citet{chow:isph:cpc:2018}.

We execute this on a GPU using OpenCL. The code is executed on an NVIDIA 1070
Ti processor with single precision. It is also executed with double precision
on a quad core Intel i5--7400 CPU running at 3.0 Ghz. The user-code is
identical as PySPH~\cite{PR:pysph:scipy16,pysph} automatically executes the
code using OpenCL or OpenMP depending on the options passed. On the GPU, the
execution for 3 seconds of simulation takes 4531 seconds and the same
simulation on the CPU with OpenMP takes 40264 seconds suggesting that the GPU
execution is around 8.9 times faster than the CPU (quad core) for this size of
problem.

The GPU implementation is not the most highly optimized but the above clearly
shows that our new algorithm can be executed effectively on a GPU without
requiring any additional effort as normally needed for ISPH
schemes~\cite{chow:isph:cpc:2018}.
%

\subsubsection{Performance of SISPH}

\begin{figure}[!h]
  \centering
  \includegraphics[width=0.6\linewidth]{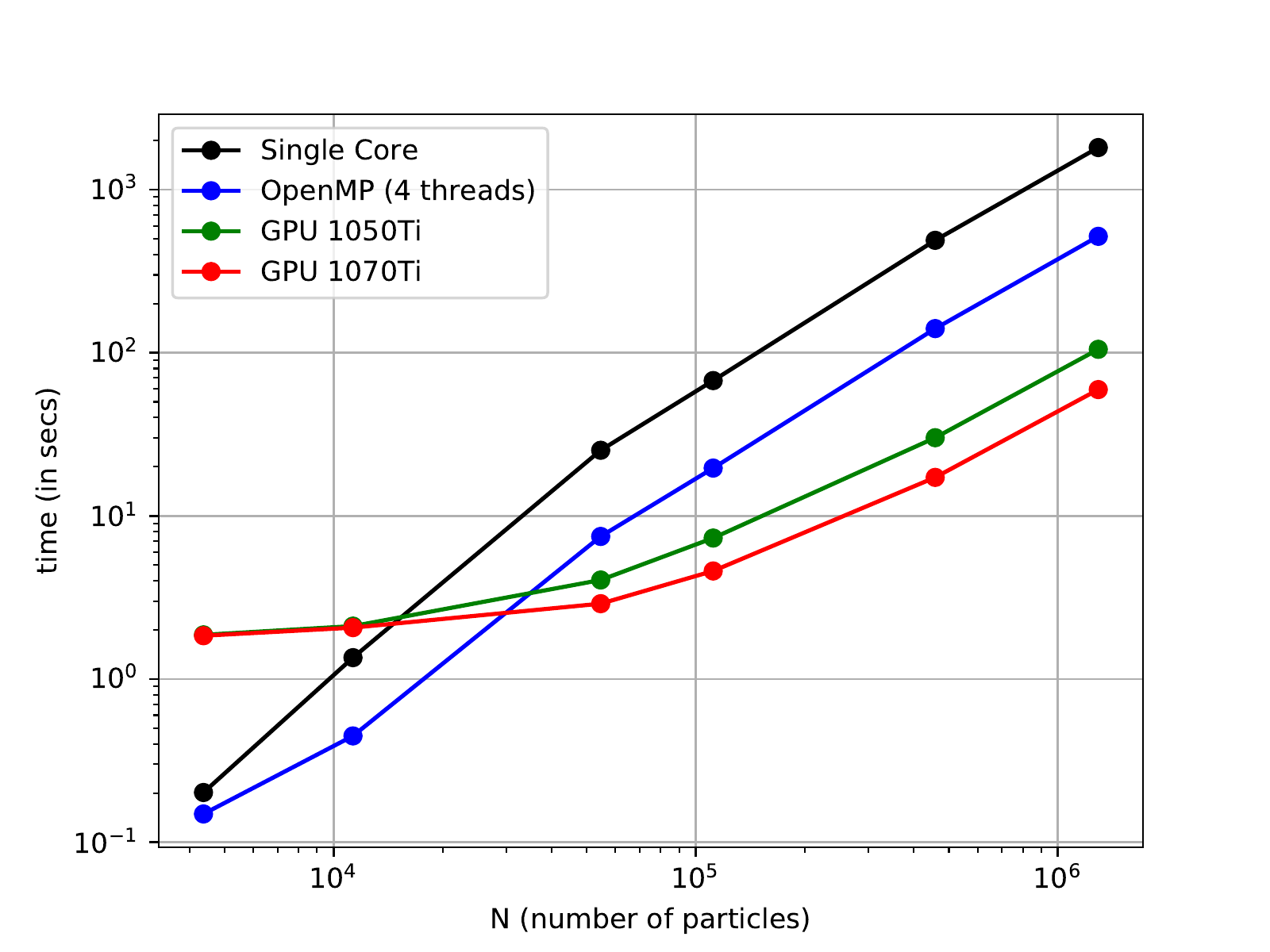}
  \caption{Log-log plot showing the performance of SISPH on different
platforms, no of particles are shown on the x-axis and time taken is shown on
the y-axis.}
\label{fig:db3d:perf}
\end{figure}

\begin{figure}[!h]
  \centering
  \includegraphics[width=0.6\linewidth]{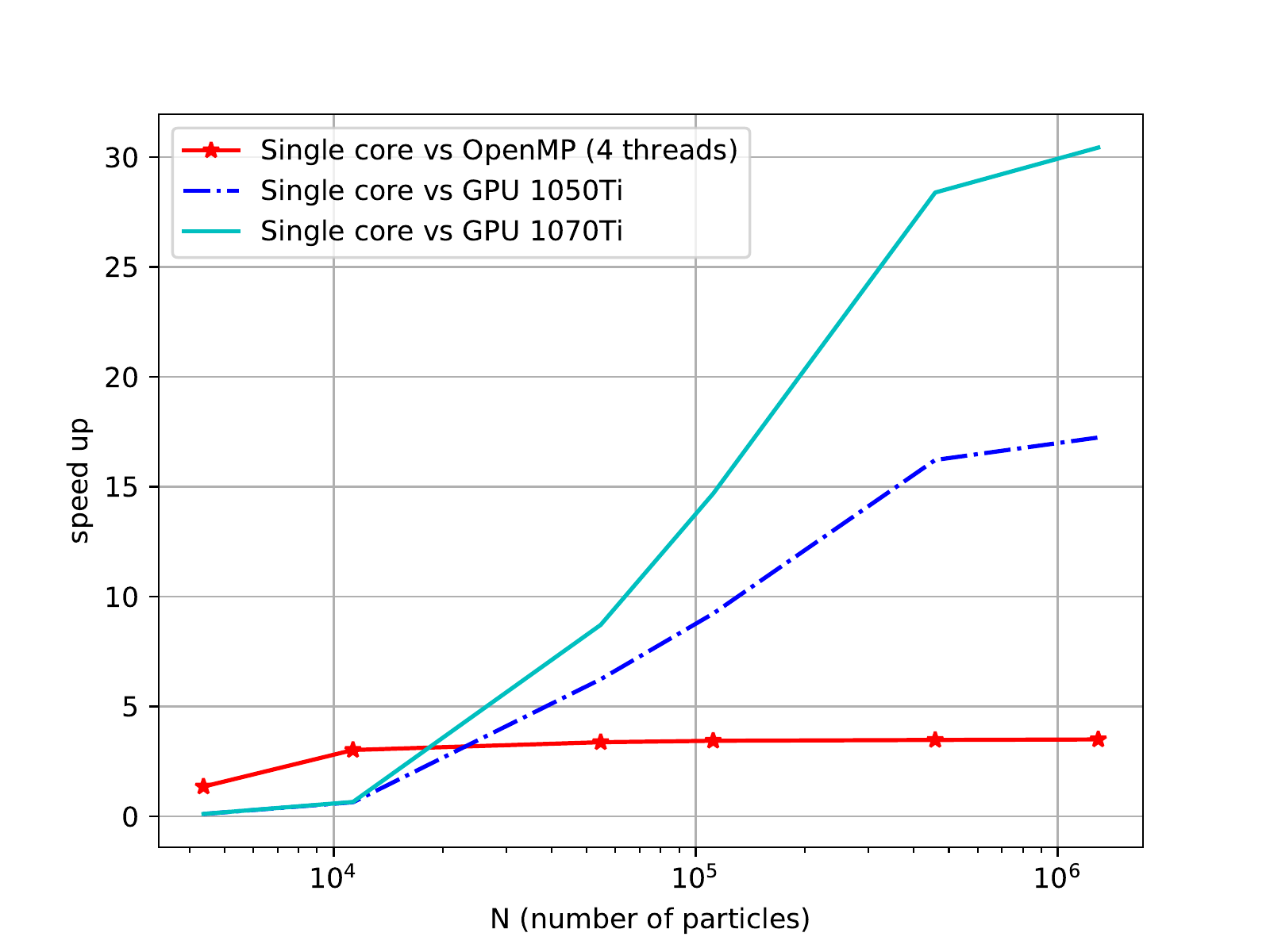}
  \caption{Plot showing the scale up of SISPH on multi-core, NVIDIA 1050Ti, and
      1070Ti GPU vs performance on single core.}
\label{fig:db3d:speedup}
\end{figure}

Here we show the performance of SISPH scheme on single-core, and multi-core
CPU (Intel i5-7400) with double precision, and NVIDIA 1050Ti, 1070Ti GPUs with
single precision. In Fig.~\ref{fig:db3d:perf} the increase in number of
particles, $N$, vs time taken by SISPH on various platforms is show, CPUs
outperform GPUs when $N$ is small as can be expected. For sufficiently large
$N$, at around 100k particles, we can start to see the improved performance of
the GPUs (albeit with single precision), and for 1M particles the time is
reduced by an order of magnitude, also the ``scale up'' is linear.
Fig.~\ref{fig:db3d:speedup} shows speed up of SISPH on multi-core and GPUs
compared to that of single core CPU, while the performance on multi-core (with
4 cores) expectedly peaks around 4 times that of single core, performance of
GPUs increases with number of particles and reaches the peak of around 16
times that of single CPU core for 1050Ti and around 32 times for 1070Ti,
showing the increase in performance with the proposed scheme.


\subsection{Flow past a circular cylinder}
\label{sec:cylinder}

\begin{figure}
\centering
\begin{subfigure}{0.6\linewidth}
\centering
\begin{tikzpicture}
  \newcommand \D {0.4};
  \draw[<->] (0, 1) --  node[below] {$5D$} (5*\D, 1);

  \draw[<->] (5*\D, 1) --  node[below] {$10D$} (5*\D+10*\D, 1);

  \draw[<->] (12*\D, 0) --  node[right] {$7.5D$} (12*\D, 7.5*\D);
  \draw[<->] (12*\D, 0) --  node[right] {$7.5D$} (12*\D, -7.5*\D+0.2);

  \draw [fill=gray] ( 5*\D, 0) circle (\D/2);
  \draw[<->] (6*\D,-\D/2) --  node[right] {$D$} (6*\D,\D/2);
  \draw[gray] (5*\D,-\D/2) -- (6*\D, -\D/2);
  \draw[gray] (5*\D,\D/2) -- (6*\D, \D/2);

  \draw[->] (-1*\D, 0) --  node[below] {\footnotesize $inlet$} (0, 0);
  \draw[->] (15*\D, 0) --  node[below] {\footnotesize $outlet$} (16*\D, 0);

  \draw [fill=gray] ( 0, 7.5*\D) rectangle ++(15*\D, 0.2);
  \draw [fill=gray] ( 0, -7.5*\D) rectangle ++(15*\D, 0.2);
\end{tikzpicture}
\end{subfigure}
\begin{subfigure}{0.35\linewidth}
  \centering
  \includegraphics[width=1.0\linewidth]{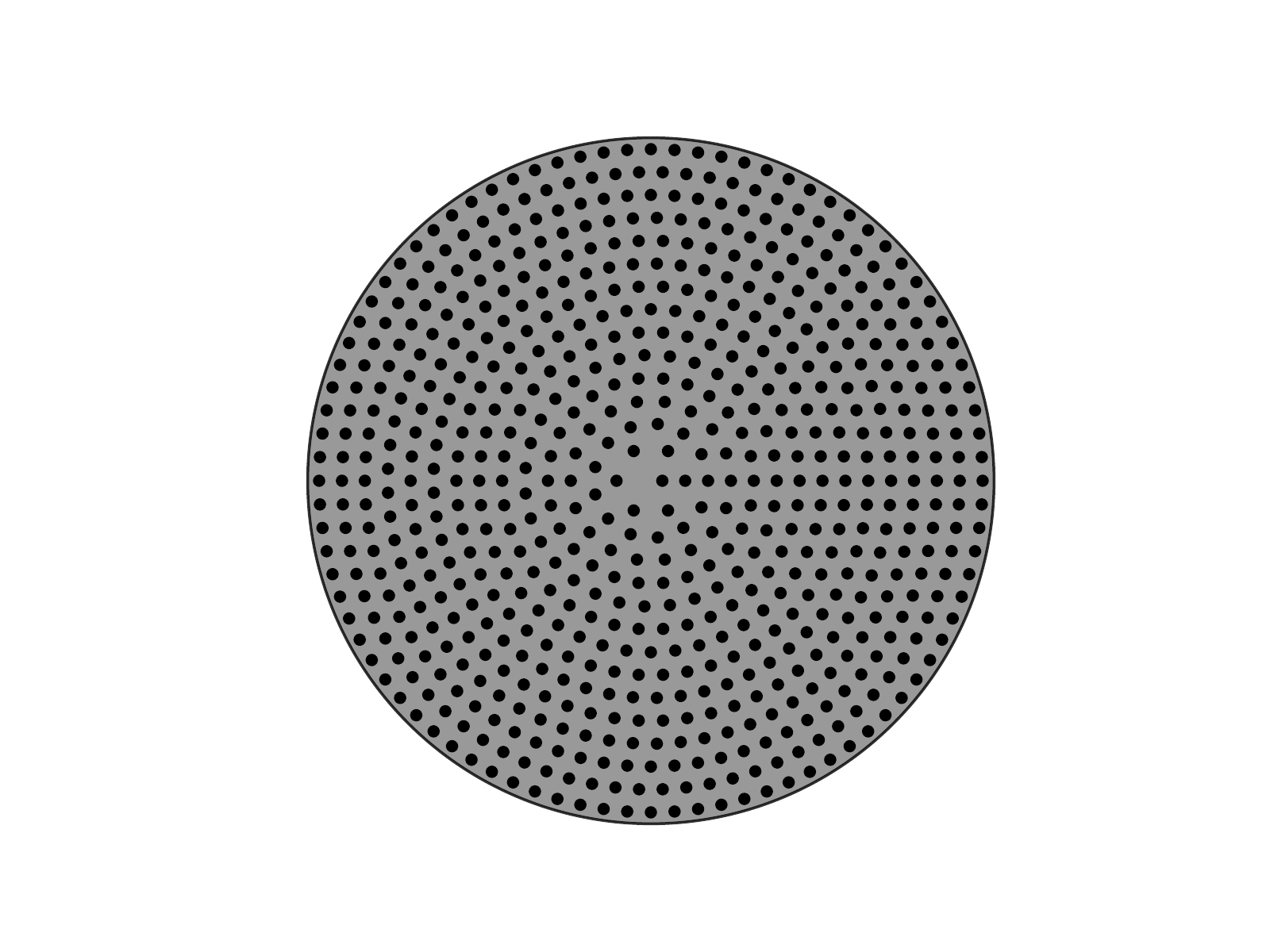}
\end{subfigure}
\caption{The figure on the left shows the setup for flow past cylinder, the
right figure shows the particle distribution on the cylinder.}
\label{fig:fpc:diagram}
\end{figure}

As a final example, the flow past a circular cylinder is considered.
A cylinder of diameter $D=2m$ is considered, the problem setup is as shown in
Fig.~\ref{fig:fpc:diagram}.  Note that the particles discretizing the cylinder
are placed at a distance of $\Delta x/2$ from the outer circumference in order
to effectively simulate the correct location of the boundary of the
cylinder. The cylinder is placed a distance of $5D$ from the inlet and the
outlet is placed $10D$ from the center of the cylinder. The inlet is set to a
uniform velocity, $U$, along the x-axis of 1m/s. The sides of the wind-tunnel
are set to slip walls. These walls are placed at a distance of $7.5D$ from the
center of the cylinder. The viscosity is set to ensure a Reynolds number, $Re
= UD/\nu = 200$. No artificial viscosity is used. The flow is started
impulsively and simulated for $200$ seconds. A quintic spline kernel is used
with $h/\Delta x = 1.2$. Boundary conditions are implemented as discussed in
Sections~\ref{sec:bc} and~\ref{sec:io-bc}. The convergence tolerance is set as
$\epsilon=0.01$. The particle spacing is chosen as $\Delta x = D/30$, which
results in 200k fluid particles. The ``symm'' form is used for the
calculation of pressure gradient. Since the resolution is low, the geometry
curvature is captured by placing particles along the circumference of the
cylinder such that the volume occupied by the particles is approximately
constant as done in~\cite{pawan:inlet-outlet:2019}. For this case we
encountered particle voiding when we use the TVF reference pressure of
$2\max(p_i)$. This was resolved by choosing $p_{\text{ref}} = \rho c^2$, where
$c$ is $10 |\ten{U}|_{\text{max}}$.

\begin{figure}[!h]
  \centering
  \begin{subfigure}{0.48\linewidth}
    \includegraphics[width=1.0\linewidth]{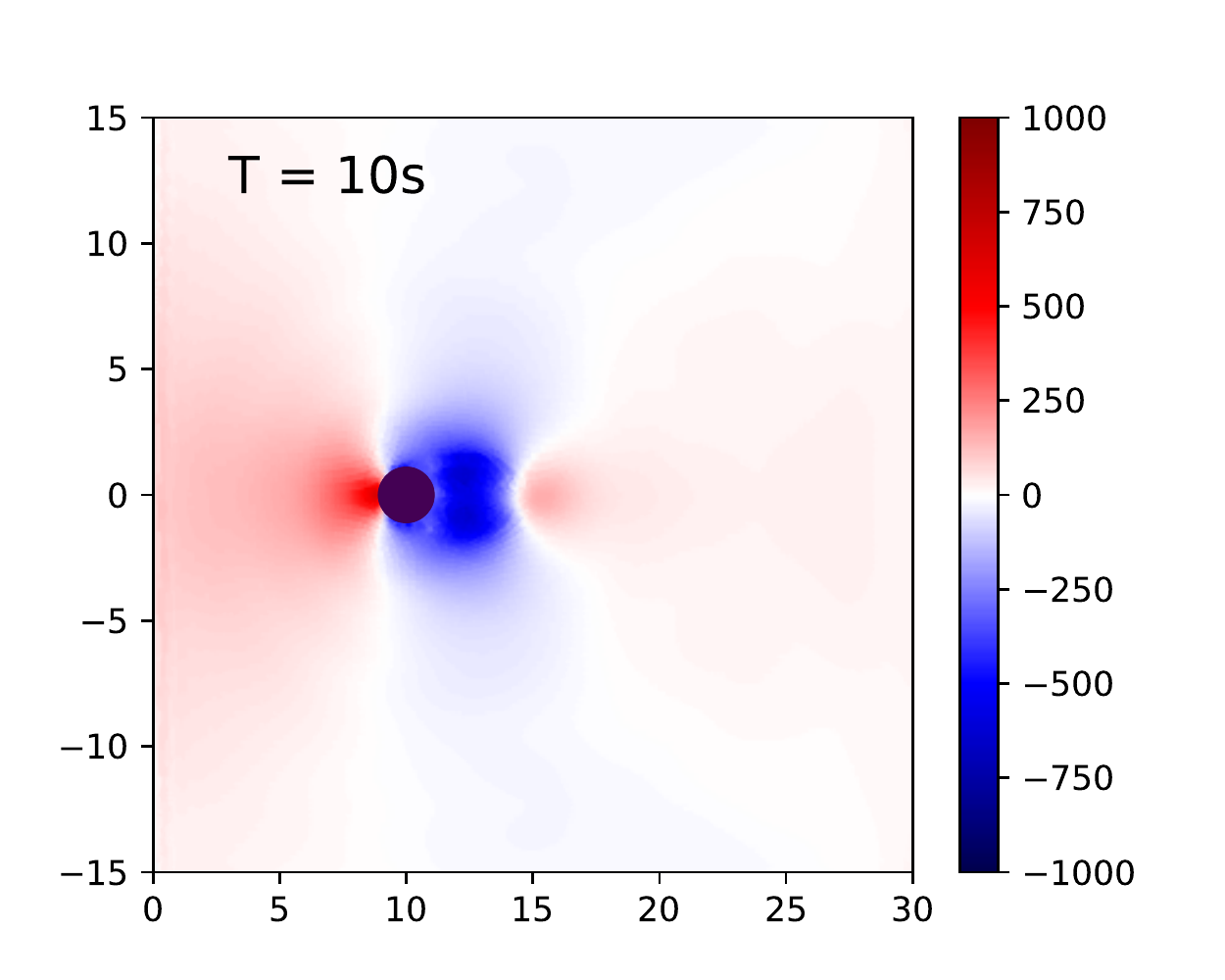}
  \end{subfigure}
  \begin{subfigure}{0.48\linewidth}
    \includegraphics[width=1.0\linewidth]{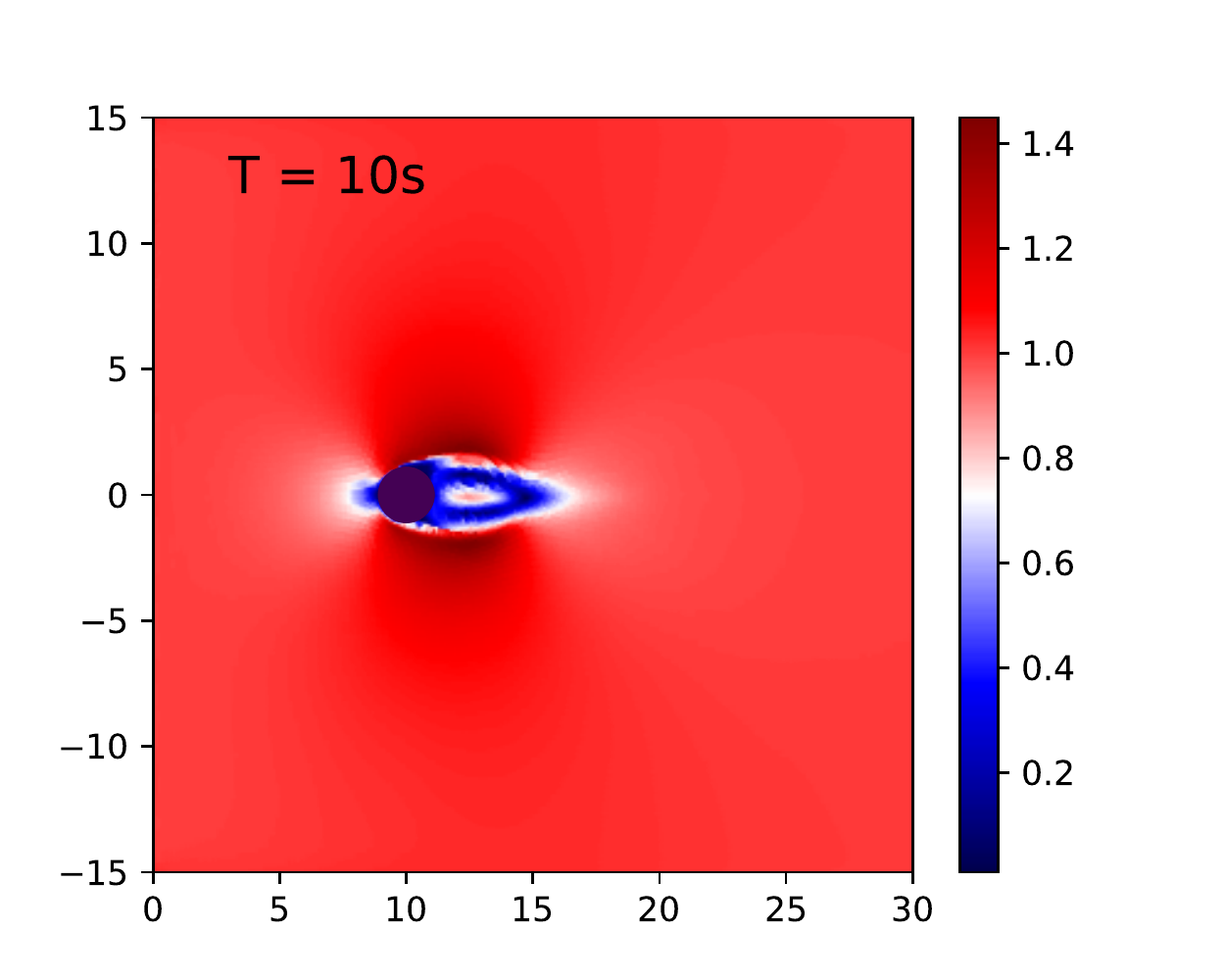}
  \end{subfigure}
\\
  \begin{subfigure}{0.48\linewidth}
    \includegraphics[width=1.0\linewidth]{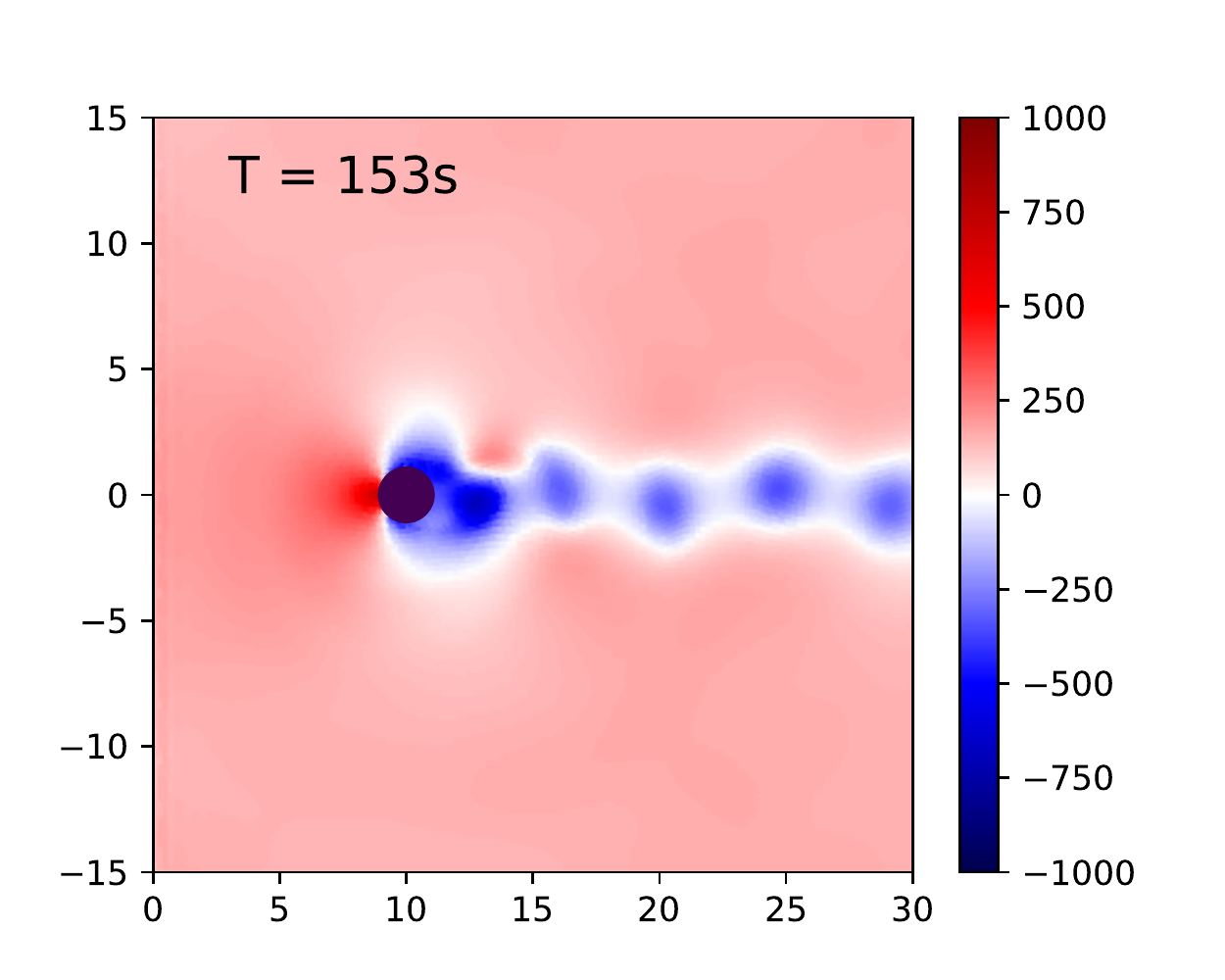}
  \end{subfigure}
  \begin{subfigure}{0.48\linewidth}
    \includegraphics[width=1.0\linewidth]{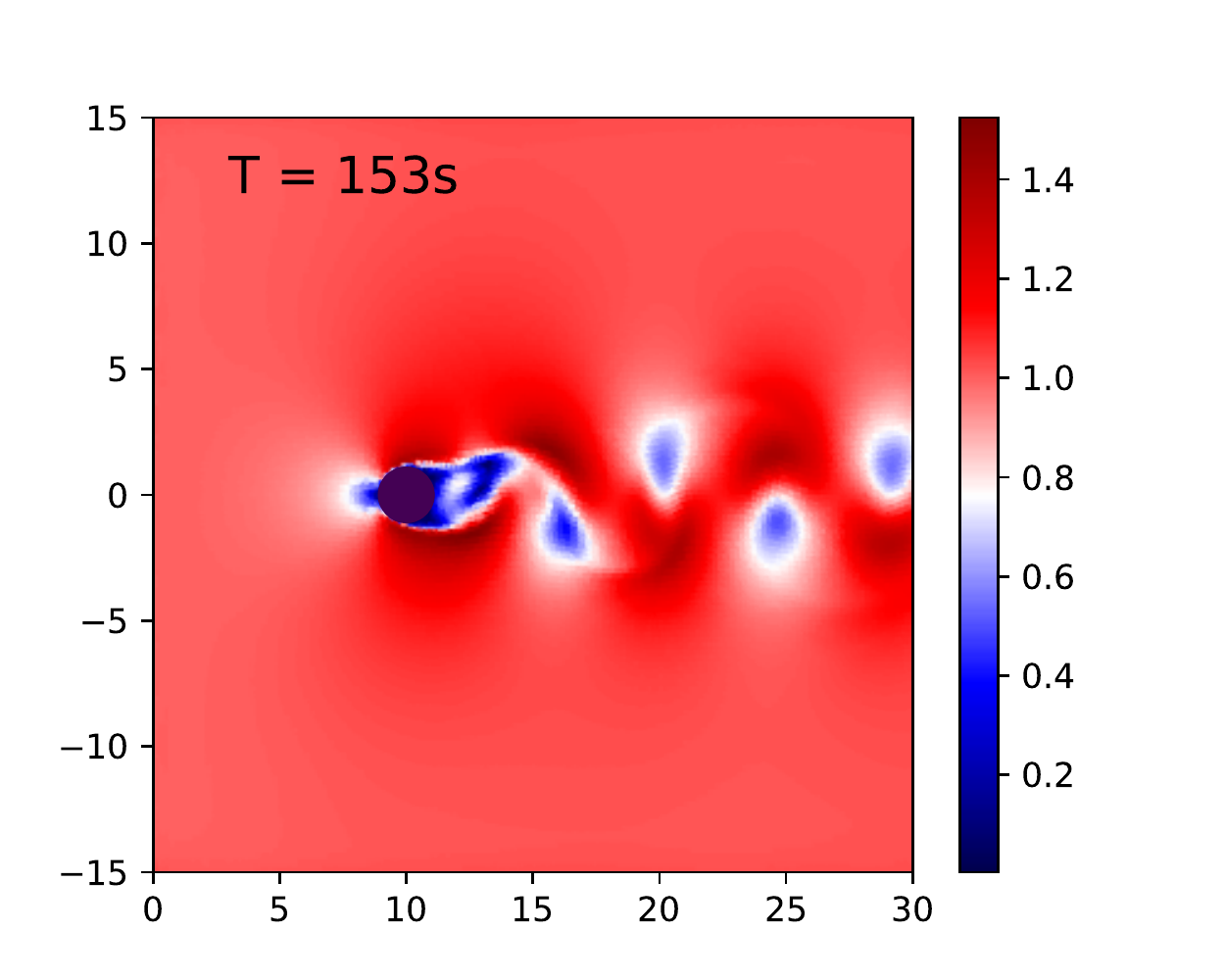}
  \end{subfigure}
\\
  \begin{subfigure}{0.48\linewidth}
    \includegraphics[width=1.0\linewidth]{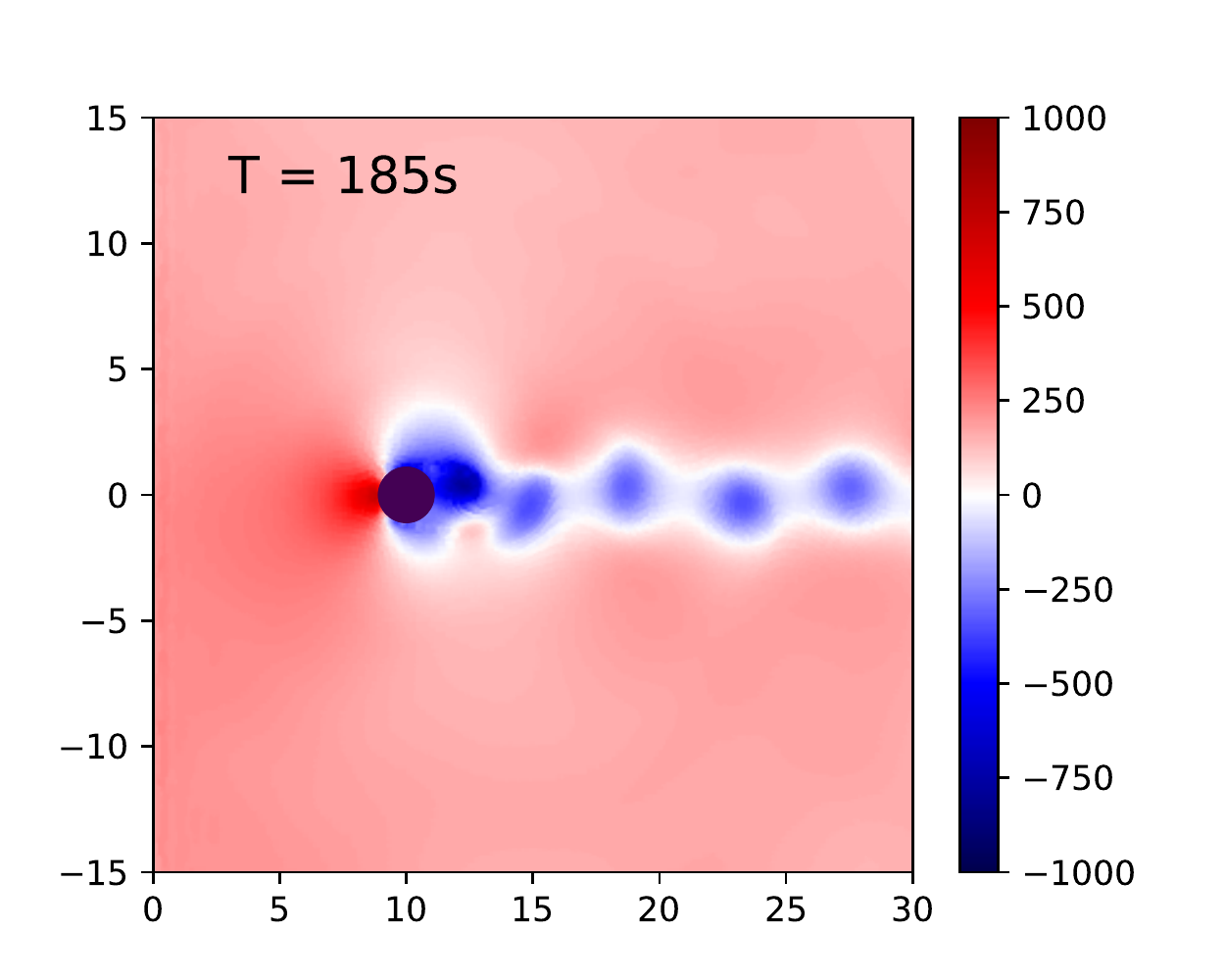}
    \subcaption{Pressure contour}
  \end{subfigure}
  \begin{subfigure}{0.48\linewidth}
    \includegraphics[width=1.0\linewidth]{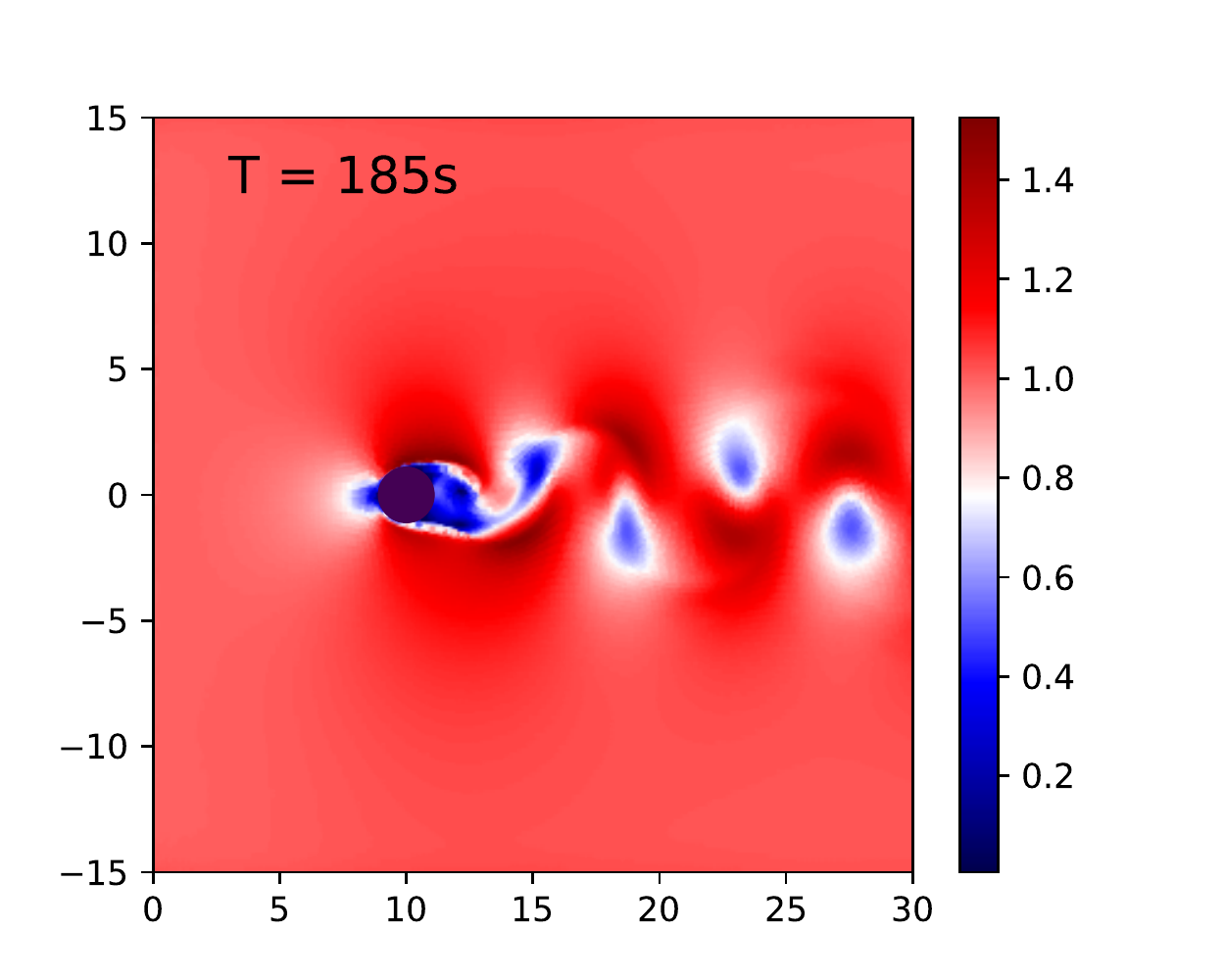}
    \subcaption{Velocity contour}
  \end{subfigure}
  \caption{Plots showing the pressure and velocity contours at times = (10s,
    153s \& 185s), as simulated by SISPH while employing
  ``symm'' form of pressure gradient~\eqref{eq:mom:sph:press:symm}. A $\Delta
  x = D/30$ is used here resulting in 200k fluid particles.}
\label{fig:fpc}
\end{figure}
In Fig.~\ref{fig:db3d}, we plot the velocity and pressure contour by taking
the snapshot at times 10s, 153s and 185s. The pressure and velocity contour
shows an excellent match with the result presented in~\cite{MARRONE2013456,
  open_bc:tafuni:cmame:2018}. The pressure and velocity variations have much
less noise than the WCSPH schemes even with a low resolution of particles as
compared to~\cite{pawan:inlet-outlet:2019}.
\begin{figure}[!h]
  \centering
  \includegraphics[width=1.0\linewidth]{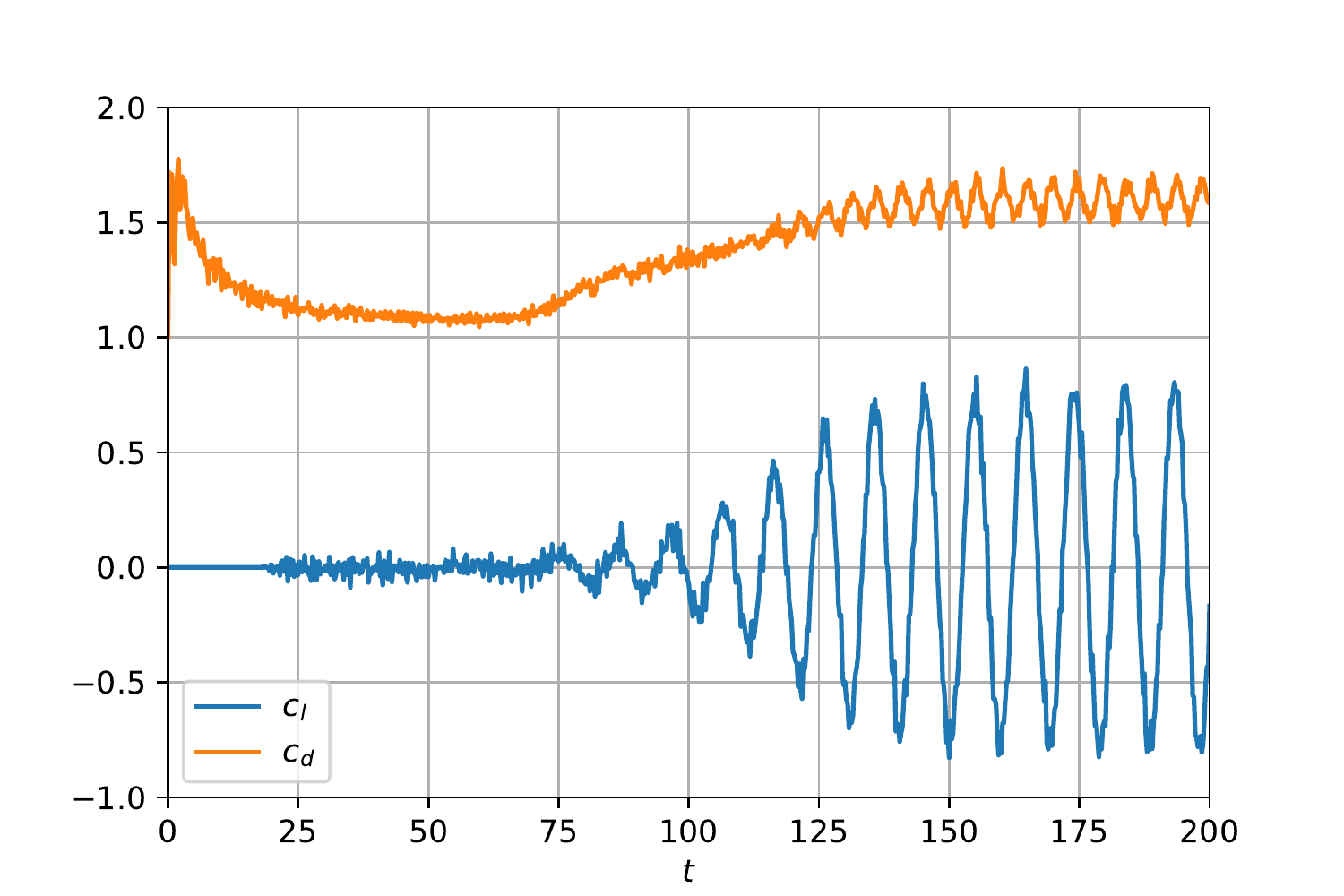}
  \caption{The coefficients of lift, $c_l$, and drag, $c_d$, vs time are
plotted here. After the initial fluctuation, the average coefficient of drag,
$c_d = 1.609$ is observed, the maximum lift coefficient is, $c_l =
0.804$.}
\label{fig:fpc:cl_cd}
\end{figure}

The drag and lift coefficients, $c_d, c_l$ are plotted as a function of time
in Fig.~\ref{fig:fpc:cl_cd}. It must be noted that unlike WCSPH schemes, the
data presented has not gone through any filtering to remove high frequency
oscillations. In order to calculate the forces on the cylinder, we compute
acceleration on the solid due to the pressure and velocity gradients using the
following equation,
\begin{equation}
  \frac{\ten{f}_{\text{solid}}}{m_{\text{solid}}} = - \nabla p + \nu \nabla
  \cdot (\nabla \ten{u})
\end{equation}
where $\ten{f}_{\text{solid}}$ is the force on the cylinder particles,
using SPH discretization~\cite{Adami2013, edac-sph:cf:2019}, the above
equation is written as,
\begin{equation}
  \label{eq:tvf-momentum}
  \begin{split}
    \ten{f}_{i, \text{solid}} = \sum_j \left( V_i^2 +
    V_j^2 \right) & \left[ -\tilde{p}_{ij} \nabla W_{ij}  + \tilde{\eta}_{ij}
    \frac{\ten{u}_{ij}}{(r_{ij}^2 + \eta h_{ij}^2)} \nabla W_{ij}\cdot
    \ten{r}_{ij} \right]
  \end{split}
\end{equation}
where,
\begin{equation}
  \label{eq:tvf-p-ij}
  \tilde{p}_{ij} = \frac{\rho_j p_i + \rho_i p_j}{\rho_i + \rho_j},
\end{equation}
\begin{equation}
  \label{eq:tvf-eta-ij}
  \tilde{\eta}_{ij} = \frac{2 \eta_i \eta_j}{\eta_i + \eta_j},
\end{equation}
\begin{equation}
  \ten{u}_{ij} = \ten{u}_{w\,i} - \ten{u}_j
\end{equation}
where $\ten{u}_{w\,i}$ is the velocity of the solid wall particles as obtained
from equation~\eqref{eq:vg-no-normal}, $\eta_i = \rho_i \nu_i$.

We obtain an average coefficient of drag, $c_d=1.609$ (once the vortex
shedding has been established) and the maximum coefficient of lift, $c_l=
0.804$.  In the Fig.~\ref{fig:fpc:cl_cd}, $c_d$ and $c_l$ show clear
oscillations without any change in the amplitude, which is consistent with the
results of others. With this we have demonstrated the suitability and
robustness of the proposed scheme for a variety of problems.

\section{Conclusions and Future work}
\label{sec:conclusions}

In this paper we introduce a simple, iterative, incompressible SPH scheme that
is based on the original projection formulation of
\citet{sph:psph:cummins-rudman:jcp:1999}. The method is matrix-free, fast, and
suitable for execution on GPUs. The formulation is simple and easy to
implement. We introduce a novel technique to ensure a homogeneous distribution
of particles. In addition we introduce a modified solid wall boundary
condition and show how we can implement simple inlet and outlet boundary
conditions. We demonstrate the accuracy and efficiency of the new scheme with
a suite of benchmark problems in two and three dimensions involving internal
and external flows.

\section*{References}
\bibliographystyle{model6-num-names}
\bibliography{references}

\begin{thebibliography}{39}
\providecommand{\natexlab}[1]{#1}
\providecommand{\url}[1]{\texttt{#1}}
\providecommand{\href}[2]{#2}
\providecommand{\path}[1]{#1}
\providecommand{\DOIprefix}{doi:}
\providecommand{\ArXivprefix}{arXiv:}
\providecommand{\URLprefix}{URL: }
\providecommand{\Pubmedprefix}{pmid:}
\providecommand{\doi}[1]{\href{http://dx.doi.org/#1}{\path{#1}}}
\providecommand{\Pubmed}[1]{\href{pmid:#1}{\path{#1}}}
\providecommand{\BIBand}{and}
\providecommand{\bibinfo}[2]{#2}
\ifx\xfnm\undefined \def\xfnm[#1]{\unskip,\space#1}\fi
\makeatletter\def\@biblabel#1{#1.}\makeatother
\bibitem[{Lucy(1977)}]{lucy77}
\bibinfo{author}{Lucy\xfnm[ L.B.]}.
\newblock \bibinfo{title}{{A numerical approach to testing the fission
  hypothesis}}.
\newblock \emph{\bibinfo{journal}{The Astronomical Journal}}
  \bibinfo{year}{1977};\bibinfo{volume}{82}(\bibinfo{number}{12}):\bibinfo{pages}{1013--1024}.
\bibitem[{Gingold and Monaghan(1977)}]{monaghan-gingold-stars-mnras-77}
\bibinfo{author}{Gingold\xfnm[ R.A.]}, \bibinfo{author}{Monaghan\xfnm[ J.J.]}.
\newblock \bibinfo{title}{Smoothed particle hydrodynamics: Theory and
  application to non-spherical stars}.
\newblock \emph{\bibinfo{journal}{Monthly Notices of the Royal Astronomical
  Society}}
  \bibinfo{year}{1977};\bibinfo{volume}{181}:\bibinfo{pages}{375--389}.
\bibitem[{Monaghan(1994)}]{sph:fsf:monaghan-jcp94}
\bibinfo{author}{Monaghan\xfnm[ J.J.]}.
\newblock \bibinfo{title}{Simulating free surface flows with {SPH}}.
\newblock \emph{\bibinfo{journal}{Journal of Computational Physics}}
  \bibinfo{year}{1994};\bibinfo{volume}{110}:\bibinfo{pages}{399--406}.
\bibitem[{Cummins and Rudman(1999)}]{sph:psph:cummins-rudman:jcp:1999}
\bibinfo{author}{Cummins\xfnm[ S.J.]}, \bibinfo{author}{Rudman\xfnm[ M.]}.
\newblock \bibinfo{title}{An {SPH} projection method}.
\newblock \emph{\bibinfo{journal}{Journal of Computational Physics}}
  \bibinfo{year}{1999};\bibinfo{volume}{152}:\bibinfo{pages}{584--607}.
\bibitem[{Shao and Lo(2003)}]{isph:shao:lo:awr:2003}
\bibinfo{author}{Shao\xfnm[ S.]}, \bibinfo{author}{Lo\xfnm[ E.Y.]}.
\newblock \bibinfo{title}{Incompressible {SPH} method for simulating newtonian
  and non-newtonian flows with a free surface}.
\newblock \emph{\bibinfo{journal}{Advances in Water Resources}}
  \bibinfo{year}{2003};\bibinfo{volume}{26}(\bibinfo{number}{7}):\bibinfo{pages}{787
  -- 800}.
\newblock \URLprefix
  \url{http://www.sciencedirect.com/science/article/pii/S0309170803000307}.
  \DOIprefix\doi{https://doi.org/10.1016/S0309-1708(03)00030-7}.
\bibitem[{Hu and Adams(2007)}]{isph:hu-adams:jcp:2007}
\bibinfo{author}{Hu\xfnm[ X.]}, \bibinfo{author}{Adams\xfnm[ N.]}.
\newblock \bibinfo{title}{{An incompressible multi-phase {SPH} method}}.
\newblock \emph{\bibinfo{journal}{Journal of Computational Physics}}
  \bibinfo{year}{2007};\bibinfo{volume}{227}(\bibinfo{number}{1}):\bibinfo{pages}{264--278}.
\newblock \URLprefix
  \url{http://linkinghub.elsevier.com/retrieve/pii/S0021999107003300}.
  \DOIprefix\doi{10.1016/j.jcp.2007.07.013}.
\bibitem[{Chow et~al.(2018)Chow, Rogers, Lind and Stansby}]{chow:isph:cpc:2018}
\bibinfo{author}{Chow\xfnm[ A.D.]}, \bibinfo{author}{Rogers\xfnm[ B.D.]},
  \bibinfo{author}{Lind\xfnm[ S.J.]}, \bibinfo{author}{Stansby\xfnm[ P.K.]}.
\newblock \bibinfo{title}{Incompressible {SPH} ({ISPH}) with fast poisson
  solver on a {GPU}}.
\newblock \emph{\bibinfo{journal}{Computer Physics Communications}}
  \bibinfo{year}{2018};\bibinfo{volume}{226}:\bibinfo{pages}{81 -- 103}.
\newblock \URLprefix
  \url{http://www.sciencedirect.com/science/article/pii/S0010465518300092}.
  \DOIprefix\doi{10.1016/j.cpc.2018.01.005}.
\bibitem[{Hosseini et~al.(2007)Hosseini, Manzari and Hannani}]{hosseini-2007}
\bibinfo{author}{Hosseini\xfnm[ M.]}, \bibinfo{author}{Manzari\xfnm[ M.]},
  \bibinfo{author}{Hannani\xfnm[ S.]}.
\newblock \bibinfo{title}{A fully explicit three-step {SPH} algorithm for
  simulation of non-newtonian fluid flow}.
\newblock \emph{\bibinfo{journal}{International Journal of Numerical Methods
  for Heat and Fluid Flow}} \bibinfo{year}{2007};\bibinfo{volume}{17}.
\newblock \DOIprefix\doi{10.1108/09615530710777976}.
\bibitem[{Rafiee and Thiagarajan(2009)}]{rafiee-2009}
\bibinfo{author}{Rafiee\xfnm[ A.]}, \bibinfo{author}{Thiagarajan\xfnm[ K.P.]}.
\newblock \bibinfo{title}{An {SPH} projection method for simulating
  fluid-hypoelastic structure interaction}.
\newblock \emph{\bibinfo{journal}{Computer Methods in Applied Mechanics and
  Engineering}}
  \bibinfo{year}{2009};\bibinfo{volume}{198}(\bibinfo{number}{33}):\bibinfo{pages}{2785
  -- 2795}.
\newblock \URLprefix
  \url{http://www.sciencedirect.com/science/article/pii/S0045782509001522}.
  \DOIprefix\doi{https://doi.org/10.1016/j.cma.2009.04.001}.
\bibitem[{Barcarolo(2013)}]{bar13}
\bibinfo{author}{Barcarolo\xfnm[ D.A.]}.
\newblock \bibinfo{title}{{Improvement of the precision and the efficiency of
  the {SPH} method: theoretical and numerical study}}.
\newblock \bibinfo{type}{Theses}; {Ecole Centrale de Nantes (ECN)};
  \bibinfo{year}{2013}.
\newblock \URLprefix \url{https://tel.archives-ouvertes.fr/tel-00904198}.
\bibitem[{Barcarolo et~al.(2014)Barcarolo, Le~Touz{\'e} and
  de~Vuyst}]{barcarolo-2014}
\bibinfo{author}{Barcarolo\xfnm[ D.]}, \bibinfo{author}{Le~Touz{\'e}\xfnm[
  D.]}, \bibinfo{author}{de~Vuyst\xfnm[ F.]}.
\newblock \bibinfo{title}{Validation of a new fully-explicit incompressible
  smoothed particle hydrodynamics method}.
\newblock \emph{\bibinfo{journal}{Blucher Mechanical Engineering Proceedings}}
  \bibinfo{year}{2014};\bibinfo{volume}{1}(\bibinfo{number}{1}).
\bibitem[{Nomeritae et~al.(2016)Nomeritae, Daly, Grimaldi and
  Bui}]{nomeritae-eisph-2016}
\bibinfo{author}{Nomeritae\xfnm[]}, \bibinfo{author}{Daly\xfnm[ E.]},
  \bibinfo{author}{Grimaldi\xfnm[ S.]}, \bibinfo{author}{Bui\xfnm[ H.H.]}.
\newblock \bibinfo{title}{Explicit incompressible {SPH} algorithm for
  free-surface flow modelling: A comparison with weakly compressible schemes}.
\newblock \emph{\bibinfo{journal}{Advances in Water Resources}}
  \bibinfo{year}{2016};\bibinfo{volume}{97}:\bibinfo{pages}{156 -- 167}.
\newblock \URLprefix
  \url{http://www.sciencedirect.com/science/article/pii/S0309170816304456}.
  \DOIprefix\doi{https://doi.org/10.1016/j.advwatres.2016.09.008}.
\bibitem[{{Marrone} et~al.(2011){Marrone}, {Antuono}, {Colagrossi},
  {Colicchio}, {Le Touz{\'e}} and {Graziani}}]{marrone-deltasph:cmame:2011}
\bibinfo{author}{{Marrone}\xfnm[ S.]}, \bibinfo{author}{{Antuono}\xfnm[ M.]},
  \bibinfo{author}{{Colagrossi}\xfnm[ A.]}, \bibinfo{author}{{Colicchio}\xfnm[
  G.]}, \bibinfo{author}{{Le Touz{\'e}}\xfnm[ D.]},
  \bibinfo{author}{{Graziani}\xfnm[ G.]}.
\newblock \bibinfo{title}{{$\delta$-SPH} model for simulating violent impact
  flows}.
\newblock \emph{\bibinfo{journal}{Computer Methods in Applied Mechanics and
  Engineering}}
  \bibinfo{year}{2011};\bibinfo{volume}{200}:\bibinfo{pages}{1526--1542}.
\newblock \DOIprefix\doi{10.1016/j.cma.2010.12.016}.
\bibitem[{Basser et~al.(2017)Basser, Rudman and Daly}]{basser-2017}
\bibinfo{author}{Basser\xfnm[ H.]}, \bibinfo{author}{Rudman\xfnm[ M.]},
  \bibinfo{author}{Daly\xfnm[ E.]}.
\newblock \bibinfo{title}{{SPH} modelling of multi-fluid lock-exchange over and
  within porous media}.
\newblock \emph{\bibinfo{journal}{Advances in Water Resources}}
  \bibinfo{year}{2017};\bibinfo{volume}{108}:\bibinfo{pages}{15 -- 28}.
\newblock \URLprefix
  \url{http://www.sciencedirect.com/science/article/pii/S0309170817304049}.
  \DOIprefix\doi{https://doi.org/10.1016/j.advwatres.2017.07.011}.
\bibitem[{Ihmsen et~al.(2014)Ihmsen, Cornelis, Solenthaler, Horvath and
  Teschner}]{iisph:ihmsen:tvcg-2014}
\bibinfo{author}{Ihmsen\xfnm[ M.]}, \bibinfo{author}{Cornelis\xfnm[ J.]},
  \bibinfo{author}{Solenthaler\xfnm[ B.]}, \bibinfo{author}{Horvath\xfnm[ C.]},
  \bibinfo{author}{Teschner\xfnm[ M.]}.
\newblock \bibinfo{title}{Implicit incompressible {SPH}}.
\newblock \emph{\bibinfo{journal}{{IEEE} Trans Vis Comput Graph}}
  \bibinfo{year}{2014};\bibinfo{volume}{20}(\bibinfo{number}{3}):\bibinfo{pages}{426--435}.
\newblock \URLprefix \url{https://doi.org/10.1109/TVCG.2013.105}.
  \DOIprefix\doi{10.1109/TVCG.2013.105}.
\bibitem[{Xu et~al.(2009)Xu, Stansby and Laurence}]{acc_stab_xu:jcp:2009}
\bibinfo{author}{Xu\xfnm[ R.]}, \bibinfo{author}{Stansby\xfnm[ P.]},
  \bibinfo{author}{Laurence\xfnm[ D.]}.
\newblock \bibinfo{title}{Accuracy and stability in incompressible sph ({ISPH})
  based on the projection method and a new approach}.
\newblock \emph{\bibinfo{journal}{Journal of Computational Physics}}
  \bibinfo{year}{2009};\bibinfo{volume}{228}(\bibinfo{number}{18}):\bibinfo{pages}{6703--6725}.
\newblock \DOIprefix\doi{10.1016/j.jcp.2009.05.032}.
\bibitem[{Lind et~al.(2012)Lind, Xu, Stansby and
  Rogers}]{diff_smoothing_sph:lind:jcp:2009}
\bibinfo{author}{Lind\xfnm[ S.]}, \bibinfo{author}{Xu\xfnm[ R.]},
  \bibinfo{author}{Stansby\xfnm[ P.]}, \bibinfo{author}{Rogers\xfnm[ B.]}.
\newblock \bibinfo{title}{Incompressible smoothed particle hydrodynamics for
  free-surface flows: A generalised diffusion-based algorithm for stability and
  validations for impulsive flows and propagating waves}.
\newblock \emph{\bibinfo{journal}{Journal of Computational Physics}}
  \bibinfo{year}{2012};\bibinfo{volume}{231}(\bibinfo{number}{4}):\bibinfo{pages}{1499
  -- 1523}.
\newblock \DOIprefix\doi{10.1016/j.jcp.2011.10.027}.
\bibitem[{Skillen et~al.(2013)Skillen, Lind, Stansby and
  Rogers}]{fickian_smoothing_sph:skillen:cmame:2013}
\bibinfo{author}{Skillen\xfnm[ A.]}, \bibinfo{author}{Lind\xfnm[ S.]},
  \bibinfo{author}{Stansby\xfnm[ P.K.]}, \bibinfo{author}{Rogers\xfnm[ B.D.]}.
\newblock \bibinfo{title}{Incompressible smoothed particle hydrodynamics
  {(SPH)} with reduced temporal noise and generalised fickian smoothing applied
  to body–water slam and efficient wave–body interaction}.
\newblock \emph{\bibinfo{journal}{Computer Methods in Applied Mechanics and
  Engineering}} \bibinfo{year}{2013};\bibinfo{volume}{265}:\bibinfo{pages}{163
  -- 173}.
\newblock \DOIprefix\doi{10.1016/j.cma.2013.05.017}.
\bibitem[{Adami et~al.(2013)Adami, Hu and Adams}]{Adami2013}
\bibinfo{author}{Adami\xfnm[ S.]}, \bibinfo{author}{Hu\xfnm[ X.]},
  \bibinfo{author}{Adams\xfnm[ N.]}.
\newblock \bibinfo{title}{{A transport-velocity formulation for smoothed
  particle hydrodynamics}}.
\newblock \emph{\bibinfo{journal}{Journal of Computational Physics}}
  \bibinfo{year}{2013};\bibinfo{volume}{241}:\bibinfo{pages}{292--307}.
\newblock \URLprefix
  \url{http://linkinghub.elsevier.com/retrieve/pii/S002199911300096X}.
  \DOIprefix\doi{10.1016/j.jcp.2013.01.043}.
\bibitem[{Zhang et~al.(2017)Zhang, Hu and Adams}]{zhang_hu_adams17}
\bibinfo{author}{Zhang\xfnm[ C.]}, \bibinfo{author}{Hu\xfnm[ X.Y.T.]},
  \bibinfo{author}{Adams\xfnm[ N.A.]}.
\newblock \bibinfo{title}{A generalized transport-velocity formulation for
  smoothed particle hydrodynamics}.
\newblock \emph{\bibinfo{journal}{Journal of Computational Physics}}
  \bibinfo{year}{2017};\bibinfo{volume}{337}:\bibinfo{pages}{216--232}.
\bibitem[{Adami et~al.(2012)Adami, Hu and Adams}]{Adami2012}
\bibinfo{author}{Adami\xfnm[ S.]}, \bibinfo{author}{Hu\xfnm[ X.]},
  \bibinfo{author}{Adams\xfnm[ N.]}.
\newblock \bibinfo{title}{{A generalized wall boundary condition for smoothed
  particle hydrodynamics}}.
\newblock \emph{\bibinfo{journal}{Journal of Computational Physics}}
  \bibinfo{year}{2012};\bibinfo{volume}{231}(\bibinfo{number}{21}):\bibinfo{pages}{7057--7075}.
\newblock \URLprefix
  \url{http://linkinghub.elsevier.com/retrieve/pii/S002199911200229X}.
  \DOIprefix\doi{10.1016/j.jcp.2012.05.005}.
\bibitem[{Ramachandran(2016)}]{PR:pysph:scipy16}
\bibinfo{author}{Ramachandran\xfnm[ P.]}.
\newblock \bibinfo{title}{{PySPH}: a reproducible and high-performance
  framework for smoothed particle hydrodynamics}.
\newblock In: \bibinfo{editor}{Benthall\xfnm[ S.]},
  \bibinfo{editor}{Rostrup\xfnm[ S.]}, eds.
  \emph{\bibinfo{booktitle}{{P}roceedings of the 15th {P}ython in {S}cience
  {C}onference}}. \bibinfo{year}{2016}:\unskip \bibinfo{pages}{127 -- 135}.
\newblock \DOIprefix\doi{10.25080/Majora-629e541a-011}.
\bibitem[{Ramachandran et~al.(10  )Ramachandran, Puri et~al.}]{pysph}
\bibinfo{author}{Ramachandran\xfnm[ P.]}, \bibinfo{author}{Puri\xfnm[ K.]},
  et~al.
\newblock \bibinfo{title}{{PySPH}: a python-based {SPH} framework}.
\newblock \bibinfo{year}{2010--}.
\newblock \URLprefix \url{http://pypi.python.org/pypi/PySPH/}.
\bibitem[{Ramachandran(2018)}]{pr:automan:2018}
\bibinfo{author}{Ramachandran\xfnm[ P.]}.
\newblock \bibinfo{title}{automan: A python-based automation framework for
  numerical computing}.
\newblock \emph{\bibinfo{journal}{Computing in Science \& Engineering}}
  \bibinfo{year}{2018};\bibinfo{volume}{20}(\bibinfo{number}{5}):\bibinfo{pages}{81--97}.
\newblock \URLprefix
  \url{doi.ieeecomputersociety.org/10.1109/MCSE.2018.05329818}.
  \DOIprefix\doi{10.1109/MCSE.2018.05329818}.
\bibitem[{Monaghan(1992)}]{monaghan92}
\bibinfo{author}{Monaghan\xfnm[ J.J.]}.
\newblock \bibinfo{title}{Smoothed particle hydrodynamics}.
\newblock \emph{\bibinfo{journal}{Annual Review of Astronomy and Astrophysics}}
  \bibinfo{year}{1992};\bibinfo{volume}{30}(\bibinfo{number}{1}):\bibinfo{pages}{543--574}.
\newblock \URLprefix \url{https://doi.org/10.1146/annurev.aa.30.090192.002551}.
  \DOIprefix\doi{10.1146/annurev.aa.30.090192.002551}.
  \href{http://arxiv.org/abs/https://doi.org/10.1146/annurev.aa.30.090192.002551}{\tt
  arXiv:https://doi.org/10.1146/annurev.aa.30.090192.002551}.
\bibitem[{Monaghan(2005)}]{monaghan-review:2005}
\bibinfo{author}{Monaghan\xfnm[ J.J.]}.
\newblock \bibinfo{title}{{Smoothed Particle Hydrodynamics}}.
\newblock \emph{\bibinfo{journal}{{Reports on Progress in Physics}}}
  \bibinfo{year}{2005};\bibinfo{volume}{68}:\bibinfo{pages}{1703--1759}.
\bibitem[{Nair and Tomar(2015)}]{nair2015}
\bibinfo{author}{Nair\xfnm[ P.]}, \bibinfo{author}{Tomar\xfnm[ G.]}.
\newblock \bibinfo{title}{Volume conservation issues in incompressible smoothed
  particle hydrodynamics}.
\newblock \emph{\bibinfo{journal}{Journal of Computational Physics}}
  \bibinfo{year}{2015};\bibinfo{volume}{297}:\bibinfo{pages}{689 -- 699}.
\newblock \URLprefix
  \url{http://www.sciencedirect.com/science/article/pii/S0021999115003769}.
  \DOIprefix\doi{https://doi.org/10.1016/j.jcp.2015.05.042}.
\bibitem[{Chorin(1968)}]{chorin1968}
\bibinfo{author}{Chorin\xfnm[ A.J.]}.
\newblock \bibinfo{title}{Numerical solution of the navier-stokes equations}.
\newblock \emph{\bibinfo{journal}{Mathematics of Computation}}
  \bibinfo{year}{1968};\bibinfo{volume}{22}(\bibinfo{number}{104}):\bibinfo{pages}{745--762}.
\newblock \URLprefix \url{http://www.jstor.org/stable/2004575}.
\bibitem[{van~der Vorst(1992)}]{Vor92}
\bibinfo{author}{van~der Vorst\xfnm[ H.A.]}.
\newblock \bibinfo{title}{Bi-cgstab: A fast and smoothly converging variant of
  bi-cg for the solution of nonsymmetric linear systems}.
\newblock \emph{\bibinfo{journal}{SIAM Journal on Scientific and Statistical
  Computing}}
  \bibinfo{year}{1992};\bibinfo{volume}{13}(\bibinfo{number}{2}):\bibinfo{pages}{631--644}.
\newblock \URLprefix \url{https://doi.org/10.1137/0913035}.
  \DOIprefix\doi{10.1137/0913035}.
  \href{http://arxiv.org/abs/https://doi.org/10.1137/0913035}{\tt
  arXiv:https://doi.org/10.1137/0913035}.
\bibitem[{Negi et~al.(2019)Negi, Ramachandran and
  Haftu}]{pawan:inlet-outlet:2019}
\bibinfo{author}{Negi\xfnm[ P.]}, \bibinfo{author}{Ramachandran\xfnm[ P.]},
  \bibinfo{author}{Haftu\xfnm[ A.]}.
\newblock \bibinfo{title}{An improved non-reflecting outlet boundary condition
  for weakly-compressible {SPH}}.
\newblock \emph{\bibinfo{journal}{ArXiV}} \bibinfo{year}{2019};.
\bibitem[{Ramachandran and Puri(2019)}]{edac-sph:cf:2019}
\bibinfo{author}{Ramachandran\xfnm[ P.]}, \bibinfo{author}{Puri\xfnm[ K.]}.
\newblock \bibinfo{title}{Entropically damped artificial compressibility for
  {SPH}}.
\newblock \emph{\bibinfo{journal}{Computers and Fluids}}
  \bibinfo{year}{2019};\bibinfo{volume}{179}(\bibinfo{number}{30}):\bibinfo{pages}{579--594}.
\newblock \DOIprefix\doi{10.1016/j.compfluid.2018.11.023}.
\bibitem[{Jones et~al.(01  )Jones, Oliphant, Peterson et~al.}]{scipy}
\bibinfo{author}{Jones\xfnm[ E.]}, \bibinfo{author}{Oliphant\xfnm[ T.]},
  \bibinfo{author}{Peterson\xfnm[ P.]}, et~al.
\newblock \bibinfo{title}{{SciPy}: Open source scientific tools for {Python}}.
\newblock \bibinfo{year}{2001--}.
\newblock \URLprefix \url{http://www.scipy.org/}.
\bibitem[{Ghia et~al.(1982)Ghia, Ghia and Shin}]{ldc:ghia-1982}
\bibinfo{author}{Ghia\xfnm[ U.]}, \bibinfo{author}{Ghia\xfnm[ K.N.]},
  \bibinfo{author}{Shin\xfnm[ C.T.]}.
\newblock \bibinfo{title}{{High-Re} solutions for incompressible flow using the
  {Navier-Stokes} equations and a multigrid method}.
\newblock \emph{\bibinfo{journal}{Journal of Computational Physics}}
  \bibinfo{year}{1982};\bibinfo{volume}{48}:\bibinfo{pages}{387--411}.
\bibitem[{Colagrossi(2005)}]{colagrossi-phdthesis:2005}
\bibinfo{author}{Colagrossi\xfnm[ A.]}.
\newblock \bibinfo{title}{A meshless lagrangian method for free-surface and
  interface flows with fragmentation}.
\newblock \emph{\bibinfo{journal}{These, Universita di Roma}}
  \bibinfo{year}{2005};\URLprefix \url{http://hdl.handle.net/10805/688}.
\bibitem[{{Ramachandran} et~al.(2019){Ramachandran}, {Muta} and
  {Mokkapati}}]{prabhu:dtsph:2019}
\bibinfo{author}{{Ramachandran}\xfnm[ P.]}, \bibinfo{author}{{Muta}\xfnm[ A.]},
  \bibinfo{author}{{Mokkapati}\xfnm[ R.]}.
\newblock \bibinfo{title}{{Dual-Time Smoothed Particle Hydrodynamics for
  Incompressible Fluid Simulation}}.
\newblock \emph{\bibinfo{journal}{arXiv e-prints}}
  \bibinfo{year}{2019};:\bibinfo{eid}{arXiv:1904.00861}\href{http://arxiv.org/abs/1904.00861}{\tt
  arXiv:1904.00861}.
\bibitem[{Lee et~al.(2010)Lee, Violeau, Issa and
  Ploix}]{lee_violeau:db3d:jhr2010}
\bibinfo{author}{Lee\xfnm[ E.S.]}, \bibinfo{author}{Violeau\xfnm[ D.]},
  \bibinfo{author}{Issa\xfnm[ R.]}, \bibinfo{author}{Ploix\xfnm[ S.]}.
\newblock \bibinfo{title}{Application of weakly compressible and truly
  incompressible {SPH} to 3-d water collapse in waterworks}.
\newblock \emph{\bibinfo{journal}{Journal of Hydraulic Research}}
  \bibinfo{year}{2010};\bibinfo{volume}{48}(\bibinfo{number}{sup1}):\bibinfo{pages}{50--60}.
\newblock \URLprefix \url{https://doi.org/10.1080/00221686.2010.9641245}.
  \DOIprefix\doi{10.1080/00221686.2010.9641245}.
  \href{http://arxiv.org/abs/https://doi.org/10.1080/00221686.2010.9641245}{\tt
  arXiv:https://doi.org/10.1080/00221686.2010.9641245}.
\bibitem[{Koshizuka and Oka(1996)}]{koshizuka_oka_mps:nse:1996}
\bibinfo{author}{Koshizuka\xfnm[ S.]}, \bibinfo{author}{Oka\xfnm[ Y.]}.
\newblock \bibinfo{title}{Moving-particle semi-implicit method for
  fragmentation of incompressible fluid}.
\newblock \emph{\bibinfo{journal}{Nuclear Science and Engineering}}
  \bibinfo{year}{1996};\bibinfo{volume}{123}:\bibinfo{pages}{421--434}.
\bibitem[{Marrone et~al.(2013)Marrone, Colagrossi, Antuono, Colicchio and
  Graziani}]{MARRONE2013456}
\bibinfo{author}{Marrone\xfnm[ S.]}, \bibinfo{author}{Colagrossi\xfnm[ A.]},
  \bibinfo{author}{Antuono\xfnm[ M.]}, \bibinfo{author}{Colicchio\xfnm[ G.]},
  \bibinfo{author}{Graziani\xfnm[ G.]}.
\newblock \bibinfo{title}{An accurate {SPH} modeling of viscous flows around
  bodies at low and moderate reynolds numbers}.
\newblock \emph{\bibinfo{journal}{Journal of Computational Physics}}
  \bibinfo{year}{2013};\bibinfo{volume}{245}:\bibinfo{pages}{456 -- 475}.
\newblock \URLprefix
  \url{http://www.sciencedirect.com/science/article/pii/S0021999113001885}.
  \DOIprefix\doi{https://doi.org/10.1016/j.jcp.2013.03.011}.
\bibitem[{Tafuni et~al.(2018)Tafuni, Dom\'{i}nguez, Vacondio and
  Crespo}]{open_bc:tafuni:cmame:2018}
\bibinfo{author}{Tafuni\xfnm[ A.]}, \bibinfo{author}{Dom\'{i}nguez\xfnm[ J.]},
  \bibinfo{author}{Vacondio\xfnm[ R.]}, \bibinfo{author}{Crespo\xfnm[ A.J.C.]}.
\newblock \bibinfo{title}{A versatile algorithm for the treatment of of open
  boundary conditions in smoothed particle hydrodynamics {GPU} models}.
\newblock \emph{\bibinfo{journal}{Computer methods in applied mechanical
  engineering}}
  \bibinfo{year}{2018};\bibinfo{volume}{342}:\bibinfo{pages}{604--624}.
\newblock \DOIprefix\doi{10.1016/j.cma.2018.08.004}.

\end{thebibliography}

\end{document}